\documentclass[preprint,11pt]{elsarticle}

\usepackage[top=2.5cm,bottom=2.5cm,left=2.0cm,right=2.0cm]{geometry}

\usepackage{amssymb}
\usepackage{amsmath,bm}

\usepackage[bbgreekl]{mathbbol} % for mathbb sigma

\usepackage{tikz}
\usepackage{graphicx}
\usepackage{float}
\usepackage{multirow}
\usepackage{accents}

\newcommand{\bv}{{\bf v}}
\newcommand{\bu}{{\bf u}}
\newcommand{\bg}{{\bf g}}
\newcommand{\be}{{\bf e}}
\newcommand{\bn}{{\bf n}}
\newcommand{\bz}{{\bf z}}
\newcommand{\bx}{{\bf x}}
\newcommand{\by}{{\bf y}}
\newcommand{\bq}{{\bf q}}

\newcommand{\bw}{{\bf w}}
\newcommand{\bbf}{{\bf f}}
\newcommand{\bQ}{{\bf Q}}

\newcommand{\mRT}{{\bf{R_T}}}
\newcommand{\PIT}{{\mathbb{\Pi}_T}}

\newcommand{\BBR}{{\mathbb{R}} }

\def\3bar{{|\hspace{-.02in}|\hspace{-.02in}|}}

\newcommand\Div{\mathrm{div}}

\newtheorem{thm}{Theorem}[section]
\newtheorem{lem}[thm]{Lemma}
\newdefinition{remark}[thm]{Remark}
\newproof{pf}{Proof}
\newdefinition{definition}{Definition}[section]
\newdefinition{algorithm}{Algorithm}[section]
\newdefinition{proposition}{Proposition}[section]

\usepackage{todonotes}

\journal{SISC}%Computers \& Mathematics with Applications}

\begin{document}

%\maketitle

\begin{frontmatter}
\title{A Pressure-Robust Weak Galerkin Finite Element Method for Navier-Stokes Equations}
\author{Lin Mu\fnref{label1}}
\ead{linmu@uga.edu}

\fntext[label1]{Department of Mathematics, University of Georgia, Athens, GA 30602}

\begin{abstract}
In this paper, we develop and analyze a novel numerical scheme for the steady incompressible Navier-Stokes equations by the weak Galerkin methods. The divergence-preserving velocity reconstruction operator is employed in the discretization of momentum equation. By employing the velocity construction operator, our algorithm can achieve pressure-robust, which means, the velocity error is independent of the pressure and the irrotational body force. Error analysis is established to show the optimal rate of convergence. Numerical experiments are presented to validate the theoretical conclusions.
\end{abstract}


\end{frontmatter}

%\begin{AMS}
%Primary, 65N15, 65N30; Secondary, 35B45, 35J50, 35J35
%\end{AMS}
%\pagestyle{myheadings}

%%%%%%%%%%
% Introduction
%%%%%%%%%%
\section{Introduction}

In this paper, we consider the following incompressible Navier-Stokes (NS) equation which seeks velocity $\bu$ and pressure $p$ satisfying
\begin{eqnarray}
-\nu\Delta \bu+(\nabla\times\bu)\times \bu +\nabla p &=& \bbf,\mbox{ in }\Omega\label{eq:pde-1}\\
\nabla\cdot\bu &=& 0,\mbox{ in }\Omega\label{eq:pde-2}\\
\bu &=& 0,\mbox{ on }\partial\Omega\label{eq:pde-bc},
\end{eqnarray}
where $\Omega$ is a polygonal or polyhedral domain in $\mathbb{R}^d$ ($d=2,3$) and $\nu>0$ is the viscosity of the fluid. Here, $\nabla\cdot$, $\nabla\times$, and $\times$ denote, respectively, the divergence operator, curl operator, and the cross product of two vectors. The weak formulation of the rotational NS equations (\ref{eq:pde-1})-(\ref{eq:pde-bc}) seeks $\bu\in[H_0^1(\Omega)]^d$ and $p\in L_0^2(\Omega)$ such that, for all $\bv\in [H_0^1(\Omega)]^d$ and $q\in L^2_0(\Omega)$,
\begin{eqnarray}
\nu(\nabla\bu,\nabla\bv)+%((\nabla\times\bu)\times\bu,\bv)
(\nabla\bu\bu,\bv)-(\nabla\bu\bv,\bu)
-(\nabla\cdot\bv,p) &=& (\bbf,\bv),\label{eq:weakform-1}\\
(\nabla\cdot\bu,q)&=& 0.\label{eq:weakform-2}
\end{eqnarray}
Here $[H_0^1(\Omega)]^d=\{\bu\in[H^1(\Omega)]^d:\bu|_{\partial\Omega}=0\}$ and $L_0^2(\Omega)=\{p\in L^2(\Omega):\int_{\Omega} p = 0\}$ with $[H^1(\Omega)]^d$ being the space of square integrable vector-valued functions whose first derivatives are also square integrable and $L^2(\Omega)$ being the space of square integrable functions. Here for any vector $\by = (y_j)_{1\le j\le d}$, $\nabla\bv\by = \sum_{j=1}^dy_j\partial_j\bv.$

% Stable FE pairs
Due to the importance of the NS equations, the development of efficient and high order NS solver has drawn great attentions. In the finite element (FE) community, in order to carry out a meaningful simulation, the stable velocity and pressure pairs are required to satisfy the inf-sup condition\cite{B1973,B1974}. Many FE pairs have been proposed in previous work and we refer the readers to \cite{BBF2013,GR1979,GR1986} for details. The recent development of discontinuous Galerkin method \cite{PBP2009,CKS2007,KVV2006,LK1999,CS1998} provides new approach for employing discontinuous functions in solving the incompressible problems numerically,  which seems to be more suitable for the convection dominant NS equations. Recently, the hybrid high-order method \cite{DK2018,BDD2019}, virtual element\cite{GMS2018}, hybridizable discontinuous Galerkin method\cite{RW2018,QS2016}, and weak Galerkin finite element \cite{ZL2019,LLC2018} have been proposed to solve NS equations.

% pressure-robust
Besides the stable FE pairs, there is another computational issue demanding extra research efforts. According to the Helmholtz decomposition, arbitrary $L^2$ vector fields can be decomposed into a divergence-free component and an irrotational part. Define the Helmholtz projector $\text{P}$ as the divergence-free part in the vector field, and then the irrotational component of the vector filed will give $\text{P}(\nabla\phi) = 0$. However, this property usually does not hold on the discrete level. %This may lead to numerical issues and sometimes will produce poor mass conservation. 
A lack of this orthogonality may indeed result in poor approximations of the velocity field, whose error estimate has an adverse dependence on the pressure. Restoring the discrete $L^2$-orthogonality between irrotational and discretely divergence-free vector fields is the key for designing mass conserved numerical scheme. The goal of this paper is to develop an optimally convergent discretization method for problem (\ref{eq:pde-1})-(\ref{eq:pde-bc}), which is robust with respect to large irrotational body forces. In the designed scheme, we expect that the velocity error estimate is independent of the pressure, which is also called pressure robust. 

In order to achieve the pressure robust property, the divergence-free mixed finite element on the unstructured tetrahedral grids was presented by Zhang\cite{Z2005}, and later developed into various divergence free schemes \cite{GN2014,GN2014_2,LS2016,Q1994,Z2009,Z2011}. Grad-div stabilization was proposed to overcome the poor mass conservation in the Stokes simulation \cite{OlshanskiiOlshanskii}. Besides, the divergence free hybridizable discontinuous Galerkin\cite{RW2018,QS2016,LS2016,CCQ2017,CCS2005,CS2014,CC2012} and virtually divergence free numerical scheme\cite{CW2019,BLV2018,BLV2017} have been proposed accordingly. As the remark in \cite{ALM2018} that the pressure-robust discretizations need not to be divergence-free.
Such robustness can be achieved by employing the velocity reconstruction operator, which is first proposed by Linke in \cite{Linke2012,JLMNR2017}.  The author proposed to project the velocity into the H(div) space and use this reconstruction only in the body force assembling for Stokes equations for fixing the classical finite element schemes.
Late on, the reconstruction operator has been used in \cite{ALM2018,BrenneckeLinkeMerdon,Linke2014,LinkeMatthiesTobiska,LinkeMerdon2016,PietroErn,QD2020} to solve Stokes equations. The pressure robust scheme for NS equation has been proposed in the convective and rotational formulation for the time-dependent NS equations solved by Taylor-Hood, MINI, and Crouzeix-Raviart finite element. Then, velocity reconstruction operator have been utilize \cite{QD2020} in the hybrid high-order method for providing a irrotational force robust scheme. The robustness is achieved by using high order gradient reconstruction in the rotation formulation, where $\text{P}^{2k+2}(T)^{d\times d}$ is needed to produce a non-dissipative simulation.

% My contribution
Our approach, also inspired by the velocity reconstruction operator, to address this problem is to utilize the divergence preserving velocity construction operator in the stable weak Galerkin finite element formulation. Weak Galerkin finite element method was proposed by Wang and Ye in \cite{wy}. By using discontinuous functions and introducing weak gradient, weak divergence and other weak derivatives, weak Galerkin finite element  methods have been applied to solve different equations including Stokes equations, Biharmonic equations, Maxwells' equations, Navier-Stokes equations \cite{MuYe2017,MuWangWangYe,MuWangYeZhang,HMY2019} and etc. The flexibility in meshing and high order convergence in approximation make the methods attractable in many applications.  {In this manuscript, we shall investigate a novel pressure-robust weak Galerkin numerical scheme, which modifies the classical weak Galerkin schemes with the minimal effort to achieve the pressure independence.} %The recently introduced notion pressure-robust scheme for solving (Navier-)Stokes equations lead to a priori error estimates for the discrete velocity that depend only on velocity $\bu$ but not on pressure $p$ \cite{LinkeMerdon2016}, and thus is called pressure robust numerical scheme. 
Due to the feature of pressure-independent, this scheme is capable of delivering better simulation than other classical methods when the problem with small viscosity values or large irrotational body force.
%Our method is inspired by this divergence preserving velocity construction operator. 
Unlike the prior work for Stokes equation, besides modifying the right hand side body force term, we also need to modify the trilinear form in the discretization in order to achieve the robustness in convergence analysis.

%The contribution of this work consists in extending the construction of? to the fully nonlinear Navier-Stokes problem.  

%The problem considered here is related to recent works pointing out the relevance of restoring at the discrete level the $L^2$-orthogonality between irrotational and discretely divergence-free vector fields. A lack of this orthogonality property may indeed result in poor approximations of the velocity field, whose error estimate has an adverse dependence on the pressure. 

%Based on the fact that $\nabla\cdot\bu = 0$, the weak formulation () can be written as follows,
%\begin{eqnarray}
%\nu(\nabla\bu,\nabla\bv)+\frac{1}{2}((\bu\cdot\nabla\bu,\bv)-(\bu\cdot\nabla\bv,\bu))-(\nabla\cdot\bv,p) &=& (\bbf,\bv),\label{weakform-1}\\
%(\nabla\cdot\bu,q) &=& 0.\label{weakform-2}
%\end{eqnarray}
%Here $\dfrac{1}{2}((\bu\cdot\nabla\bu,\bv)-(\bu\cdot\nabla\bv,\bu))$ is skew symmetric. This is a classical stabilization technique for solving Naiver-Stokes equations.

% Organization
This rest of the paper is organized as follows. In Section~\ref{Sect:FEM}, we introduce some preliminaries and notations for finite element spaces, and then develop the weak Galerkin Algorithm. %A variational formula is presented in Section~\ref{Sect:Scheme}. 
The wellposedness is established in Section~\ref{Sect:Wellpose} and the error estimates are analyzed in Section~\ref{Sect:Error}. Section~\ref{Sect:NumTest} contributes to provide several numerical tests for validating the proposed numerical scheme. Finally, conclusions and future work are discussed in Section~\ref{Sect:Con}.

\section{Finite Element Scheme}\label{Sect:FEM}

\subsection{Finite Element Space}
%In this section, we introduce the WG scheme. 
We adopt standard definitions for the Sobolev spaces $W^{s,r}$ and their associated inner products $(\cdot,\cdot)_{D}$, norms $\|\cdot\|_{W^{s,r}(D)}$, and seminorms $|\cdot|_{W^{s,r}(D)}$ for $s\ge 0$, integer $r$. When $s=0$, we denote $L^r(D) := W^{0,r}(D)$; when $r=2$, we denote $H^s(D):=W^{s,2}(D)$. If $r = 2$, we shall denote $\|\cdot\|_{s,D}$ and $|\cdot|_{s,D}$ as norm and seminorm. Furthermore, if $s=0$, we shall drop the subscript $s$. If moreover, $D=\Omega$ (domain), we shall drop the subscript $D$. 
Let $\mathcal{T}_h$ be a partition of the domain $\Omega$ consisting of triangles/rectangles in two dimensions or tetrahedrons/cubes in three dimensions. Denote by $\mathcal{E}_h$ the set of all edges or flat faces in $\mathcal{T}_h$ and let $\mathcal{E}_h^0=\mathcal{E}_h\backslash\partial\Omega$ be the set of all interior edges or flat faces. For every element $T\in\mathcal{T}_h$, we denote by $h_T$ its diameter and define the mesh size $h=\max_{T\in\mathcal{T}_h}h_T$ for $\mathcal{T}_h$. In addition, throughout this paper, we use $C$ to denote generic constant that is independent of mesh size $h$ and the functions in the estimates. We also employ the notation $a\lesssim b$ for $a\le Cb.$

On the mesh $\mathcal{T}_h$, we define WG finite element space $V_h^0$ for the velocity as follows,
$$V_h^0 = \{\bv=\{\bv_0,\bv_b\}\in V_h,\bv_b = 0\text{ on }\partial\Omega\},$$
where
$$V_h=\{\bv=\{\bv_0,\bv_b\}:\bv_0|_T\in[\text{P}_k(T)]^d,\bv_b|_e\in[\text{P}_k(e)]^d,e\subset\partial T,\forall T\in\mathcal{T}_h\}.$$
Here $\text{P}_k$ denotes the space of polynomials of degree no more than $k$ with $k\ge 0$. We would like to emphasize that there is only single value $\bv_b$ defined on each edge $e\in\mathcal{E}_h$. For the pressure, we define the following finite element space,
$$W_h=\{q\in L_0^2(\Omega),q|_T\in\text{P}_{k}(T),\forall T\in\mathcal{T}_h\}. $$

\subsection{Definition of Weak Derivatives}
We define the weak derivatives that are used to derive the WG scheme. For $\bv\in V_h$ and $T\in\mathcal{T}_h$, we define weak gradient $\nabla_w\bv\in  [\mathbb{RT}_{k}(T)]^{d\times d}$ as the unique polynomial satisfying the following equation
\begin{eqnarray}
(\nabla_w\bv,\tau)_T = -(\bv_0,\nabla\cdot\tau)_T+\langle\bv_b,\tau\cdot\bn\rangle_{\partial T}, \forall \tau\in [\mathbb{RT}_{k}(T)]^{d\times d},
\end{eqnarray}
and define weak divergence $\nabla_w\cdot\bv\in\text{P}_{k}(T)$ as the unique polynomial satisfying
\begin{eqnarray}
(\nabla_w\cdot\bv,q)_T = -(\bv_0,\nabla q)_T+\langle\bv_b,q\bn\rangle_{\partial T},\forall q\in\text{P}_{k}(T),
\end{eqnarray}
where $(\cdot,\cdot)_T = (\cdot,\cdot)_{L^2(T)}$ and $\langle\cdot,\cdot\rangle_{\partial T} = \int_{\partial T}\cdot ds$. Here, the notation $[\mathbb{RT}_k(T)]^{d\times d}$ denotes a tensor with each column as a function belongs vector space $[\mathbb{RT}_k(T)]^d$.

Next, we define the following broken inner product,
\begin{eqnarray*}
(v,w)_{\mathcal{T}_h} &=& \sum_{T\in\mathcal{T}_h}(v,w)_T = \sum_{T\in\mathcal{T}_h}\int_T vw d\bx,\\
\langle v,w\rangle_{\partial\mathcal{T}_h}&=&\sum_{T\in\mathcal{T}_h}\langle v,w\rangle_{\partial T} = \sum_{T\in\mathcal{T}_h}\int_{\partial T}vwds.
\end{eqnarray*}

We furnish finite element space $V_h$ with the discrete $H^1$-like semi-norm such that, for all $\bv\in V_h$
\begin{eqnarray}
\3bar\bv\3bar:=\left(\sum_{T\in\mathcal{T}_h}\|\nabla_w\bv\|_T^2\right)^{1/2}.\label{eq:3bar}
\end{eqnarray}
Moreover, we shall define the equivalence semi-norm
\begin{eqnarray}
\3bar\bv\3bar_1:=\left(\sum_{T\in\mathcal{T}_h}\big(\|\nabla\bv_0\|_T^2+h_T^{-1}\|\bv_0-\bv_b\|_{\partial T}^2\big)\right)^{1/2}.\label{eq:3bar-1}
\end{eqnarray}
The equivalence between (\ref{eq:3bar}) and (\ref{eq:3bar-1}) is shown in following Lemma.
\begin{lem}(\cite{MuWangWangYe} Lemma 3.2) 
For any $\bv_h=\{\bv_0,\bv_b\}\in V_h^0$, we have the equivalence between the following two norms: 
\begin{eqnarray}
\3bar\bv\3bar_{1}
\lesssim \3bar\bv\3bar
\lesssim \3bar\bv\3bar_{1}.\label{eq:3bar-3bar1}
\end{eqnarray}
\end{lem}

\subsection{Divergence-preserving Velocity Reconstruction}
We denote the following two spaces:
\begin{eqnarray*}
{\bf H}(\text{div};\Omega) &=& \{\bw\in{\bf L}^2(\Omega):{\nabla\cdot}\bw\in L^2(\Omega)\},\\
{\bf H}_0(\text{div};\Omega) &=&\{\bw\in{\bf H}(\text{div};\Omega):\bv\cdot\bn|_{\partial\Omega} = 0\},
\end{eqnarray*}
where $\bn$ denotes the outward unit normal of $\partial\Omega$. We define the operator $\mRT: V_h\to [\mathbb{RT}_k(T)]^d\subset {\bf H}_0(\text{div};\Omega)$ such that, for all $\bv=\{\bv_0,\bv_b\}\in V_h$,
\begin{eqnarray}
\int_T\mRT\bv\cdot \bw dT &=& \int_T \bv_0\cdot \bw dT,\quad\forall \bw\in [\text{P}^{k-1}(T)]^d,\label{eq:Def-mRT-1}\\
%\int_{e}\mRT(\bv)\cdot\bn q ds &=&\int_{e} \bv_b\cdot\bn q ds,\quad \forall q\in\text{P}^k(e), e\in\partial T.\label{eq:Def-mRT-2}
\mRT\bv\cdot\bn_e  &=& \bv_b\cdot\bn_e,\quad \forall e\subset\partial T.\label{eq:Def-mRT-2}
\end{eqnarray}
where it is understood that (\ref{eq:Def-mRT-1}) is not needed in the case of $k = 0.$
Classically, the above relation (\ref{eq:Def-mRT-1})-(\ref{eq:Def-mRT-2}) identify $\mRT$ uniquely. The next lemma demonstrates the properties of the re-construction operator $\mRT$.

\begin{lem}\label{Lem:RT} (\cite{M2020} Lemma 3.2)
The operator $\mRT$ is divergence-preserving, i.e., for all $\bv\in V_h$, the following holds:
\begin{eqnarray}
{\nabla\cdot}(\mRT\bv) = \nabla_w\cdot \bv,\label{eq:DivCW-1}
\end{eqnarray}
and $\mRT\bv|_e\cdot\bn$ only depends on $\bv_b|_e\cdot\bn$. 
Besides, for all $\bv\in V_h$, the following bound holds:
\begin{eqnarray}
\|\mRT\bv-\bv_0\|_T\lesssim \sum_{e\in\partial T}h_e^{1/2}\|(\bv_0-\bv_b)\cdot\bn\|_e.\label{eq:Esti-RTv}
\end{eqnarray}
\end{lem}

%\begin{remark}
%In the definition (\ref{eq:weak-grad}), it is not necessary to require $\bbpsi\in [\mathbb{RT}_k(T)]^{2\times 2}\subset [{\bf H}(\text{div};\Omega)]^2$. The basis function for $\bbpsi$ as $\left[ [\text{P}_k(T)]^2+\bx\widehat{\text{P}_k(T)}\right]^2$. Here $\widehat{\text{P}_k(T)}$ denotes a homogeneous polynomial with degree $k$ on $T$.
%\end{remark}

\subsection{Finite element scheme}\label{Sect:Scheme}
%\subsection{Bilinear Forms and Trilinear Form}
Now, we are ready to introduce the following bilinear forms and trilinear form as follows for $\bv,\bw,\bz\in V_h$ and $q\in W_h$
\begin{eqnarray*}
a(\bv,\bw) &=& \nu(\nabla_w\bv,\nabla_w\bw)_{\mathcal{T}_h},\\
b(\bv,q) &=& (\nabla_w\cdot\bv,q)_{\mathcal{T}_h},\\
c(\bv,\bw,\bz) &=& (\nabla_w\bv \mRT\bw,\mRT\bz)_{\mathcal{T}_h} - (\nabla_w\bv \mRT\bz,\mRT\bw)_{\mathcal{T}_h}.
\end{eqnarray*}

Then, the WG discretization of problem (\ref{eq:pde-1})-(\ref{eq:pde-bc}) is summarized blow.
\begin{algorithm}\label{alg:WG}
Find $\bu_h\in V_h^0$ and $p_h\in W_h$ such that
\begin{eqnarray}
a(\bu_h,\bv)+c(\bu_h,\bu_h,\bv)-b(\bv,p_h) &=& (\bbf,\mRT\bv)_{\mathcal{T}_h},\forall \bv\in V_h^0,\label{scheme:wg1}\\
b(\bu_h,q) &=& 0,\forall q\in W_h\label{scheme:wg2}.
\end{eqnarray}
\end{algorithm}

In comparison, we shall also cite the following algorithm. 
\begin{algorithm}\label{alg:WG-2}
(Classical WG scheme) Find $\hat{\bu}_h\in V_h^0$ and $\hat{p}_h\in W_h$ such that
\begin{eqnarray}
a(\hat{\bu}_h,\bv)+\hat{c}(\hat{\bu}_h,\bu_h,\bv)-b(\bv,\hat{p}_h) &=& (\bbf,\bv_0)_{\mathcal{T}_h},\forall \bv\in V_h^0,\label{scheme:wg2-1}\\
b(\hat{\bu}_h,q) &=& 0,\forall q\in W_h\label{scheme:wg2-2}.
\end{eqnarray}
where $\hat{c}(\bw,\bv,\bz) = (\nabla_w\bw\bv_0,\bz_0)_{\mathcal{T}_h}-(\nabla_w\bw\bz_0,\bv_0)_{\mathcal{T}_h}.$
\end{algorithm}

\begin{remark}
The classical WG algorithm\cite{HMY2019} employs the trilinear term 
$$\hat{c}(\bv,\bw,\bz) = \frac{1}{2}(\bv_0\nabla_w\bw,\bz_0)_{\mathcal{T}_h}-\frac{1}{2}(\bv_0\nabla_w\bz,\bw_0)_{\mathcal{T}_h},$$
which is different as the one defined in Algorithm~\ref{alg:WG-2}. Here, in order to compare our proposed new algorithm with the classical WG scheme, we shall modify the scheme in \cite{HMY2019} to Algorithm~\ref{alg:WG-2}. 
For the sake of simplicity, we shall only present the results for Algorithm~\ref{alg:WG} and comment that all the results for existence and uniqueness can be extended to Algorithm~\ref{alg:WG-2} naturally without much difficulties. The convergence results for Algorithm~\ref{alg:WG-2} can be derived follow the similar techniques in \cite{HMY2019}.
\end{remark}

\begin{remark}
It is noted that, Algorithm~\ref{alg:WG} and Algorithm~\ref{alg:WG-2} share the same bilinear form $a(\cdot,\cdot)$ and $b(\cdot,\cdot)$,  but differ at the trilinear form and the body force assembling. Later we shall investigate the advantages of modifying the convective term and body force as Algorithm~\ref{alg:WG}.
\end{remark}
% =================
% Existence and uniqueness
% ====================
\section{Existence and uniqueness of the WG solution}\label{Sect:Wellpose}
In this section, we discuss the well-posedness of the WG scheme (\ref{scheme:wg1})-(\ref{scheme:wg2}). The main theoretical tool we use is the Leray-Schauder fixed point theorem. 
%First, we introduce the following norm for the WG finite element space $V_h^0$,
%\begin{eqnarray}
%\3bar\bv\3bar^2 = \sum_{T\in\mathcal{T}_h}\|\nabla_w\bv\|_T^2+\sum_{T\in\mathcal{T}_h}h^{-1}\|\bv_0-\bv_b\|_{\partial T}^2.
%\end{eqnarray}
%First, we present some lemmas that are used in our analysis. 
First, we shall introduce several operators:
\begin{itemize}
\item Let $\pi_T^\ell/\bm{\pi}_T^\ell$ be the local $L^2$-projection onto $\text{P}_{\ell}(T)/[\text{P}_{\ell}(T)]^d$. Denote $\pi_h$ as the $L^2$-projection with $\pi_h|_T:=\pi_T^{k}$.
%If $\ell = k$, we shall denote $\pi_h:=\pi_h^k$.
\item Let $\bQ_0$ be the $L^2$-projection onto $[\text{P}_k(T)]^d$ and $\bQ_b$ be the $L^2$-projection onto $[\text{P}_k(e)]^d$, respectively. Define the $L^2$-projection $\bQ_h$ of $\bu$ to the WG finite element space $V_h$ by $\bQ_h\bu=\{\bQ_0\bu,\bQ_b\bu\}$.
\item Denote by ${\bf R}_h$ a local Raviart-Thomas-Necelec interpolator. By construction, we have
\begin{eqnarray}
{\bf R}_h\mRT\bQ_h\bu = {\bf R}_h\bu.
\end{eqnarray}
%projection such that ${\bf R}_h\bq\in {\bf H}(\text{div},\Omega)$, and on each $T\in\mathcal{T}_h$, one has ${\bf R}_h\bq\in [\mathbb{RT}_k(T)]^{d}$ 
%and the following equation is satisfied:
%\begin{eqnarray*}
%(\nabla\cdot\bq, v)_T = (\nabla\cdot\BBR_h\bq, v)_T,\ \forall v\in \text{P}_k(T).
%\end{eqnarray*}
\item Denote by $\BBR_h$ a projection such that $\BBR_h\mathbb{q}\in [{\bf H}(\text{div},\Omega)]^d$, and on each $T\in\mathcal{T}_h$, one has $\BBR_h\mathbb{q}\in [\mathbb{RT}_k(T)]^{d\times d}$ and the following equation is satisfied:
\begin{eqnarray}
(\nabla\cdot\mathbb{q},\bv)_T = (\nabla\cdot\BBR_h\mathbb{q},\bv)_T,\ \forall\bv\in[\text{P}_k(T)]^2.\label{eq:Rh-operator}
\end{eqnarray}
\item Denote $\PIT$ as the local $L^2$-projection to $[\mathbb{RT}_k(T)]^{d\times d}$.
%\begin{eqnarray}
%\nabla_w\bQ_h\bu= \PIT(\nabla\bu).
%\end{eqnarray}
\end{itemize}
%Let $\pi_h^\ell$ be the local $L^2$-projection onto $\text{P}_{\ell}(T)$. 
%If $\ell = k$, we shall denote $\pi_h:=\pi_h^k$. and $\mathbb{Q}_h^\ell$ be two local $L^2$-projections onto $\text{P}_{\ell}(T)$ and $[\text{P}_{\ell}(T)]^{d\times d}$, respectively. 
%The $L^2$-projection $\bQ_h$ of $\bu$ to the WG finite element space $V_h$ is defined by $\bQ_h\bu=\{\bQ_0\bu,\bQ_b\bu\}$ with $\bQ_0$ being the $L^2$-projection onto $[\text{P}_k(T)]^d$ and $\bQ_b$ being the $L^2$-projection onto $[\text{P}_k(e)]^d$. Denote by $\BBR_h$ a projection such that $\BBR_h\bq\in [{\bf H}(\text{div},\Omega)]^d$, and on each $T\in\mathcal{T}_h$, one has $\BBR_h\bq\in [\mathbb{RT}_k(T)]^{d\times d}$ and the following equation is satisfied:
%\begin{eqnarray*}
%(\nabla\cdot\bq,\bv_0)_T = (\nabla\cdot\BBR_h\bq,\bv_0)_T,\ \forall\bv_0\in[\text{P}_k(T)]^2.
%\end{eqnarray*}
%For any $\mathbb{\tau}\in[{\bf H}(\text{div},\Omega)]^d$, we have
%\begin{eqnarray}
%\sum_{T\in\mathcal{T}_h}(-\nabla\cdot\tau,\bv_0)_T = \sum_{T\in\mathcal{T}_h}(\BBR_h\tau,\nabla_w\bv)_T,\ \forall\bv = \{\bv_0,\bv_b\}\in V_h.\label{eq:Rh}
%\end{eqnarray}
%Besides, we also define ${\bf r}_h$ as a $L^2$ projection from $[L^2(T)]^{d\times d}$ to $[\mathbb{RT}_k(T)]^{d\times d}$. 

%Denote the $L^2$-projection in the finite element space $V_h$ is given by $\bQ_h\bv=\{\bQ_0\bv,\bQ_b\bv\}$ for $\bv\in [H^1(\Omega)]^2$. 
 The following three lemmas have been proved in \cite{MuYe2017,WY2016}.

\begin{lem}\label{lem:bilinearA}
For any $\bv,\bw\in V_h$, we have 
\begin{eqnarray*}
\left|a(\bv,\bw)\right|&\le& \nu\3bar\bv\3bar\3bar\bw\3bar,\\
a(\bv,\bv) &=& \nu\3bar\bv\3bar^2.
\end{eqnarray*}
\end{lem}

\begin{lem}\label{Lem:GradDiv}
The projection operators $\bQ_h$, $\PIT$, and $\pi_h$ satisfy the following commutative properties: %\cite{WangYeStokes}
\begin{eqnarray*}
\nabla_w(\bQ_h\bv) &=& \PIT(\nabla\bv),\quad \forall \bv\in [H^1(\Omega)]^d,\label{eq:grad-commute}\\
\nabla_w\cdot(\bQ_h\bv) &=& \mathbb{\pi}_h(\nabla\cdot\bv),\quad \forall\bv\in {\bf H}(\Div;\Omega).
\end{eqnarray*}
\end{lem}

%\begin{lem}
%The projection operators $Q_h,\bQ_h$, and $\mathbb{Q}_h$ satisfy the following commutative properties
%\begin{eqnarray}
%\nabla_w(\bQ_h\bv) &=& \Pi_h(\nabla\bv), \ \forall\bv\in [H^1(\Omega)]^d,\\
%\nabla_w\cdot(\bQ_h\bv) &=& \mathbb{Q}_h(\nabla\cdot\bv),\ \forall\bv\in H(\text{div},\Omega),
%\end{eqnarray}
%where $H(\text{div},\Omega)$ is the space of square integrable vector-valued functions whose divergence is also square integrable.
%\end{lem}

%\begin{lem}
%For all $\bv\in V_h$, it holds
%\begin{eqnarray}
%\|\nabla_w\bv\|\le C\|\bv\|
%\end{eqnarray}
%\end{lem}

\begin{lem}\label{lem:InfSup}
There exists a positive constant $\beta$ independent of $h$ such that  %\cite{MuYe2017}
\begin{eqnarray*}
\sup_{\bv\in V_h^0}\frac{b(\bv,\rho)}{\3bar\bv\3bar}\ge \beta\|\rho\|,\ \forall\rho\in W_h.
\end{eqnarray*}
\end{lem}

%The following can be approved for $\bv=\{\bv_0,\bv_b\}\in V_h$,
%\begin{eqnarray*}
%\sum_{T\in\mathcal{T}_h}\|\nabla\bv_0\|_T^2\le C\3bar\bv\3bar^2.
%\end{eqnarray*}
%Define another norm $\3bar\bv\3bar_1$ on $V_h$ for $\bv=\{\bv_0,\bv_b\}$ as follows,
%\begin{eqnarray}
%\3bar\bv\3bar_1^2=\sum_{T\in\mathcal{T}_h}\|\nabla\bv_0\|_T^2+\sum_{e\in\mathcal{E}_h}h_e^{-1}\|\ljump\bv_0\rjump\|_e^2.
%\end{eqnarray}
%Then, it follows that,
%\begin{eqnarray}
%\3bar\bv\3bar_1\le C\3bar\bv\3bar,\forall\bv\in V_h.
%\end{eqnarray}

%The following lemma has been proved in \cite{WWY2008}, which gives a upper bound of the trilinear form $C_\text{skew}(\cdot,\cdot,\cdot)$ by the $\3bar\cdot\3bar_1$ norm.

\begin{lem}
It holds for all $r\in[1,6]$ and all $\bv\in V_h^0$,
\begin{eqnarray}
\|\bv_0\|_{L^r(\Omega)^d}&\lesssim& \3bar\bv\3bar.\label{eq:Embed}\\
\|\mRT(\bv)\|_{L^r(\Omega)^d}&\lesssim& \3bar\bv\3bar.\label{eq:reconstructionEmbed}
\end{eqnarray}
where the hidden constant is independent of both $h$ and $\bv.$
\end{lem}
\begin{pf}
By the discrete Sobolev embeddings in $V_h^0$, it holds (\cite{DD2017}, Proposition 5.4)
\begin{eqnarray*}
\|\bv_0\|_{L^4(\Omega)^d}\le \|\nabla\bv_0\|\lesssim \3bar\bv\3bar_1,
\end{eqnarray*}
where the hidden constant is independent of both $h$ and $\bv$, but possibly depends on $\Omega$, $k$, $r$, and the mesh regularity parameter. Then (\ref{eq:Embed}) follows (\ref{eq:3bar-3bar1}) and the above inequality. 

By the estimate (\cite{QD2020}, Proposition 3) and the equivalence (\ref{eq:3bar-3bar1}),
\begin{eqnarray*}
\|\mRT(\bv)\|_{L^r(\Omega)^d}\lesssim\3bar\bv\3bar_1\lesssim \3bar\bv\3bar,
\end{eqnarray*}
and thus completes the proof.
\end{pf}

%\begin{lem}
%For all $\bv\in H^1(\Omega)^d$, we have the following boundedness property of the $L^2$-projection operator $\bQ_0$:
%\begin{eqnarray}
%\|\nabla(\bQ_0\bv-\bv)\|\le C_I\|\nabla\bv\|.
%\end{eqnarray}
%\end{lem}

%\begin{lem} 
%Let $\sigma,\bv,\bw\in V_h$ and $\3bar\cdot\3bar_1$ be defined in (?). Then we have
%\begin{eqnarray}
%\|\bw\|_{L^4(\Omega)}&\le& C\3bar\bw\3bar_1.%,\\
%%C_\text{skew}(\sigma,\bv,\bw)&\le& C\3bar\sigma\3bar_1\3bar\bv\3bar_1\3bar\bw\3bar_1.
%\end{eqnarray} 
%\end{lem}

Next, we are ready to prove the properties for the trilinear term $c(\cdot,\cdot,\cdot)$.
\begin{lem}
For $\bv=\{\bv_0,\bv_b\}$, $\bw=\{\bw_0,\bw_b\}$, and $\bz=\{\bz_0,\bz_b\}$ in $V_h$, we have
\begin{eqnarray}
c(\bv,\bw,\bw) &=& 0,\label{eq:skew-C}\\
\left|c(\bv,\bw,\bz)\right|&\le& \mathcal{N}_h\3bar\bv\3bar\3bar\bw\3bar\3bar\bz\3bar,\label{eq:Bound-C}
\end{eqnarray}
where $\mathcal{N}_h$ is a constant independent of $h$.
\end{lem}
\begin{pf}
\noindent (A). Proof of (\ref{eq:skew-C}). From the definition of $c(\cdot,\cdot,\cdot)$, one has
\begin{eqnarray*}
c(\bv,\bw,\bw) = (\nabla_w\bv\mRT\bw,\mRT\bw)_{\mathcal{T}_h} - (\nabla_w\bv\mRT\bw,\mRT\bw)_{\mathcal{T}_h},
\end{eqnarray*}
which implies (\ref{eq:skew-C}).

\noindent (B). Proof of (\ref{eq:Bound-C}).
It follows from H\"{o}lder inequality with exponent $(2,4,4)$, Schwartz inequality, %definition of $\3bar\cdot\3bar$, %inequality (\ref{eq:reconstructionEmbed}),
\begin{eqnarray*}
&&\left|\sum_{T\in\mathcal{T}_h}\int_T\nabla_w\bv\mRT\bw\cdot\mRT\bz\right|\\
&\le& \sum_{T\in\mathcal{T}_h}\|\nabla_w\bv\|_{L^2(T)^{3\times 3}}\|\mRT\bw\|_{L^4(T)^3}\|\mRT\bz\|_{L^4(T)^3}\\
%%%
&\le& \left(\sum_{T\in\mathcal{T}_h}\|\nabla_w\bv\|_{L^2(T)^{3\times 3}}^2\right)^{1/2}
\left(\sum_{T\in\mathcal{T}_h}\|\mRT\bw\|_{L^4(T)^3}^2\right)^{1/2}
\left(\sum_{T\in\mathcal{T}_h}\|\mRT\bz\|_{L^4(T)^3}^2\right)^{1/2}\\
%%%
&\le& \3bar\bv\3bar\|\mRT\bw\|_{L^4(\Omega)^3}\|\mRT\bz\|_{L^4(\Omega)^3}\\
%%%
&\lesssim& \3bar\bv\3bar \3bar\bw\3bar \3bar\bz\3bar,
\end{eqnarray*}
where, we used the embeddings (\ref{eq:reconstructionEmbed}) with $r = 4.$ Similarly, switch $\bw$ and $\bz$ in the above estimate, we complete the conclusion in (\ref{eq:Bound-C}).
\end{pf}

With all these preparations, now we are ready to apply the Leray-Schauder fixed point theorem to the WG Algorithm~\ref{alg:WG} and shows the existence and uniqueness of the solution. To this end, we introduce a discrete divergent free subspace $D_h$ of $V_h$ as follows:
\begin{eqnarray*}
D_h=\{\bv\in V_h^0:\nabla_w\cdot\bv = 0\}.
\end{eqnarray*}
%By Helmholtz-Hodege decomposition, we can decompose the body force to divergence free and curl free parts,
%\begin{eqnarray*}
%\bbf = \bg+\nabla\phi,
%\end{eqnarray*}
%where $\bg$ is the curl of a function in $H(curl,\Omega)$ with tangent trace vanishing on $\partial\Omega$, $\phi_1\in H^1(\Omega)$ denotes the irrotational part. For any $\bv\in V_h^0$, we have $(\bbf,\mRT\bv) = (\bg,\mRT\bv)+(\nabla\phi,\mRT\bv) = (\bg,\mRT\bv)+\langle\phi,\mRT\bv\cdot\bn\rangle-(\phi,\nabla\cdot\mRT\bv) = (\bg,\mRT\bv)-(\phi,\nabla_w\cdot\bv)$.
Then the WG formulation can be reformulated as seeking $\bu_h\in D_h$ such that 
\begin{eqnarray}
a(\bu_h,\bv)+c(\bu_h,\bu_h,\bv) = (\bbf,\mRT\bv)_{\mathcal{T}_h},\forall\bv\in D_h.\label{eq:scheme-wg3}
\end{eqnarray} 
Let $F:D_h\to D_h$ be a nonlinear map so that for each $\bw\in D_h,\tilde{\bu}_h:=F(\bw)\in D_h$ is given as the solution of the following linear problem:
\begin{eqnarray}
a(\tilde{\bu}_h,\bv)+c(\bw,\tilde{\bu}_h,\bv) = (\bbf,\mRT\bv)_{\mathcal{T}_h},\forall\bv\in D_h.\label{eq:defineF}
\end{eqnarray} 
The map $F$ is clearly continuous and, therefore, compact in the finite dimensional space $D_h$. If $\lambda>0$ and $\bw$ satisfies $F(\bw) = \lambda\bw$, then from above, we have
\begin{eqnarray*}
\lambda a(\bw,\bv)+\lambda c(\bw,\bw,\bv) = (\bbf,\mRT\bv)_{\mathcal{T}_h},\forall\bv\in D_h.
\end{eqnarray*} 
By choosing $\bv = \bw$ in above equation, we obtain that
\begin{eqnarray}
\lambda ( a(\bw,\bw)+ c(\bw,\bw,\bw) )= (\bbf,\mRT\bw)_{\mathcal{T}_h},\forall\bv\in D_h.
\end{eqnarray} 
It now follows from definition of $\3bar\cdot\3bar$-norm and (\ref{eq:skew-C}), 
\begin{eqnarray}
\lambda \nu\3bar\bw\3bar^2 =  \left|(\bbf,\mRT\bw)_{\mathcal{T}_h}\right|.
\end{eqnarray}
By introducing a mesh-dependent norm
\begin{eqnarray*}
\|\bbf\|_{*,h} = \sum_{\bv\in D_h}\frac{(\bbf,\mRT\bv)_{\mathcal{T}_h}}{\3bar\bv\3bar},
\end{eqnarray*}
and therefore, 
\begin{eqnarray*}
\lambda\le \frac{\|\bbf\|_{*,h}}{\nu\3bar\bw\3bar}.
\end{eqnarray*}
Thus, $\lambda<1$ holds true for any $\bw$ being on the boundary of the ball in $D_h$ centered at the origin with radius $\rho>\dfrac{\|\bbf\|_{*,h}}{\nu}$. Consequently, the Leray-Schauder fixed point theorem implies that the nonlinear map $F$ defined by (\ref{eq:defineF}) has a fixed point $\bu_h$ such that,
\begin{eqnarray*}
F(\bu_h) = \bu_h,
\end{eqnarray*}
in any ball centered at the origin with radius $\rho>\dfrac{\|\bbf\|_{*,h}}{\nu}$. The fixed point $\bu_h$ also is a solution of the finite element scheme (\ref{scheme:wg1})-(\ref{scheme:wg2}), which in turn provides a solution of the original WG Algorithm~\ref{alg:WG}. There can be summarized in the following theorem.
\begin{thm}
The finite element discretization scheme (\ref{scheme:wg1})-(\ref{scheme:wg2}) has at least one solution $\bu_h\in D_h$. Moreover, all the solutions of (\ref{eq:scheme-wg3}) satisfy the following estimates:
\begin{eqnarray}
\3bar\bu_h\3bar\le \frac{\|\bbf\|_{*,h}}{\nu}.\label{eq:boundU}
\end{eqnarray}
\end{thm}
\begin{pf}
Not that $\bu_h\in D_h$ is a solution of (\ref{eq:scheme-wg3}) if and only if it is a fixed-point of the nonlinear map $F$ has at least one fixed point in the ball of $D_h$ centered at the origin with radius $\rho>\dfrac{\|\bbf\|_{*,h}}{\nu}$, then the finite element scheme (\ref{eq:scheme-wg3}) must have a solution and all the solutions mush satisfy the estimate (\ref{eq:boundU}).
\end{pf}

Next, we show the uniqueness of the solution (\ref{eq:scheme-wg3}). Let $\bu_h$ and $\bar{\bu}_h\in D_h$ be two solutions of the finite element scheme (\ref{eq:scheme-wg3}). Since both of them satisfy the nonlinear equation (\ref{eq:scheme-wg3}), let $\phi_h = \bu_h-\bar{\bu}_h$, for all $\bv\in D_h$, we have,
\begin{eqnarray*}
a(\phi_h,\bv)+c(\bu_h,\bu_h,\bv)-c(\bar{\bu}_h,\bar{\bu}_h,\bv) = 0.
\end{eqnarray*}
Observe that,
\begin{eqnarray*}
c(\bu_h,\bu_h,\bv)-c(\bar{\bu}_h,\bar{\bu}_h,\bv) = c(\phi_h,{\bu}_h,\bv)+c(\bar{\bu}_h,\phi_h,\bv).
\end{eqnarray*}
Thus, for any $\bv\in D_h$, we have
\begin{eqnarray*}
a(\phi_h,\bv)+c(\bar{\bu}_h,\phi_h,\bv) = -c(\phi_h,\bar{\bu}_h,\bv).
\end{eqnarray*}
Letting $\bv=\phi_h$, from (\ref{eq:Bound-C}), the fact that $c(\bar{\bu}_h,\phi_h,\phi_h)=0$, and (\ref{eq:Bound-C}), we obtain,
\begin{eqnarray*}
\nu\3bar\phi_h\3bar^2=\left|c(\phi_h,\bar{\bu}_h,\phi_h)\right|\le \mathcal{N}_h\3bar\bar{\bu}_h\3bar\3bar\phi_h\3bar^2.
\end{eqnarray*}
Note that $\bar{\bu}_h$ is a solution of (\ref{eq:scheme-wg3}), therefore, substituting back into the right-hand side of (\ref{eq:boundU}) yields,
\begin{eqnarray}
\nu\3bar\phi_h\3bar^2\le \frac{\mathcal{N}_h\|\bbf\|_{*,h}}{\nu}\3bar\phi_h\3bar^2,
\end{eqnarray}
which implies the uniqueness of the solutions under certain conditions. We summarize the result in the following theorem.
\begin{thm}
Let $\mathcal{N}_h$ be defined in (\ref{eq:boundU}). If $\dfrac{\mathcal{N}_h\|\bbf\|_{*,h}}{\nu^2}<1$ holds, then the WG finite element scheme (\ref{eq:scheme-wg3}) has at most one solution in the discrete divergence-free subspace $D_h.$
\end{thm}

% =================
% Error Estimate
% =================
\section{Main Results}\label{Sect:Error}
In this section, we discuss the convergence results of the WG scheme. We first derive the error equations and then analyze the error estimates. 
Let $\bu_h=\{\bu_0,\bu_b\}\in V_h$ and $p_h\in W_h$ be the solutions to the WG scheme (\ref{scheme:wg1})-(\ref{scheme:wg2}). Let $\bu$ and $p$ be the exact solutions of (\ref{eq:pde-1})-(\ref{eq:pde-bc}). Recall that $Q_h\bu=\{Q_0\bu,Q_b\bu\}$. Similarly, the pressure $p$ is projected onto $W_h$ by $\pi_hp$ with $\pi_hp|_T = \pi_T^kp$. Then the errors $\be_h$ and $\epsilon_h$ for velocity and pressure are defined as follows,
\begin{eqnarray}
\be_h=\{\be_0,\be_b\}=\bQ_h\bu-\bu_h=\{\bQ_0\bu-\bu_0,\bQ_b\bu-\bu_b\},\ \epsilon_h=\pi_hp-p_h.\label{eq:error-define}
\end{eqnarray}
%For vectors $\bu=(u_1,u_2)^\top$ and $\bv=(v_1,v_2)^\top$, we define the following notations,
%\begin{eqnarray*}
%\bu^\top \bv = \begin{pmatrix}
%u_1v_1 & u_1v_2\\ u_2v_1 & u_2v_2
%\end{pmatrix},\ 
%\nabla\bv=\begin{pmatrix}
%\dfrac{\partial v_1}{\partial x} &\dfrac{\partial v_2}{\partial x}\\[6pt]
%\dfrac{\partial v_1}{\partial y} &\dfrac{\partial v_2}{\partial y}
%\end{pmatrix}
%\end{eqnarray*}

\subsection{Error Equation}
%\subsection{Useful Equalities and Inequalities}
First we cite the following integration by parts formula and the estimates for projection operators, which will be used in the error estimate. %The details and prove can be found in \cite{QD2020}.
\begin{lem}
Let $D$ denote a simply connected open polyhedral subset of $\Omega$. For all $\bv,\bw,\bz\in [H^1(D)]^d$, it holds (\cite{QD2020}, Proposition 1)
\begin{eqnarray}
\int_D (\nabla\times\bw)\times\bv\cdot\bz d\bx = \int_D\nabla\bw\bv\cdot\bz d\bx-\int_D\nabla\bw\bz\cdot\bv.\label{eq:integration}
\end{eqnarray}
\end{lem}

%\begin{lem}
%Let $\mathcal{T}_h$ be a finite element partition of $\Omega$ satisfying the shape regularity assumptions and $\bw\in[H^{r+1}(\Omega)]^d$ and $\rho\in H^r(\Omega)$ with $1\le r\le k.$ Then, for $0\le s\le 1$, we have \cite{WY2016}
%\begin{eqnarray} 
%\sum_{T\in\mathcal{T}_h}h_T^{2s}\|\bw-\bQ_0\bw\|_{T,s}^2&\lesssim& h^{2(r+1)}\|\bw\|_{r+1}^2,\label{eq:Q-1}\\
%\sum_{T\in\mathcal{T}_h}h_T^{2s}\|\nabla\bw-\tilde{\pi}_h(\nabla\bw)\|_{T,s}^2&\lesssim& h^{2(r)}\|\bw\|_{r+1}^2,\label{eq:Q-2}\\
%\sum_{T\in\mathcal{T}_h}h_T^{2s}\|\rho-\pi_h\rho\|_{T,s}^2&\lesssim& h^{2(r)}\|\rho\|_{r}^2.\label{eq:Q-3}
%\end{eqnarray}
%%Here $C$ denotes a generic constant independent of the mesh size $h$ and the functions in the estimates.
%\end{lem}

\begin{lem}
For $\bv\in W^{s,r}(T)$ and all $m\in\{0,\cdots,s\}$,
\begin{eqnarray}
&&|\bv-\bQ_0\bv|_{W^{m,r}(T)}\lesssim h_T^{s-m}|\bv|_{W^{s,r}(T)}.
\end{eqnarray}
If $s\ge 1$ and $m\le s-1$, it holds,
\begin{eqnarray}
&&h_T^{1/r}|\bv-\bQ_0\bv|_{W^{m,r}(\partial T)}\lesssim h_T^{s-m}|\bv|_{W^{s,r}(T)}.\label{eq:trace_L4}
\end{eqnarray}
\end{lem}
In the following estimate, we shall use the case with $r = 2$ and $r = 4.$ %When $r = 2,$ one has $L^2(T)=W^{0,2}(T)$.

\begin{lem}
For $v\in H^{k+1}(\Omega)$, $\bv\in [H^{k+2}(\Omega)]^d$, we have (\cite{wy}, Lemma 7.3),
\begin{eqnarray}
%\|{\bf R}_h(\nabla v)-\nabla_w(Q_h v)\|&\le& Ch^{k+1}\|\bv\|_{k+2},\\
%\|\nabla\bv-\nabla_w\bQ_h\bv\|&\le& Ch^{k+1}\|\bv\|_{k+2},\\
\|v-\pi_h v\|&\le& Ch^{k+1}\|v\|_{k+1},\\
\|\bv-{\bf R}_h\bv\|_{L^4(T)^d}&\le& Ch_T^{k+1}\|\bv\|_{k+2},\label{Ineq:Rh}\\
\|\mathbb{R}_h(\nabla \bv)-\nabla_w(\bQ_h\bv)\|&\le& Ch^{k+1}\|\bv\|_{k+2},\label{Ineq:RRh}\\
\|\nabla\bv-\nabla_w\bQ_h\bv\|&\le& Ch^{k+1}\|\bv\|_{k+2}.
\end{eqnarray}
\end{lem}
Then due to the norm equivalence we have the following lemma. 
\begin{lem}\label{lem:Interpolation}
For all $\bv\in H^1(\Omega)^d$, we have
\begin{eqnarray}
\3bar\bQ_h\bv\3bar\le C_I\|\bv\|_1
\end{eqnarray}
with real number $C_I>0$ independent of both $h$ and $\bv.$
\end{lem}

%\begin{lem}
%For all $\bv\in V_h$, we have
%\begin{eqnarray}
%\|\mRT\bv-\bv\|_{0,T}%\le Ch_e^{1/2}\|(\bv_b-\bv_0)\|_{\partial T}
%\le Ch_T\3bar\bv\3bar_1.
%\end{eqnarray}
%\end{lem}

\begin{lem}
Let $\be_h$ and $\epsilon_h$ be defined in (\ref{eq:error-define}). Then we have
\begin{eqnarray}
a(\be_h,\bv) +c(\be_h,\bQ_h\bu,\bv)+c(\bu_h,\be_h,\bv)-b(\epsilon_h,\bv) &=& \ell_\bu(\bv)+\phi_\bu(\bv),\forall\bv\in V_h^0,\label{eq:error-u}\\
b(\be_h,q) &=& 0,\forall q\in W_h,\label{eq:error-p}
\end{eqnarray}
where
\begin{eqnarray*}
%\ell_\bu(\bv) &=& \nu\langle(\mathbb{R}_h(\nabla\bu)-\nabla\bu)\cdot\bn,\bv_0-\bv_b\rangle_{\partial\mathcal{T}_h}-(\nu\Delta\bu-\nu\bm{\pi}_T^{k-1}(\Delta\bu),\bv_0-\mRT\bv)_{\mathcal{T}_h},\\
\ell_\bu(\bv) &=& \nu(\mathbb{R}_h(\nabla\bu)-\PIT(\nabla\bu)),\nabla_w\bv)_{\mathcal{T}_h}-(\nu\Delta\bu-\nu\bm{\pi}_T^{k-1}(\Delta\bu),\bv_0-\mRT\bv)_{\mathcal{T}_h},\\
\phi_\bu(\bv)&=& c(\bQ_h\bu,\bQ_h\bu,\bv)-(\nabla\times\bu\times\bu,\mRT\bv)_{\mathcal{T}_h}.
\end{eqnarray*}
\end{lem}
\begin{pf}

%
%Upper bound. 
%By integration by parts and the property in (), we drive
%\begin{eqnarray*}
%(\nu\Delta\bu,\bv_0)_{\mathcal{T}_h} &=& (\nu\nabla\cdot\nabla\bu,\bv_0)_{\mathcal{T}_h} = -(\nu\nabla\bu,\nabla\bv_0)_{\mathcal{T}_h}+\langle\nu\nabla\bu\cdot\bn,\bv_0-\bv_b\rangle_{\partial\mathcal{T}_h}\\
%%%%
%&=& - (\nu\mathbb{R}_h(\nabla\bu),\nabla\bv_0)_{\mathcal{T}_h}+\langle\nu\nabla\bu\cdot\bn,\bv_0-\bv_b\rangle_{\partial\mathcal{T}_h}\\
%%%%
%&=& - (\nu\nabla_w\bQ_h\bu,\nabla\bv_0)_{\mathcal{T}_h}+\langle\nu\nabla\bu\cdot\bn,\bv_0-\bv_b\rangle_{\partial\mathcal{T}_h}\\
%%%%
%&=&(\nu\nabla\cdot\nabla_w\bQ_h\bu,\bv_0)_{\mathcal{T}_h}-\langle\nu\nabla_w\bQ_h\bu\cdot\bn,\bv_b\rangle_{\partial\mathcal{T}_h}+\langle\nu\mathbb{R}_h(\nabla\bu)\cdot\bn,\bv_b-\bv_0\rangle_{\partial\mathcal{T}_h}\\
%&&+\langle\nu\nabla\bu\cdot\bn,\bv_0-\bv_b\rangle_{\mathcal{T}_h}\\
%%%%
%&=&-\nu(\nabla_w\bQ_h\bu,\nabla_w\bv)_{\mathcal{T}_h}+\nu\langle(\mathbb{R}_h(\nabla\bu)-\nabla\bu)\cdot\bn,\bv_0-\bv_b\rangle_{\partial\mathcal{T}_h}
%\end{eqnarray*}
By integration by parts, the property in (\ref{eq:Rh-operator}), adding and subtracting $(\nu\mathbb{R}_h(\nabla\bu),\nabla_w\bv)_{\mathcal{T}_h}$, and Lemma~\ref{Lem:GradDiv}, we drive
\begin{eqnarray*}
(\nu\Delta\bu,\bv_0)_{\mathcal{T}_h} &=& (\nu\nabla\cdot\nabla\bu,\bv_0)_{\mathcal{T}_h} = -(\nu\mathbb{R}_h(\nabla\bu),\nabla_w\bv)_{\mathcal{T}_h}\\
%%%
&=& -(\nu\PIT(\nabla\bu),\nabla_w\bv)_{\mathcal{T}_h}-
 (\nu\mathbb{R}_h(\nabla\bu)-\nu\PIT(\nabla\bu),\nabla_w\bv)_{\mathcal{T}_h}\\
%%%
&=& - (\nu\nabla_w\bQ_h\bu,\nabla_w\bv)_{\mathcal{T}_h}-
 (\nu\mathbb{R}_h(\nabla\bu)-\nu\PIT(\nabla\bu),\nabla_w\bv)_{\mathcal{T}_h}.
\end{eqnarray*}
Multiply both sides of the equation (\ref{eq:pde-1}) by $\mRT\bv$, by adding and subtracting $(\nu\Delta\bu,\bv_0)$, the above equation, the fact $(\bv_0-\mRT\bv,\bq)_T =0$ for any $\bq\in[\text{P}_{k-1}(T)]^d$, integration by parts, (\ref{eq:DivCW-1}), and the fact $\langle p,\mRT\bv\cdot\bn\rangle_{\partial\mathcal{T}_h} = 0$, it follows,
\begin{eqnarray*}
(\bbf,\mRT\bv)_{\mathcal{T}_h} &=& -(\nu\Delta\bu,\mRT\bv)_{\mathcal{T}_h}+(\nabla\times\bu\times\bu,\mRT\bv)_{\mathcal{T}_h}+(\nabla p,\mRT\bv)_{\mathcal{T}_h}\\
%%%
&=& -(\nu\Delta\bu,\bv_0)_{\mathcal{T}_h} + (\nu\Delta\bu,\bv_0-\mRT\bv)_{\mathcal{T}_h}+(\nabla\times\bu\times\bu,\mRT\bv)_{\mathcal{T}_h}+(\nabla p,\mRT\bv)_{\mathcal{T}_h}\\
%%%
%&=& \nu(\nabla_w\bQ_h\bu,\nabla_w\bv)_{\mathcal{T}_h}-\nu\langle(\mathbb{R}_h(\nabla\bu)-\nabla\bu)\cdot\bn,\bv_0-\bv_b\rangle_{\partial\mathcal{T}_h}\\
%&&+(\nu\Delta\bu-\nu\bm{\pi}_T^{k-1}\Delta\bu,\bv_0-\mRT\bv)_{\mathcal{T}_h}\\
%&&+(\nabla\times\bu\times\bu,\mRT\bv)_{\mathcal{T}_h} - (p,\nabla\cdot\mRT\bv)_{\mathcal{T}_h}+\langle p,\mRT\bv\cdot\bn\rangle_{\partial\mathcal{T}_h}\\
&=& \nu(\nabla_w\bQ_h\bu,\nabla_w\bv)_{\mathcal{T}_h}-\nu(\mathbb{R}_h(\nabla\bu)-\PIT(\nabla\bu),\nabla_w\bv)_{\mathcal{T}_h}\\
&&+(\nu\Delta\bu-\nu\bm{\pi}_T^{k-1}\Delta\bu,\bv_0-\mRT\bv)_{\mathcal{T}_h}\\
&&+(\nabla\times\bu\times\bu,\mRT\bv)_{\mathcal{T}_h} - (p,\nabla\cdot\mRT\bv)_{\mathcal{T}_h}+\langle p,\mRT\bv\cdot\bn\rangle_{\partial\mathcal{T}_h}\\
%%%
&=& \nu(\nabla_w\bQ_h\bu,\nabla_w\bv)_{\mathcal{T}_h}+(\nabla\times\bu\times\bu,\mRT\bv)_{\mathcal{T}_h} - (\pi_hp,\nabla_w\cdot\bv)_{\mathcal{T}_h}\\
%&&-\nu\langle(\mathbb{R}_h(\nabla\bu)-\nabla\bu)\cdot\bn,\bv_0-\bv_b\rangle_{\partial\mathcal{T}_h}+(\nu\Delta\bu-\nu\bm{\pi}_T^{k-1}\Delta\bu,\bv_0-\mRT\bv)_{\mathcal{T}_h}.
&&-\nu(\mathbb{R}_h(\nabla\bu)-\PIT(\nabla\bu),\nabla_w\bv)_{\mathcal{T}_h}
+(\nu\Delta\bu-\nu\bm{\pi}_T^{k-1}\Delta\bu,\bv_0-\mRT\bv)_{\mathcal{T}_h}.
\end{eqnarray*}
Thus, by adding $c(\bQ_h\bu,\bQ_h\bu,\bv)$ to both sides and moving all the other terms to the right hand side, one obtain
\begin{eqnarray*}
&&\nu(\nabla_w\bQ_h\bu,\nabla_w\bv)+c(\bQ_h\bu,\bQ_h\bu,\bv) - (\pi_hp,\nabla_w\cdot\bv)_{\mathcal{T}_h} 
= (\bbf,\mRT\bv)_{\mathcal{T}_h}+\ell_\bu(\bv)+\phi_\bu(\bv).\label{eq:trueU}
\end{eqnarray*}
Furthermore, since
\begin{eqnarray*}
c(\bQ_h\bu,\bQ_h\bu,\bv)-c(\bu_h,\bu_h,\bv) &=& c(\bQ_h\bu-\bu_h,\bQ_h\bu,\bv)+c(\bu_h,\bQ_h\bu,\bv)-c(\bu_h,\bu_h,\bv)\\
&=&c(\be_h,\bQ_h\bu,\bv)+c(\bu_h,\be_h,\bv),
\end{eqnarray*}
and then subtracting (\ref{scheme:wg1}) from (\ref{eq:trueU}) implies,
\begin{eqnarray*}
a(\be_h,\bv) + c(\be_h,\bQ_h\bu,\bv)+c(\bu_h,\be_h,\bv)-b(\epsilon_h,\nabla_w\cdot\bv) = \ell_\bu(\bv)+\phi_\bu(\bv).
\end{eqnarray*}
By multiplying (\ref{eq:pde-2}) by $q\in W_h$ and Lemma~\ref{Lem:GradDiv},
\begin{eqnarray*}
0 = (\nabla\cdot\bu,q)_{\mathcal{T}_h} = (\pi_h(\nabla\cdot\bu),q)_{\mathcal{T}_h} = (\nabla_w\cdot\bQ_h\bu_h,q)_{\mathcal{T}_h},
\end{eqnarray*}
and subtracting (\ref{scheme:wg2}) from above equation completes the proof.
%here
%\begin{eqnarray*}
%\ell_\bu(\bv) = \nu\langle(\Pi_h\nabla\bu-\nabla\bu)\cdot\bn,\bv_0-\bv_b\rangle-(\nu\Delta\bu-\nu\bQ_0^{k-1}\Delta\bu,\bv_0-\mRT\bv)
%\end{eqnarray*}

\end{pf}

\begin{lem}
Assume $\bw\in [H^{k+2}(\Omega)]^d$, we have the following estimates true for $\bv\in V_h$,
\begin{eqnarray}
\left|\ell_{\bw}(\bv)\right|\le C\nu h^{k+1}\|\bw\|_{k+2}\3bar\bv\3bar.\label{eq:estimate-l}
\end{eqnarray}
\end{lem}
\begin{pf}
%Multiply both sides of the equation (\ref{eq:pde-1}) by $\mRT\bv$
%\begin{eqnarray*}
%\mathcal{EE}_h(\bv)&:=&(\bbf,\mRT\bv) - a(\bQ_h\bu,\bv) - c(\bQ_h\bu,\bQ_h\bu,\bv) - b(\bv,\mathbb{Q}_hp)\\
%&=&a(\be_h,\bv)+c(\bu_h,\bu_h,\bv)- c(\bQ_h\bu,\bQ_h\bu,\bv)+b(\bv,\epsilon_h)\\
%%%%
%&=&a(\be_h,\bv)+c(\be_h,\bQ_h\bu,\bv)+b(\bv,\epsilon_h)
%\end{eqnarray*}
%By taking $q = \epsilon_h$, we have
%\begin{eqnarray*}
%b(\be_h,\epsilon_h) = 0.
%\end{eqnarray*}
%Thus,
%\begin{eqnarray*}
%\mathcal{E}_h(\be_h) &=& a(\be_h,\be_h)+c(\be_h,\bQ_h\bu,\be_h) = \nu\3bar\be_h\3bar^2+c(\be_h,\bQ_h\bu,\be_h)\\
%&\ge& (\nu-\mathcal{N}_h\3bar\bQ_h\bu\3bar)\3bar\be_h\3bar^2\\
%%%%
%&\ge&(1-\alpha)\nu C_a\3bar\be_h\3bar^2,
%\end{eqnarray*}
%and here we have used the boundedness property
%\begin{eqnarray}
%\3bar\bQ_h\bu\3bar\le C_I\|\bu\|_1\le\nu^{-1}C_IC_P\|\bg\|.
%\end{eqnarray}
%
%By the definition of $\nabla_w$ and integration by parts, we get
%\begin{eqnarray}
%(\nabla_w(\bQ_h\bw),\nabla_w\bv) = (\nabla\bw,\nabla\bv_0)-\langle\bv_0-\bv_b,\Pi_h(\nabla\bw)\cdot\bn\rangle
%\end{eqnarray}

By Cauchy-Schwarz inequality, trace inequality, definition of $\3bar\cdot\3bar$, and (\ref{Ineq:RRh}), %(\ref{eq:Q-2}), 
we have
\begin{eqnarray*}
%\nu\langle(\mathbb{R}_h(\nabla\bw)-\nabla\bw)\cdot\bn,\bv_0-\bv_b\rangle_{\partial\mathcal{T}_h}
%&\le& \sum_{T\in\mathcal{T}_h}\nu h^{1/2}\|\mathbb{R}_h(\nabla\bw)-\nabla\bw\|_{\partial T}\left(h^{-1/2}\|\bv_0-\bv_b\|_{\partial T}\right)\\
%%%%
%&\le& \nu(\|\mathbb{R}_h(\nabla\bw)-\nabla\bw\|+h\|\mathbb{R}_h(\nabla\bw)-\nabla\bw\|_1) \3bar\bv\3bar_1\\
%&\le& \nu h^{k+1}\|\bw\|_{k+2}\3bar\bv\3bar.
\nu(\mathbb{R}_h(\nabla\bw)-\nabla_w\bQ_h\bw,\nabla_w\bv)&\le& \nu\|\mathbb{R}_h(\nabla\bw)-\nabla_w\bQ_h\bw\| \|\nabla_w\bv\|\\
%%%
&\le&C\nu h^{k+1}\|\bw\|_{k+2}\3bar\bv\3bar.
\end{eqnarray*}
Moreover, by Cauchy-Schwartz inequality, property of $\bm{\pi}_T^{k-1}$, (\ref{eq:Esti-RTv}), and definition of $\3bar\cdot\3bar_1$, it implies,
\begin{eqnarray*}
(\nu\Delta\bw-\nu\bm{\pi}_T^{k-1}\Delta\bw,\bv_0-\mRT\bv)_{\mathcal{T}_h}
&\le& \|\nu\Delta\bw-\nu\bm{\pi}_T^{k-1}\Delta\bw\|\|\bv_0-\mRT\bv\|\\
&\le& \nu h^k\|\bw\|_{k+2}\bigg(\sum_{T\in\mathcal{T}_h}h\|\bv_0-\bv_b\|_{\partial T}^2\bigg)^{1/2}\\
&\lesssim& \nu h^{k+1}\|\bw\|_{k+2}\3bar\bv\3bar_{1}\\
&\lesssim& \nu h^{k+1}\|\bw\|_{k+2}\3bar\bv\3bar.
\end{eqnarray*}
Combining the above two estimates, we completes the proof.
\end{pf}

\begin{lem}
Let $\bw\in [H^{k+2}(\Omega)]^d$ and $\bv\in V_h^0$, then we have
\begin{eqnarray}
\left|\phi_\bw(\bv)\right| %= \left|(\nabla\times\bu\times\bu,\mRT\bv) - c(\bQ_h\bu,\bQ_h\bu,\bv)\right| 
\le Ch^{k+1}\|\bw\|_{k+2}\|\bw\|_2\3bar\bv\3bar.\label{eq:estimate-phi}
\end{eqnarray}
\end{lem}
\begin{pf}
For any $\bv\in V_h^0$, integration by parts (\ref{eq:integration}) gives
\begin{eqnarray*}
(\nabla\times\bw\times\bw,\mRT\bv)_{\mathcal{T}_h} &=& \sum_{T\in\mathcal{T}_h}\int_T(\nabla\times\bw)\times\bw\cdot\mRT\bv d\bx\notag\\
%%%
&=& \sum_{T\in\mathcal{T}_h}\int_T(\nabla\bw\bw\cdot\mRT\bv-\nabla\bw\mRT\bv\cdot\bw)d\bx.
\end{eqnarray*}
Thus, by the above equation and the definition of $c(\cdot,\cdot,\cdot)$ one has,
\begin{eqnarray*}
&&(\nabla\times\bw\times\bw,\mRT\bv)_{\mathcal{T}_h}-c(\bQ_h\bw,\bQ_h\bw,\bv) \notag\\
%%%
&=& \sum_{T\in\mathcal{T}_h}\int_T(\nabla\bw\bw\cdot\mRT\bv-\nabla\bw\mRT\bv\cdot\bw)d\bx\notag\\
&&-\sum_{T\in\mathcal{T}_h}\int_T(\nabla_w\bQ_h\bw\mRT\bQ_h\bw\cdot\mRT\bv-\nabla_w\bQ_h\bw\mRT\bv\cdot\mRT\bQ_h\bw)d\bx\notag\\
%%%
&=& \underbrace{\sum_{T\in\mathcal{T}_h}\int_T(\nabla\bw-\nabla_w\bQ_h\bw)\bw\cdot\mRT \bv d\bx}_{\mathcal{I}_1}
+\underbrace{\sum_{T\in\mathcal{T}_h}\int_T(\nabla_w\bQ_h\bw-\nabla\bw)\mRT\bv\cdot\bw d\bx}_{\mathcal{I}_2}\notag\\
&+&\underbrace{\sum_{T\in\mathcal{T}_h}\int_T\nabla_w\bQ_h\bw(\bw-\mRT\bQ_h\bw)\cdot\mRT\bv d\bx}_{\mathcal{I}_3}
+\underbrace{\sum_{T\in\mathcal{T}_h}\int_T\nabla_w\bQ_h\bw\mRT\bv\cdot(\mRT\bQ_h\bw-\bw) d\bx}_{\mathcal{I}_4}.
\end{eqnarray*}
Next, we shall estimate the terms $\mathcal{I}_1,\dots,\mathcal{I}_4.$

\noindent (A) Estimate of $\mathcal{I}_1$.
First, we add and subtract $\bm{\pi}^0_T\bw$ in the second factor, and then derive
\begin{eqnarray*}
\mathcal{I}_1 &=& \sum_{T\in\mathcal{T}_h}\int_T(\nabla\bw-\nabla_w\bQ_h\bw)(\bw-\bm{\pi}^0_T\bw)\cdot\mRT \bv d\bx+\sum_{T\in\mathcal{T}_h}\int_T(\nabla\bw-\nabla_w\bQ_h\bw)\bm{\pi}^0_T\bw\cdot\mRT \bv d\bx\\
&:=&\mathcal{I}_{1,1}+\mathcal{I}_{1,2}.
\end{eqnarray*}
By H\"{o}lder inequalities with exponents $(2,4,4)$, Cauchy Schwartz inequality, (\ref{Ineq:RRh}), the embedding $H^1(T)\hookrightarrow L^4(T)$, and (\ref{eq:reconstructionEmbed}), we have
\begin{eqnarray*}
\left|\mathcal{I}_{1,1}\right|&\le& \sum_{T\in\mathcal{T}_h}\|\nabla\bw-\nabla_w\bQ_h\bw\|_{L^2(T)^{d\times d}}\|\bw-\bm{\pi}^0_T\bw\|_{L^4(T)^d}\|\mRT \bv \|_{L^4(T)^d}\\
%%%
&\le& \sum_{T\in\mathcal{T}_h}\|\nabla\bw-\nabla_w\bQ_h\bw\|_{L^2(T)^{d\times d}}\sum_{T\in\mathcal{T}_h}\|\nabla(\bw-\bm{\pi}^0_T\bw)\|_{L^2(T)^d}\sum_{T\in\mathcal{T}_h}\|\mRT \bv \|_{L^4(T)^d}\\
%%%
&\le& h^{k+1}|\bw|_{k+2}\|\bw\|_2 \3bar\bv\3bar,
\end{eqnarray*}
where, in the last step, we have used the discrete Sobolev embedding (\ref{eq:reconstructionEmbed}) with $r=4$. By integration by parts, the definition of $\nabla_w$, the identity $\nabla\cdot(\bw\times\bz) = (\nabla\cdot\bz)\bw+\nabla\bw\bz$, the fact $(\bw-\bQ_0\bw,\bq)_T = 0 $ for any $\bq\in [\text{P}_{k}(T)]^d$, the fact $\langle\bw-\bQ_b\bw,\bq\rangle_{\partial T} = 0 $ for any $\bq\in [\text{P}_{k}(\partial T)]^d$, H\"{o}lder inequality with exponent (4,2,4), one obtains,
\begin{eqnarray*}
\mathcal{I}_{1,2} &=& \sum_{T\in\mathcal{T}_h}\int_T(\nabla\bw-\nabla_w\bQ_h\bw)\bm{\pi}_T^0\bw\cdot\mRT \bv d\bx\\
&=&-\sum_{T\in\mathcal{T}_h}\int_T(\bw-\bQ_0\bw)\cdot\nabla\cdot(\mRT\bv\bigotimes\bm{\pi}_T^0\bw) + \sum_{T\in\mathcal{T}_h}\int_{\partial T}(\bu-\bQ_b\bu)\cdot(\mRT\bv\bigotimes\bm{\pi}_T^0\bw)\bn ds\\
%%%
&=& -\sum_{T\in\mathcal{T}_h}\int_T(\bw-\bQ_0\bw)\cdot(\nabla\mRT\bv\bm{\pi}_T^0\bw) + \sum_{T\in\mathcal{T}_h}\int_{\partial T}(\bu-\bQ_b\bu)\cdot(\mRT\bv\bigotimes\bm{\pi}_T^0\bw)\bn ds\\
%%%
&=&\sum_{T\in\mathcal{T}_h}\int_{\partial T}(\bw-\bQ_b\bw)\cdot(\mRT\bv\bigotimes\bm{\pi}_T^0\bw)\bn ds\\
%%%
&=&\sum_{T\in\mathcal{T}_h}\int_{\partial T}(\bw-\bQ_b\bw)\cdot((\mRT\bv-\bm{\pi}_T^0\mRT\bv)\bigotimes\bm{\pi}_T^0\bw)\bn ds\\
%%%
&\le& \sum_{T\in\mathcal{T}_h} \|\bw-\bQ_b\bw\|_{L^4(\partial T)^3}\|\mRT\bv-\bm{\pi}_T^0\mRT\bv\|_{L^2(\partial T)^3}\|\bm{\pi}_T^0\bw\|_{L^4(\partial T)^3}\\
%%%
&\le& Ch^{k+1}\|\bw\|_{k+2}\3bar\bv\3bar\|\bw\|_2,
\end{eqnarray*}
where we have used the inequalities (\ref{eq:phi-1})-(\ref{eq:phi-3}) below.

By adding and subtracting $\bQ_0\bw$, property of $\bQ_b$, trace inequality (\ref{eq:trace_L4}), and the embedding $W^{k+1,4}(T)\hookrightarrow H^{k+2}(T),$ we have,
\begin{eqnarray}
\|\bw-\bQ_b\bw\|_{L^4(\partial T)^3}&\le& \|\bw-\bQ_0\bw\|_{L^4(\partial T)^3}+\|\bQ_b(\bQ_0\bw-\bw)\|_{L^4(\partial T)^3}\notag\\
&\le&\|\bw-\bQ_0\bw\|_{L^4(\partial T)^3}\notag\\
&\le& h_T^{k+1-1/4}\|\bw\|_{W^{k+1,4}(T)}\lesssim h_T^{k+1-1/4}\|\bw\|_{k+2}.\label{eq:phi-1}
\end{eqnarray}
By the property of projection $\bm{\pi}_T^0$, trace inequality, adding and subtracting $\bv_0$, inverse inequality, (\ref{eq:reconstructionEmbed}), (\ref{eq:Esti-RTv}), (\ref{eq:3bar-1}), and (\ref{eq:3bar-3bar1}), it implies
\begin{eqnarray}
\|\mRT\bv-\bm{\pi}_T^0\mRT\bv\|_{L^2(\partial T)^3}&\le& h_T^{1/2}\|\nabla\mRT\bv\|_{T}\notag\\
&\le& h_T^{1/2}\left(\|\nabla(\mRT\bv-\bv_0)\|_{T}+\|\nabla\bv_0\|_T\right)\notag\\
&\le& Ch_T^{1/2}\left(h_T^{-1}\|\mRT\bv-\bv_0\|+\|\nabla\bv_0\|_T\right)\notag\\
&\le& Ch_T^{1/2}\left(h_T^{-1/2}\|\bv_b-\bv_0\|_{\partial T}+\|\nabla\bv_0\|_T\right)\notag\\
&\le& Ch_T^{1/2}\3bar\bv\3bar.\label{eq:phi-2}
\end{eqnarray}
By $L^4$-boundedness of $\bm{\pi}_T^0$, trace inequality (\ref{eq:trace_L4}), and the embedding $W^{1,4}(T)\hookrightarrow H^2(T)$, we have
\begin{eqnarray}
\|\bm{\pi}_T^0\bw\|_{L^4(\partial T)^3}\le \|\bw\|_{L^4(\partial T)^3}\le h_T^{-1/4}\|\bw\|_{W^{1,4}(T)}\le h_T^{-1/4}\|\bw\|_{2,T}.\label{eq:phi-3}
\end{eqnarray}

% I_2
\noindent(B). Estimate of $\mathcal{I}_2$. We add $\bm{\pi}_T^0\bw$ to the third factor, and then
\begin{eqnarray}
\mathcal{I}_2 &=& \sum_{T\in\mathcal{T}_h}\int_T(\nabla_w\bQ_h\bw-\nabla\bw)\mRT\bv\cdot(\bw-\bm{\pi}_T^0\bw) d\bx+ \sum_{T\in\mathcal{T}_h}\int_T(\nabla_w\bQ_h\bw-\nabla\bw)\mRT\bv\cdot\bm{\pi}_T^0\bw d\bx\notag\\
&:=&\mathcal{I}_{2,1}+\mathcal{I}_{2,2}
\end{eqnarray}
By H\"{o}lder inequalities with exponents (2,4,4), bound (\ref{eq:reconstructionEmbed}), embedding $L^4(T)\hookrightarrow H^1(T)$, we have,
\begin{eqnarray*}
\left|\mathcal{I}_{2,1}\right| &=& \sum_{T\in\mathcal{T}_h}\|\nabla_w\bQ_h\bw-\nabla\bw\|_{L^2(T)^{d\times d}} \|\mRT\bv\|_{L^4(T)^d} \|\bw-\bm{\pi}_T^0\bw\|_{L^4(T)^d}\\
%%%
&\le& Ch^{k+1}\|\bw\|_{k+2}\3bar\bv\3bar\|\bw\|_2
\end{eqnarray*}
Then, we can rewrite
\begin{eqnarray*}
\mathcal{I}_{2,2} = \sum_{T\in\mathcal{T}_h}\int_T(\nabla_w\bQ_h\bw-\nabla\bw)\mRT\bv\cdot\bm{\pi}_T^0\bw d\bx = \sum_{T\in\mathcal{T}_h}\int_T(\nabla_w\bQ_h\bw-\nabla\bw):\bm{\pi}_T^0\bw\bigotimes\mRT\bv d\bx.
\end{eqnarray*}
By integration by parts, the fact $\bm{\pi}_T^0\bw\bigotimes\mRT\bv\in [\mathbb{RT}^k(T)]^{d\times d}$, the definition of $\nabla_w$, (\ref{eq:Def-mRT-2}), and the property of projection operator $\bQ_0$ and $\bQ_b$, we get
\begin{eqnarray*}
\mathcal{I}_{2,2} &=& -\sum_{T\in\mathcal{T}_h}\int_T(\bQ_0\bw-\bw)\cdot\nabla\cdot(\bm{\pi}_T^0\bv\bigotimes\mRT\bv) d\bx+\sum_{T\in\mathcal{T}_h}\int_{\partial T}(\bQ_b\bw-\bw)\cdot(\bm{\pi}_T^0\bw\bigotimes\mRT\bv)\bn ds\\
%%%
&=&-\sum_{T\in\mathcal{T}_h}\int_T(\bQ_0\bw-\bw)\cdot(\bm{\pi}_T^0\bw\nabla\cdot \mRT\bv) d\bx+\sum_{T\in\mathcal{T}_h}\int_{\partial T}(\bQ_b\bw-\bw)\cdot(\bm{\pi}_T^0\bw(\mRT\bv\cdot\bn)) ds\\
%%%
&=&-\sum_{T\in\mathcal{T}_h}\int_T(\bQ_0\bw-\bw)\cdot(\bm{\pi}_T^0\bw\nabla\cdot \mRT\bv) d\bx+\sum_{T\in\mathcal{T}_h}\int_{\partial T}(\bQ_b\bw-\bw)\cdot(\bm{\pi}_T^0\bw(\bv_b\cdot\bn)) ds\\
&=& 0.
\end{eqnarray*}
%%%

% III
\noindent(C). Estimate of $\mathcal{I}_3$. By the fact $\mRT\bQ_h\bw = {\bf R}_h\bw$, H\"{o}lder inequality with exponent $(2,4,4)$, (\ref{Ineq:Rh}), the embedding $W^{k+1,4}(T)\hookrightarrow H^{k+2}(T)$, and (\ref{eq:reconstructionEmbed}), it follows,
\begin{eqnarray*}
\mathcal{I}_3 &=& \sum_{T\in\mathcal{T}_h}\int_T\nabla_w\bQ_h\bw(\bw-\mRT\bQ_h\bw)\cdot\mRT\bv d\bx\\
%%%
&=& \sum_{T\in\mathcal{T}_h}\int_T\nabla_w\bQ_h\bw(\bw-{\bf R}_h\bw)\cdot\mRT\bv d\bx\\
%%%
&\le& \sum_{T\in\mathcal{T}_h}\|\nabla_w\bQ_h\bw\| _{L^2(T)^{d\times d}}\|\bw-{\bf R}_h\bw\|_{L^4(T)^d}\|\mRT\bv\|_{L^4(T)^d}\\
%%%
&\le&Ch^{k+1}\3bar\bQ_h\bw\3bar\|\bw\|_{k+2}\3bar\bv\3bar\\
&\le&Ch^{k+1}\|\bw\|_2\|\bw\|_{k+2}\3bar\bv\3bar.
\end{eqnarray*}

% IV
\noindent(D). Estimate of $\mathcal{I}_4$. By the fact $\mRT\bQ_h\bw = {\bf R}_h\bw$, H\"{o}lder inequality with exponent $(2,4,4)$,  (\ref{eq:reconstructionEmbed}), (\ref{Ineq:Rh}), and the embedding $W^{k+1,4}(T)\hookrightarrow H^{k+2}(T)$, it follows,

\begin{eqnarray*}
\mathcal{I}_4 &=& \sum_{T\in\mathcal{T}_h}\int_T\nabla_w\bQ_h\bw\mRT\bv\cdot(\mRT\bQ_h\bw-\bw) d\bx\\
%%%
&=& \sum_{T\in\mathcal{T}_h}\int_T\nabla_w\bQ_h\bw\mRT\bv\cdot({\bf R}_h\bw-\bw) d\bx\\
%%%
&\le& \sum_{T\in\mathcal{T}_h}\|\nabla_w\bQ_h\bw\|_{L^2(T)^{d\times d}}\|\mRT\bv\|_{L^4(T)^d}\|{\bf R}_h\bw-\bw\|_{L^4(T)^d}\\
%%%
&\le& Ch^{k+1}\3bar\bQ_h\bw\3bar\3bar\bv\3bar\|\bw\|_{k+2}\\
&\le& Ch^{k+1}\|\bw\|_2\3bar\bv\3bar\|\bw\|_{k+2}.
\end{eqnarray*}
Lastly, by combining all above we complete the proof.
\end{pf}

\subsection{Error Estimates}
First, by Helmholtz decomposition, we can denote $\bbf = \bg+\nabla\psi$, where $\bg$ is the curl of a function in ${\bf H}(\text{curl};\Omega)$ whose tangent trace vanishes on $\partial\Omega$ and $\psi\in H^1(\Omega)$. As shown in \cite{QD2020}, by taking $\bv=\bu$, $q=p-\psi$ in (\ref{eq:weakform-1}) and (\ref{eq:weakform-2}), one has
\begin{eqnarray*}
\nu\|\nabla\bu\|^2 = (\bg,\bu)\le \|\bg\|\|\bu\|\le C\|\bg\|\|\nabla\bu\|,
\end{eqnarray*}
where we have used Poincare theorem in the last step. Thus, the exact solution is only bounded by solenoidal part of $\bbf$:
\begin{eqnarray}
\|\nabla\bu\|\lesssim \nu^{-1}\|\bg\|.\label{eq:exact-u}
\end{eqnarray}

Now, combining all the estimates above, we have the following overall error estimates of the WG scheme (\ref{scheme:wg1})-(\ref{scheme:wg2}).
\begin{thm}\label{thm:main}
We assume the following holds, for the solenoidal component of the body force $\bbf$, 
\begin{eqnarray}
\|\bg\|\le C_3\nu^2, %\alpha\frac{\nu^2}{\mathcal{N}_hC}.
\end{eqnarray}
where the constant $C_3$ will be specified later.
Let $(\bu,p)$ are the solutions to (\ref{eq:weakform-1})-(\ref{eq:weakform-2}), $\bu\in [H_0^1(\Omega)\cap H^{k+2}(\Omega)]^d$, $k\ge 0$, and $(\bu_h,p_h)\in V_h\times W_h$ be the solutions of (\ref{scheme:wg1})-(\ref{scheme:wg2}), respectively. %If $\rho=\dfrac{\mathcal{N}_h\|\bbf\|_{*,h}}{\nu^2}<1$, 
Then the following error estimate holds
\begin{eqnarray}
\3bar\be_h\3bar&\le& Ch^{k+1}(1+\nu^{-1}\|\bu\|_2)\|\bu\|_{k+2}\\
\|\epsilon_h\|&\le&C h^{k+1}(\nu+\|\bu\|_2)\|\bu\|_{k+2}
\end{eqnarray}
\end{thm}
\begin{pf}
\noindent(A). Estimate on the velocity error $\be_h$.
Let $q = \epsilon_h$ and $\bv=\be_h$ in the error equation (\ref{eq:error-u}), and one gets
\begin{eqnarray*}
a(\be_h,\be_h)+c(\be_h,\bQ_h\bu,\bv) = \ell_\bu(\be_h)+\phi_\bu(\be_h).
\end{eqnarray*}
By (\ref{eq:Bound-C}), Lemma~\ref{lem:Interpolation}, bound (\ref{eq:exact-u}), and thus $\3bar\bQ_h\bu\3bar\le C_I\|\nabla\bu\|\le C\nu^{-1}\|\bg\|$, where the constant $C$ contains the interpolation and Poincare constants.  It follows,
\begin{eqnarray*}
\text{Left} &:=& a(\be_h,\be_h)+c(\be_h,\bQ_h\bu,\be_h) = \nu\3bar\be_h\3bar^2+c(\be_h,\bQ_h\bu,\be_h)\\
&\ge& (\nu-\mathcal{N}_h\3bar\bQ_h\bu\3bar)\3bar\be_h\3bar^2\\
%%%
&\ge& (\nu-\nu^{-1}\mathcal{N}_hC\|\bg\|)\3bar\be_h\3bar^2\\
%%%
&\ge&C_2\nu \3bar\be_h\3bar^2,
\end{eqnarray*}
here we have used $\|\bg\|\le C_3\nu^2$, where $C_3 = \dfrac{(1-C_2)}{\mathcal{N}_hC}$.
Combining (\ref{eq:estimate-l})-(\ref{eq:estimate-phi}), we obtain,
\begin{eqnarray*}
\text{Right} := \ell_\bu(\be_h)+\phi_\bu(\be_h)\le C\nu h^{k+1}\|\bu\|_{k+2}\3bar\be_h\3bar+Ch^{k+1}\|\bu\|_{k+2}\|\bu\|_2\3bar\be_h\3bar,
\end{eqnarray*}
and thus all above imply
\begin{eqnarray*}
\3bar\be_h\3bar\le Ch^{k+1}(1+\nu^{-1}\|\bu\|_2)\|\bu\|_{k+2}.
\end{eqnarray*}

\noindent(B). Estimate on the pressure error $\epsilon_h$. By inf-sup condition Lemma~\ref{lem:InfSup}, error equation (\ref{eq:error-u}), Lemma~\ref{lem:bilinearA}, (\ref{eq:estimate-l})-(\ref{eq:estimate-phi}), $\|\bg\|\lesssim \nu^2$, and $\3bar\bQ_h\bu\3bar\le C_I\|\nabla\bu\|\le C\nu^{-1}\|\bg\|\lesssim \nu$, we conclude,
\begin{eqnarray*}
\|\epsilon_h\|&\le& \sup_{\bv\in V_h^0}\frac{\left|b(\bv,\epsilon_h)\right|}{\3bar\bv\3bar}
 = \sup_{\bv\in V_h^0}\frac{\left|-\nu a(\be_h,\bv)-c(\be_h,\bQ_h\bu,\bv)+\ell_\bu(\bv)+\phi_\bu(\bv)\right|}{\3bar\bv\3bar}\\
 %%%
 &\le& \nu\3bar\be_h\3bar+\3bar\be_h\3bar\3bar\bQ_h\bu\3bar+\nu h^{k+2}\|\bu\|_{k+2}+Ch^{k+1}\|\bu\|_{k+2}\|\bu\|_2\\
 %%%
 &\le&C h^{k+1}(\nu+\|\bu\|_2)\|\bu\|_{k+2}.
\end{eqnarray*}
\end{pf}

\begin{remark}\label{thm:main-2}
In comparison, by Algorithm~\ref{alg:WG-2}, one can derive the error estimate using the similar argument. Here we omit the proof but only introduce the following error estimates. Let $(\bu,p)\in [H_0^1(\Omega)\cap H^{k+2}(\Omega)]^d\times (L_0^2(\Omega)\cap H^{k+1}(\Omega))$, $k\ge 0$, and $(\hat{\bu}_h,\hat{p}_h)\in V_h\times W_h$ be the solution of (\ref{eq:pde-1})-(\ref{eq:pde-bc}) and weak Galerkin Algorithm~\ref{alg:WG-2}, respectively. If we have $\|\bbf\|\lesssim \dfrac{\nu^2}{\mathcal{N}_hC}$, the following error estimates hold:
\begin{eqnarray}
\3bar\bQ_h\bu-\hat{\bu}_h\3bar&\le& Ch^{k+1}\left( (1+\nu^{-1}\|\bu\|_2)\|\bu\|_{k+2}+\nu^{-1}\|p\|_{k+1}\right),\label{eq:error-u-wg2}\\
\|\pi_hp-\hat{p}_h\|&\le&C h^{k+1}\left( (\nu+\|\bu\|_2)\|\bu\|_{k+2}+\|p\|_{k+1}\right).
\end{eqnarray}
\end{remark}

\begin{remark}
The error estimates in Theorem~\ref{thm:main} and Remark~\ref{thm:main-2} indicate that:
\begin{itemize}
\item Error estimate in Theorem~\ref{thm:main} is established under a data smallness condition which only involves the solenoidal part of the body force; the errors are independent of the pressure. 
\item Error estimate in Remark~\ref{thm:main-2} is established under a data smallness condition which depends of the full body force; the errors are depending on the velocity and pressure.  
\item Due to the independence of irrotational body force and pressure, Algorithm~\ref{alg:WG} shows the robustness with respect to irrotational body force and pressure.
\end{itemize}
\end{remark}

%\begin{thm}
%Assume ? are fulfilled. For $\bbf\in L^2(\Omega)^2$, let $\bu$ and $\bu_h$ be the solutions of (\ref{eq:pde-1})-(\ref{eq:pde-bc}) and (\ref{scheme:wg1})-(\ref{scheme:wg2}). Then
%\begin{eqnarray}
%\|\bu-\bu_h\|\le C\min\{\frac{h^2}{\nu},\frac{1}{\alpha+\|w\|_\infty}\}\|\bbf\|
%\end{eqnarray}
%holds with a constant $C$ independent of $\nu$, $\alpha$, $w$, $h$, and $\bbf.$
%\end{thm}

% =================
% Numerical
% =================
\section{Numerical Experiment}\label{Sect:NumTest}
In this section, we present several two dimensional numerical experiments to demonstrate the effectiveness of the proposed scheme. Test~\ref{Sect:Num-1}-Test~\ref{Sect:Num-6} will be performed on the uniform triangular mesh, with mesh size $h$. In all the numerical experiments, we use Newton's method to linearize the nonlinear discrete problem:
%Newton's method will be performed for the nonlinear solver. We substitute $\bu_h = \bu_h+\delta\bu_h$ and $p_h = p_h+\delta p_h$, into the WG algorithm:
%\begin{eqnarray*}
%\nu a(\bu_h+\delta\bu_h,\bv)+c(\bu_h+\delta\bu_h,\bu_h+\delta\bu_h,\bv)+b(\bv,p_h+\delta p_h) = (\bbf,\mRT(\bv)).
%\end{eqnarray*}
%By splitting trilinear term $c(\cdot,\cdot,\cdot)$ and ignoring the quadratic term in $\delta\bu_h$, we obtain
%\begin{eqnarray*}
%\nu a(\bu_h+\delta\bu_h,\bv)+c(\bu_h+\delta\bu_h,\bu_h,\bv)+c(\bu_h,\delta\bu_h,\bv)+b(\bv,p_h+\delta p_h) = (\bbf,\mRT(\bv)),
%\end{eqnarray*}
%and together with adding and subtracting $c(\bu_h,\bu_h,\bv)$, 
%\begin{eqnarray*}
%\nu a(\bu_h+\delta\bu_h,\bv)+c(\bu_h+\delta\bu_h,\bu_h,\bv)+c(\bu_h,\bu_h+\delta\bu_h,\bv)+b(\bv,p_h+\delta p_h) = (\bbf,\mRT(\bv))+c(\bu_h,\bu_h,\bv).
%\end{eqnarray*}
%Denote $\bu_h^{n} := \bu_h$, $\bu_h^{n+1} := \bu_h+\delta\bu_h$, $p_h^{n}:=p_h$, and $p_h^{n+1}:=p_h+\delta p_h$.
\begin{eqnarray}
\nu a(\bu_h^{n+1},\bv^{n+1})+c(\bu_h^{n},\bu_h^{n+1},\bv)+c(\bu_h^{n+1},\bu_h^n,\bv)+b(\bv,p_h^{n+1}) = (\bbf,\mRT(\bv))+c(\bu_h^n,\bu_h^n,\bv).\label{eq:Newton}
\end{eqnarray}
In the Newton solver, we shall take $(\bu_h^0,p_h^0)$ as the solution for the corresponding Stokes solver and then proceed the iteration (\ref{eq:Newton}) until stopping criterion $|(\bu_h^{n+1},p_h^{n+1})^\top-(\bu_h^{n},p_h^n)^\top|<1\text{E}$-10 or \#Iteration $>$ 1E3 is satisfied.
%Here
%\begin{eqnarray*}
%c(\bu^n,\bu^{n+1},\bv) &=& (\nabla_w\bu^n\bu^{n+1},\bv) - (\nabla_w\bu^n\bv,\bu^{n+1}) \\
%&=& \int_T\nabla_w\bu^n:(\mRT(\bu^{n+1})\bigotimes\mRT\bv)dT-\int_T\nabla_w\bu^n:(\mRT\bv\bigotimes\mRT(\bu^{n+1}))dT\\
%%%%
%&=&\int_T \begin{bmatrix}
%\nabla_w \bu_1,\nabla_w \bu_2
%\end{bmatrix}: 
%\left(
%\begin{bmatrix}
%\mRT\bu^{n+1}_1\\ \mRT\bu^{n+1}_2
%\end{bmatrix} 
%\bigotimes
%\begin{bmatrix}
%\mRT\bv_1\\ \mRT\bv_2
%\end{bmatrix} 
%\right)
%dT\\
%&-&\int_T \begin{bmatrix}
%\nabla_w \bu_1,\nabla_w \bu_2
%\end{bmatrix}: 
%\left(
%\begin{bmatrix}
%\mRT\bv_1\\ \mRT\bv_2
%\end{bmatrix} 
%\bigotimes
%\begin{bmatrix}
%\mRT\bu^{n+1}_1\\ \mRT\bu^{n+1}_2
%\end{bmatrix} 
%\right)
%dT\\
%%%%
%&=&\int_T \begin{bmatrix}
%\nabla_w \bu_1,\nabla_w \bu_2
%\end{bmatrix}: 
%\left(
%\begin{bmatrix}
%\mRT\bu^{n+1}_1\mRT\bv_1 &\mRT\bu^{n+1}_1\mRT\bv_2\\ 
%\mRT\bu^{n+1}_2\mRT\bv_1 &\mRT\bu^{n+1}_2\mRT\bv_2
%\end{bmatrix} 
%\right)dT\\
%&-&\int_T \begin{bmatrix}
%\nabla_w \bu_1,\nabla_w \bu_2
%\end{bmatrix}: 
%\left(
%\begin{bmatrix}
%\mRT\bv_1\mRT\bu^{n+1}_1 &\mRT\bv_1\mRT\bu^{n+1}_2\\ 
%\mRT\bv_2\mRT\bu^{n+1}_1 &\mRT\bv_2\mRT\bu^{n+1}_2
%\end{bmatrix} 
%\right)dT\\
%%%%
%&=&\int_T \begin{bmatrix}
%\nabla_w \bu_1,\nabla_w \bu_2
%\end{bmatrix}: 
%\left(
%\begin{bmatrix}
%0 &\mRT\bu^{n+1}_1\mRT\bv_2\\ 
%\mRT\bu^{n+1}_2\mRT\bv_1 &0
%\end{bmatrix} 
%\right)dT
%\end{eqnarray*}

\subsection{Convergence Test}\label{Sect:Num-1}
In this test, we choose smooth exact solutions and report the convergence results.
Let $\Omega = (0,1)^2$ and the exact solutions are chosen as follows
\begin{eqnarray*}
\bu = \begin{pmatrix}
10x^2y(x-1)^2(2y-1)(y-1)\\
-10xy^2(2x-1)(x-1)(y-1)^2
\end{pmatrix},\quad p = 10(2x-1)(2y-1).
\end{eqnarray*}

We perform the WG simulation by Algorithm~\ref{alg:WG} and Algorithm~\ref{alg:WG-2} for various viscosity values and WG element, including $\nu = 1,1\text{E-}2,1\text{E-}4$ and $k = 0,1,2$ and the numerical results are present ed in Table~\ref{Tab:Num-1_k0}-\ref{Tab:Num-1_k2}. 
We observe:
\begin{itemize}
\item For viscosity value $\nu=1,$ 1E-2, 1E-4, both algorithms produce an approximation with convergence rate $\mathcal{O}(h^{k+1}),\mathcal{O}(h^{k+2}),\mathcal{O}(h^{k+1})$ for the velocity errors measured in $\3bar\cdot\3bar$-norm, $\|\cdot\|$-norm, and pressure error measured in $L^2$-norm, as long as the nonlinear solver converges. 
\item For viscosity value $\nu=1$ and low order polynomial $k=0,1$, WG Algorithm~\ref{alg:WG} produce better numerical results with 8X improvement in velocity and 6X improvement in pressure than that of Algorithm~\ref{alg:WG-2}. But when $k=2$, Algorithm~\ref{alg:WG-2} produce slightly better velocity error and the same order pressure approximation comparing to that from Algorithm~\ref{alg:WG}. The velocity produced by Algorithm~\ref{alg:WG-2} has 2X improvement. This is because that by employing WG element $k = 2$, the pressure (quadratic function in this case) is fully resolved by WG Algorithm~\ref{alg:WG-2} and will not affect the velocity approximation. But the inconsistent error by Algorithm~\ref{alg:WG} is dominate.
\item For viscosity value $\nu = 1$E-2, the velocity approximation in Algorithm~\ref{alg:WG-2} is deteriorated for $k = 0,1$. One can observe that velocity error is increased by $100$ times comparing to that of $\nu = 1$. The pressure error remains the same. In contrary, by Algorithm~\ref{alg:WG}, the velocity error remains the same but the pressure error is reduced by $100$ times comparing to that from $\nu = 1.$ However, for $k=2$, as the pressure is fully resolved, Algorithm~\ref{alg:WG-2} produce the similar simulation as that from Algorithm~\ref{alg:WG}.
\item For viscosity value $\nu=$1E-4, Algorithm~\ref{alg:WG-2} with low polynomial order $k=0$ and 1 does not converge within the maximum iteration number. Contrarily, Algorithm~\ref{alg:WG-2} is still robust to produce better numerical results, which has the same order in velocity error and $1/\nu$ smaller pressure error. Moreover, one can observe that, the pressure approximation converges at the super-convergence rate $\mathcal{O}(h^{k+2})$. When the pressure variable is fully resolved by WG element $k = 2$, both of Algorithm~\ref{alg:WG} and Algorithm~\ref{alg:WG-2} can provide good numerical solution, while the error through Algorithm~\ref{alg:WG-2} is slightly smaller.
%\item As one reduces the value in $\nu$, the performance of WG Algorithm 1 is deteriorated. 
%\item It can be seen from these tables as low order polynomials, WG Algorithm~\ref{alg:WG} produces better numerical results.
\end{itemize}
All the observations agree with our theoretical conclusions in Theorem~\ref{thm:main} and Remark~\ref{thm:main-2}, and thus show the robustness of Algorithm~\ref{alg:WG}.

{
\begin{table}
\caption{Example~\ref{Sect:Num-1}: Error Profiles and Convergence Results for $k=0$.}\label{Tab:Num-1_k0}
\centering
\tabcolsep=3pt
\begin{tabular}{c|cc|cc|cc||cc|cc|cc}
\hline\hline
 &\multicolumn{6}{c||}{WG Algorithm~\ref{alg:WG-2}} & \multicolumn{6}{c}{WG Algorithm~\ref{alg:WG}} \\
$1/h$	&$\3bar\be_h\3bar$	&Rate	&$\|\be_0\|$ &Rate		&$\|\epsilon_h\|$ &Rate	&$\3bar\be_h\3bar$	&Rate	&$\|\be_0\|$ &Rate		&$\|\epsilon_h\|$ &Rate	\\ \hline\hline
%k = 0, nu = 1		
\multicolumn{13}{c}{$k = 0, \nu = 1$}\\ \hline										
16	&3.26E-1	& 	&8.72E-3	 &	&6.31E-1	&	&5.73E-2	& 	&1.10E-3	& 	&1.17E-2	&\\
32	&1.68E-1	&0.95	&2.34E-3	&1.90	&3.00E-1	&1.07	&2.89E-2	&0.99	&2.85E-4	&1.95	&5.32E-3	&1.14\\
64	&8.50E-2	&0.98	&5.99E-4	&1.97	&1.44E-1	&1.06	&1.45E-2	&1.00	&7.18E-5	&1.99	&2.57E-3	&1.05\\
128	&4.27E-2	&0.99	&1.51E-4	&1.99	&6.98E-2	&1.04	&7.23E-3	&1.00	&1.80E-5	&2.00	&1.28E-3	&1.01\\ \hline
\multicolumn{13}{c}{$k = 0, \nu = 1\text{E}-2$}\\ \hline
%k = 0, nu = 1e-2												
16	&29.4	& 	&7.58E-1	& 	&1.42	&	&5.73E-2	& 	&1.11E-3	 &	&1.73E-4	& \\
32	&16.7	&0.82	&2.56E-1	&1.56	&3.20E-01	&2.15	&9.80E-5	&0.99	&2.85E-4	&1.95	&6.24E-5	&1.47\\
64	&8.46	&0.98	&6.16E-2	&2.06	&1.43E-01	&1.16	&1.24E-5	&1.00	&7.20E-5	&1.99	&2.70E-5	&1.21\\
128	&4.25	&0.99	&1.52E-2	&2.02	&6.97E-02	&1.04	&1.56E-6	&1.00	&1.80E-5	&2.00	&1.29E-5	&1.06\\ \hline
%k = 0, nu = 1e-4	
\multicolumn{13}{c}{$k = 0, \nu = 1\text{E}-4$}\\ \hline											
16	&- &-&-&-&-&-						&6.14E-2	& 	        &1.63E-3	 &	        &4.36E-5	&\\
32	&- &-&-&-&-&-						&2.93E-2	&1.07	&3.94E-4	&2.05	&1.07E-5	&2.03\\
64	&- &-&-&-&-&-						&1.45E-2	&1.02	&9.69E-5	&2.03	&2.60E-6	&2.04\\
128	&- &-&-&-&-&-						&7.24E-3	&1.00	&2.42E-5	&2.00	&6.59E-7	&1.98\\ \hline
\end{tabular}
\end{table}

\begin{table}
\caption{Example~\ref{Sect:Num-1}: Error Profiles and Convergence Results for $k=1$.}\label{Tab:Num-1_k1}
\centering
\tabcolsep=3pt
\begin{tabular}{c|cc|cc|cc||cc|cc|cc}
\hline\hline
 &\multicolumn{6}{c||}{WG Algorithm~\ref{alg:WG-2}} & \multicolumn{6}{c}{WG Algorithm~\ref{alg:WG}} \\
$1/h$	&$\3bar\be_h\3bar$	&Rate	&$\|\be_0\|$ &Rate		&$\|\epsilon_h\|$ &Rate	&$\3bar\be_h\3bar$	&Rate	&$\|\be_0\|$ &Rate		&$\|\epsilon_h\|$ &Rate	\\ \hline\hline
\multicolumn{13}{c}{$k = 1, \nu = 1$}\\ \hline	
%k = 1, nu = 1												
16	&1.05E-2	& 	&9.21E-5	 &	&2.04E-2	&	&2.37E-3	& 	&1.98E-5	& 	&9.27E-4	&\\
32	&2.66E-3	&1.99	&1.17E-05	&2.98	&5.09E-3	&2.00	&6.00E-4	&1.98	&2.52E-6	&2.98	&2.23E-4	&2.05\\
64	&6.68E-4	&1.99	&1.47E-06	&2.99	&1.27E-3	&2.00	&1.51E-4	&1.99	&3.18E-7	&2.99	&5.44E-5	&2.04\\
128	&1.67E-4	&2.00	&1.85E-07	&2.99	&3.18E-4	&2.00	&3.78E-5	&2.00	&3.99E-8	&2.99	&1.34E-5	&2.02\\ \hline
\multicolumn{13}{c}{$k = 1, \nu = 1\text{E}-2$}\\ \hline	
%k = 1, nu = 1e-2												
16	&1.04&	 	&9.17E-3&	 	&2.04E-2&		&2.37E-3&	 	&1.98E-5&	 	&9.32E-6&\\	
32	&2.63E-1	&1.99	&1.16E-3	&2.98	&5.09E-3	&2.00	&6.00E-4	&1.98	&2.52E-6	&2.98	&2.24E-6	&2.06\\
64	&6.59E-2	&1.99	&1.47E-4	&2.99	&1.27E-3	&2.00	&1.51E-4	&1.99	&3.18E-7	&2.99	&5.44E-7	&2.04\\
128	&1.65E-2	&2.00	&1.84E-5	&2.99	&3.18E-4	&2.00	&3.78E-5	&2.00	&3.99E-8	&2.99	&1.34E-7	&2.02\\ \hline
\multicolumn{13}{c}{$k = 1, \nu = 1\text{E}-4$}\\ \hline	
%k = 1,nu = 1e-4												
16	&- &-&-&-&-&-						&2.48E-3&	 	&2.20E-5&	 	&1.14E-6	&\\
32	&- &-&-&-&-&-						&6.03E-4	&2.04	&2.55E-6	&3.11	&1.32E-7	&3.11\\
64	&- &-&-&-&-&-						&1.51E-4	&2.00	&3.18E-7	&3.00	&1.70E-8	&2.97\\
128	&- &-&-&-&-&-						&3.78E-5	&2.00	&3.99E-8	&3.00	&2.13E-9	&2.98\\ \hline
\end{tabular}
\end{table}

\begin{table}
\caption{Example~\ref{Sect:Num-1}: Error Profiles and Convergence Results for $k=2$.}\label{Tab:Num-1_k2}
\centering
\tabcolsep=3pt
\begin{tabular}{c|cc|cc|cc||cc|cc|cc}
\hline\hline
 &\multicolumn{6}{c||}{WG Algorithm~\ref{alg:WG-2}} & \multicolumn{6}{c}{WG Algorithm~\ref{alg:WG}} \\
$1/h$	&$\3bar\be_h\3bar$	&Rate	&$\|\be_0\|$ &Rate		&$\|\epsilon_h\|$ &Rate	&$\3bar\be_h\3bar$	&Rate	&$\|\be_0\|$ &Rate		&$\|\epsilon_h\|$ &Rate	\\ \hline\hline
\multicolumn{13}{c}{$k = 2, \nu = 1$}\\ \hline	
%k = 2, nu = 1												
16	&8.45E-5&	 	&2.29E-7&	 	&4.06E-5&		&9.80E-5&	 	&4.96E-7&	 	&5.97E-5&\\	
32	&1.07E-5	&2.98	&1.49E-8	&3.95	&5.02E-6	&3.01	&1.24E-5	&2.98	&3.16E-8	&3.97	&7.40E-6	&3.01\\
64	&1.35E-6	&2.99	&9.49E-10&3.97	&6.22E-7	&3.01	&1.56E-6	&2.99	&1.99E-9	&3.99	&9.17E-7	&3.01\\
128	&1.69E-7	&3.00	&5.99E-11	&3.99	&7.73E-8	&3.01	&1.96E-7	&3.00	&1.25E-10	&4.00	&1.14E-7	&3.01\\ \hline
\multicolumn{13}{c}{$k = 2, \nu = 1\text{E}-2$}\\ \hline
%k = 2, nu = 1e-2												
16	&8.45E-5&	 	&2.30E-7&	 	&4.06E-7&		&9.80E-5&	 	&4.96E-7&	 	&5.97E-7&\\	
32	&1.07E-5	&2.98	&1.49E-8	&3.95	&5.02E-8	&3.01	&1.24E-5	&2.98	&3.16E-8	&3.97	&7.40E-8	&3.01\\
64	&1.35E-6	&2.99	&9.49E-10&3.97	&6.22E-9	&3.01	&1.56E-6	&2.99	&1.99E-9	&3.99	&9.17E-9	&3.01\\
128	&1.69E-7	&3.00	&5.99E-11	&3.99	&7.73E-10 &3.01	&1.96E-7	&3.00	&1.25E-10	&4.00	&1.14E-9	&3.01\\ \hline
\multicolumn{13}{c}{$k = 2, \nu = 1\text{E}-4$}\\ \hline
%k = 2, nu = 1e-4												
16	&9.39E-5	& 	&3.85E-7    &	    &8.83E-9&		&1.01E-4&	 	&5.14E-7&	 	&2.57E-8&\\	
32	&1.10E-5	&3.09&1.75E-8	  &4.46 &5.33E-10	&4.05	&1.25E-5	&3.01	&3.17E-8	&4.02	&1.66E-9	&3.95\\
64	&1.36E-6	&3.02&9.97E-10 &4.14 &6.29E-11	&3.08	&1.57E-6	&3.00	&1.99E-9	&3.99	&1.04E-10 &3.95\\
128	&6.13E-6	&-	 &6.07E-11 &4.04 &7.7562E-12	&3.02	&1.96E-7	&3.00	&1.25E-10	&4.00	&6.48E-12 &3.99\\ \hline\hline
\end{tabular}
\end{table}
}

\subsection{No Flow Test}\label{Sect:Num-2}
In this test, we shall compare the approximation for zero velocity by WG Algorithm~\ref{alg:WG} and WG Algorithm~\ref{alg:WG-2}. Let $\Omega=(0,1)^2$ and 
the velocity field is given by
\begin{eqnarray*}
\bu = (0,0)^\top.
\end{eqnarray*}
The pressure is defined as 
$$
p = -\frac{\text{Ra}}{2}y^2+\text{Ra}\; y-\frac{\text{Ra}}{3},
$$
where $\text{Ra} = 1000$.

\begin{figure}[H]
\centering
\begin{tabular}{ccc}
\includegraphics[width=.3\textwidth,height = .2\textwidth]{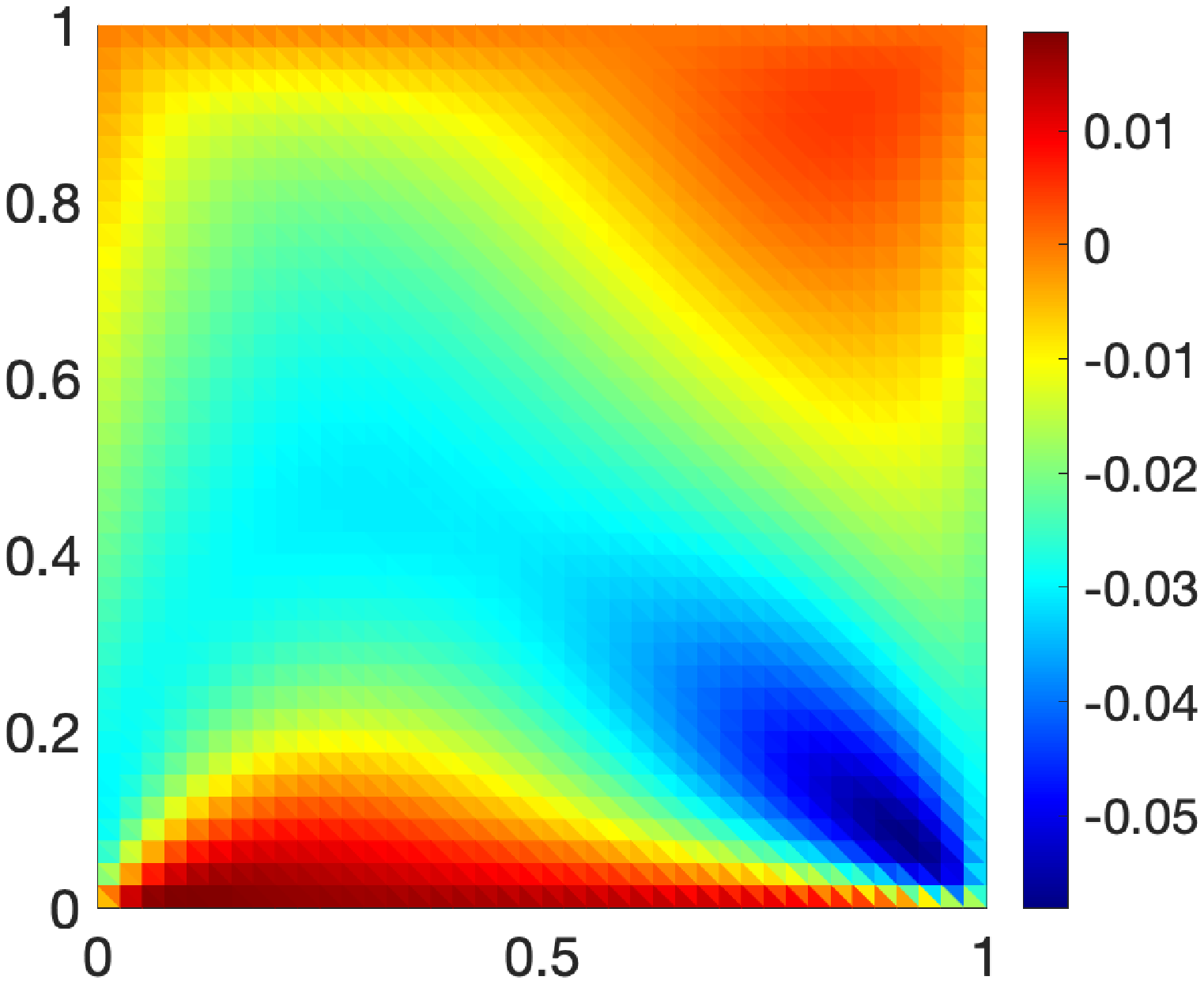}
&
\includegraphics[width=.3\textwidth,height = .2\textwidth]{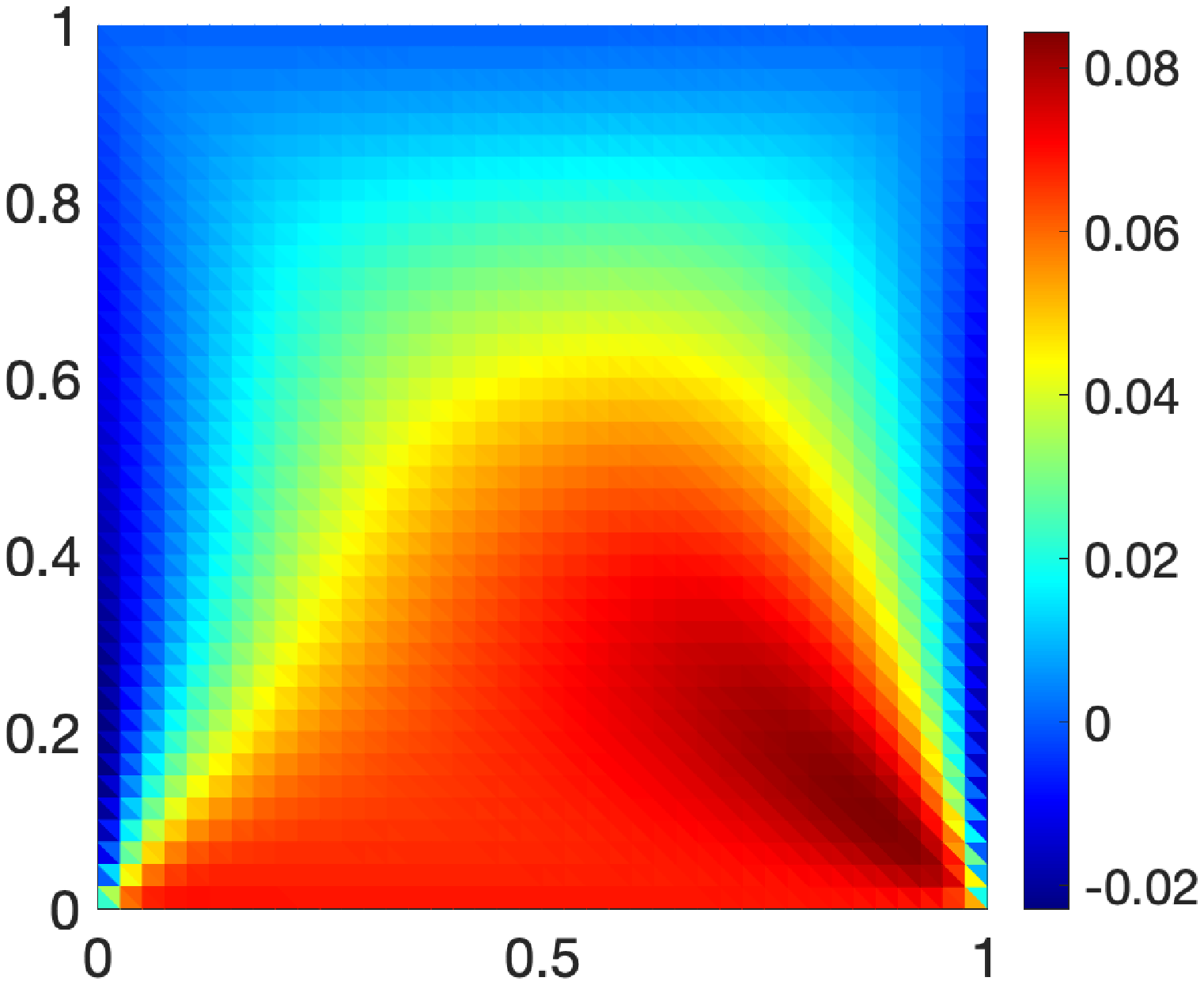}
&
\includegraphics[width=.3\textwidth,height = .2\textwidth]{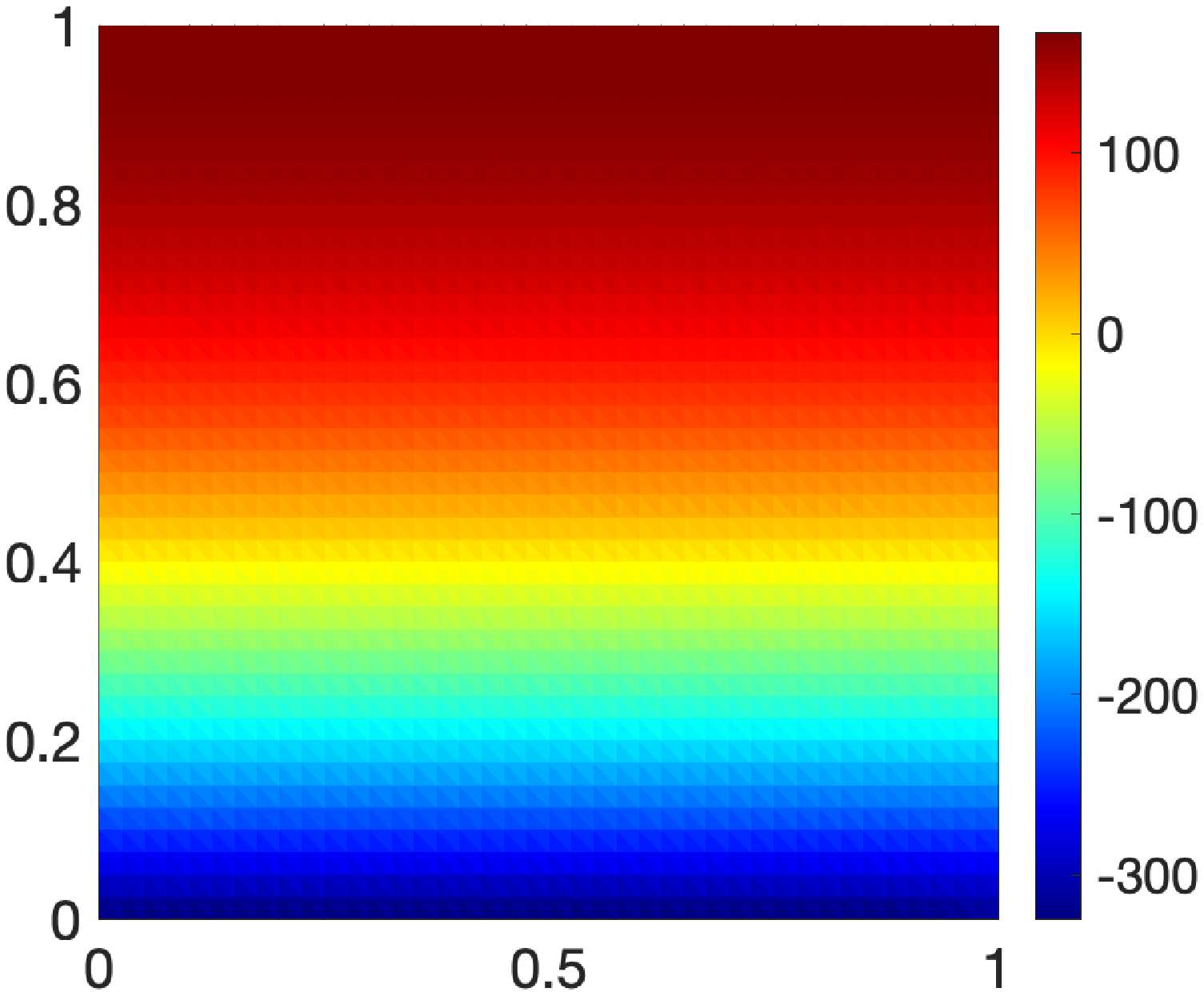}
\\
(a) & (b) & (c)\\
\includegraphics[width=.3\textwidth,height = .2\textwidth]{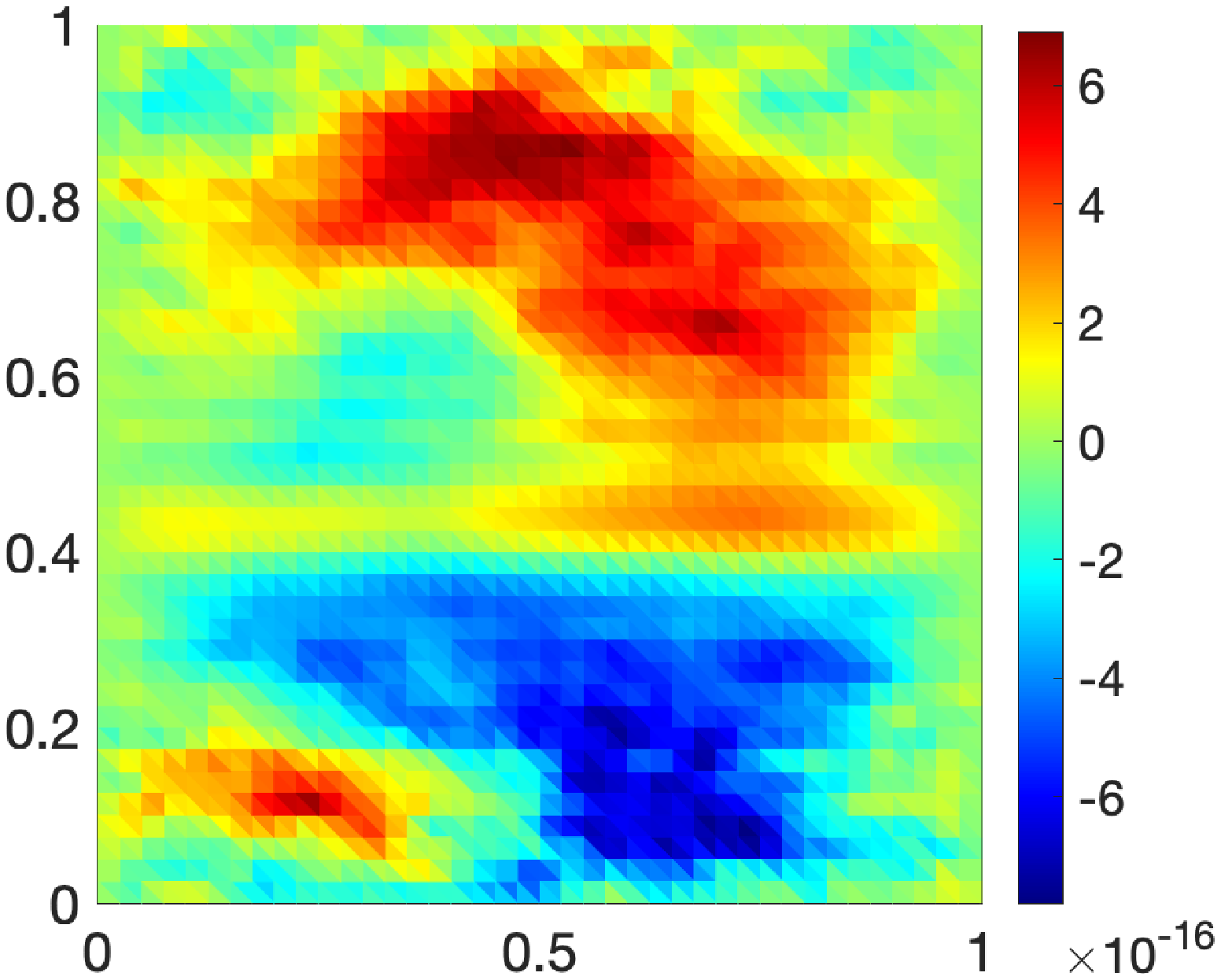}
&
\includegraphics[width=.3\textwidth,height = .2\textwidth]{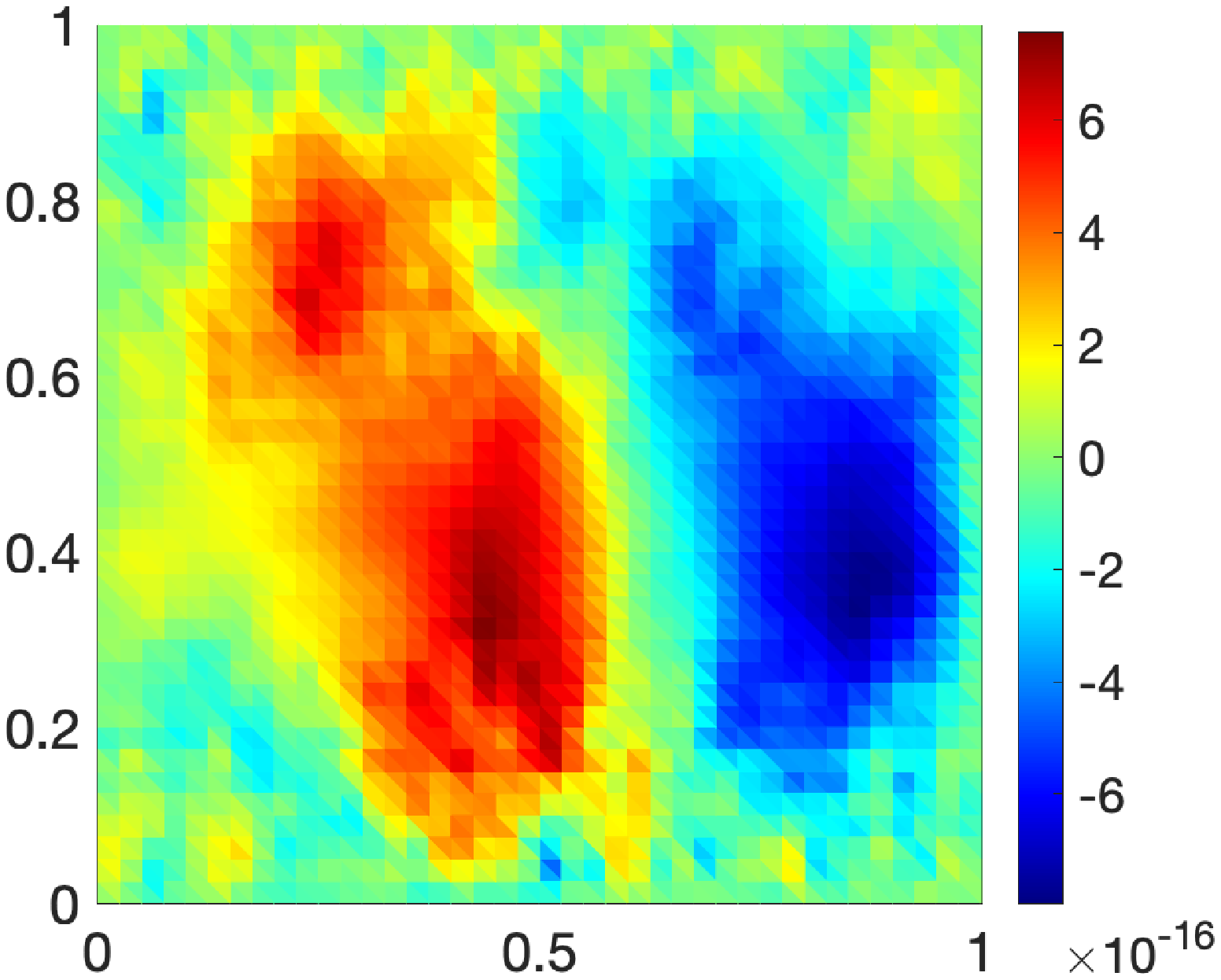}
&
\includegraphics[width=.3\textwidth,height = .2\textwidth]{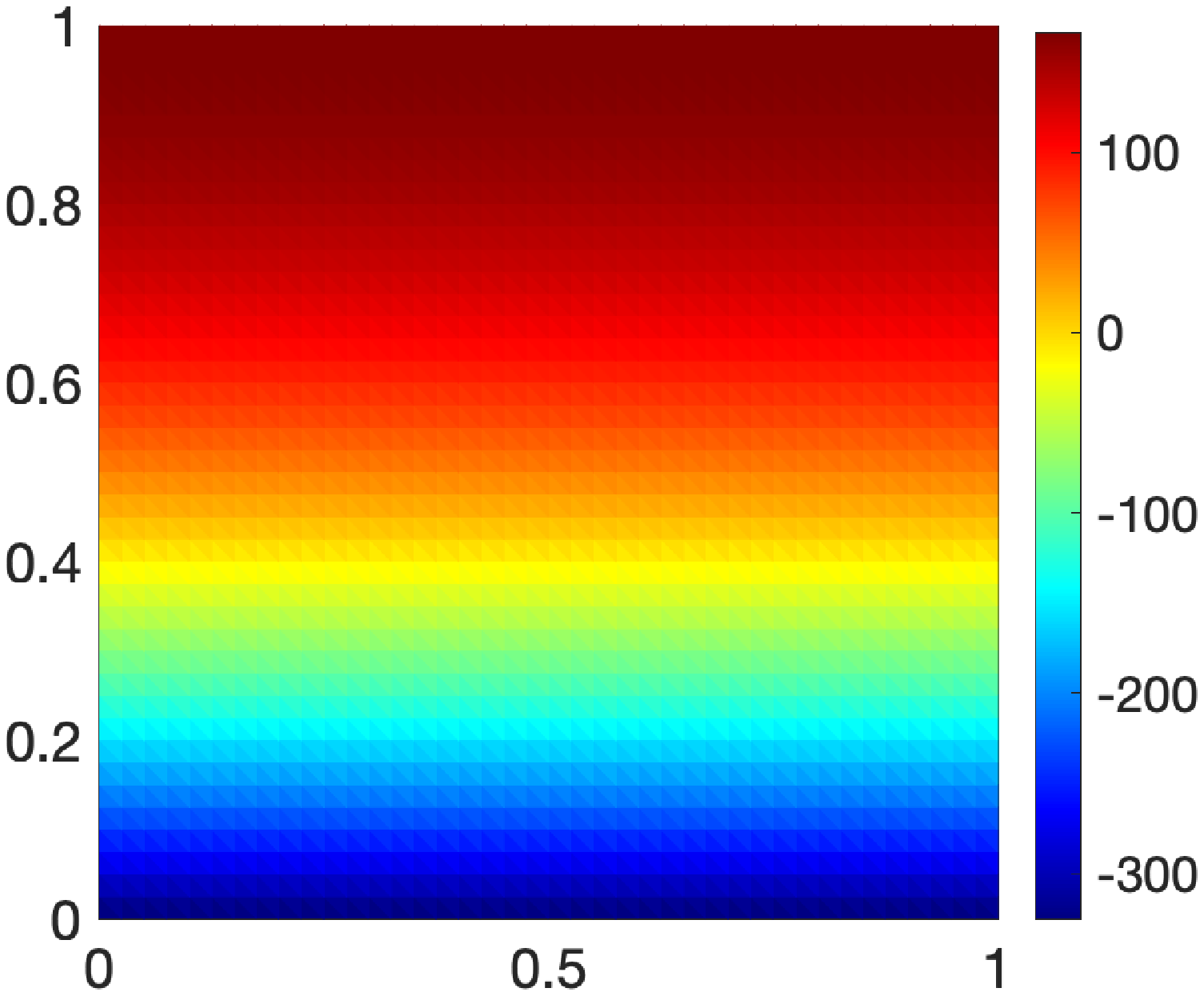}
\\
(d) & (e) &(f)
\end{tabular}
\caption{Example~\ref{Sect:Num-2}: Plots of WG approximation on mesh $h=1/40$ with $k=0$ for: (a)$\bu_1$ by Algorithm~\ref{alg:WG-2}; (a)$\bu_2$ by Algorithm~\ref{alg:WG-2}; (c)$p$ by Algorithm~\ref{alg:WG-2}; (d)$\bu_1$ by Algorithm~\ref{alg:WG}; (e)$\bu_2$ by Algorithm~\ref{alg:WG}; (c)$p$ by Algorithm~\ref{alg:WG};  }\label{Fig:Test-2}
\end{figure}

As the theoretical conclusions in Theorem~\ref{thm:main}, we expect numerical approximation gives 0 for velocity. However, as the dependence of pressure for the error estimate in (\ref{eq:error-u-wg2}), we cannot expect Algorithm~\ref{alg:WG-2} produces very well simulation for polynomial degree less than $2$. We perform the WG simulation on the mesh with size $h=1/40$ and $k=0$ and the numerical solutions are presented in Fig.~\ref{Fig:Test-2}. The pattern of numerical solutions validate the theoretical conclusions.

\subsection{L-shape Benchmark}\label{Sect:Num-3}
\begin{figure}[H]
\centering
\includegraphics[width=0.45\textwidth]{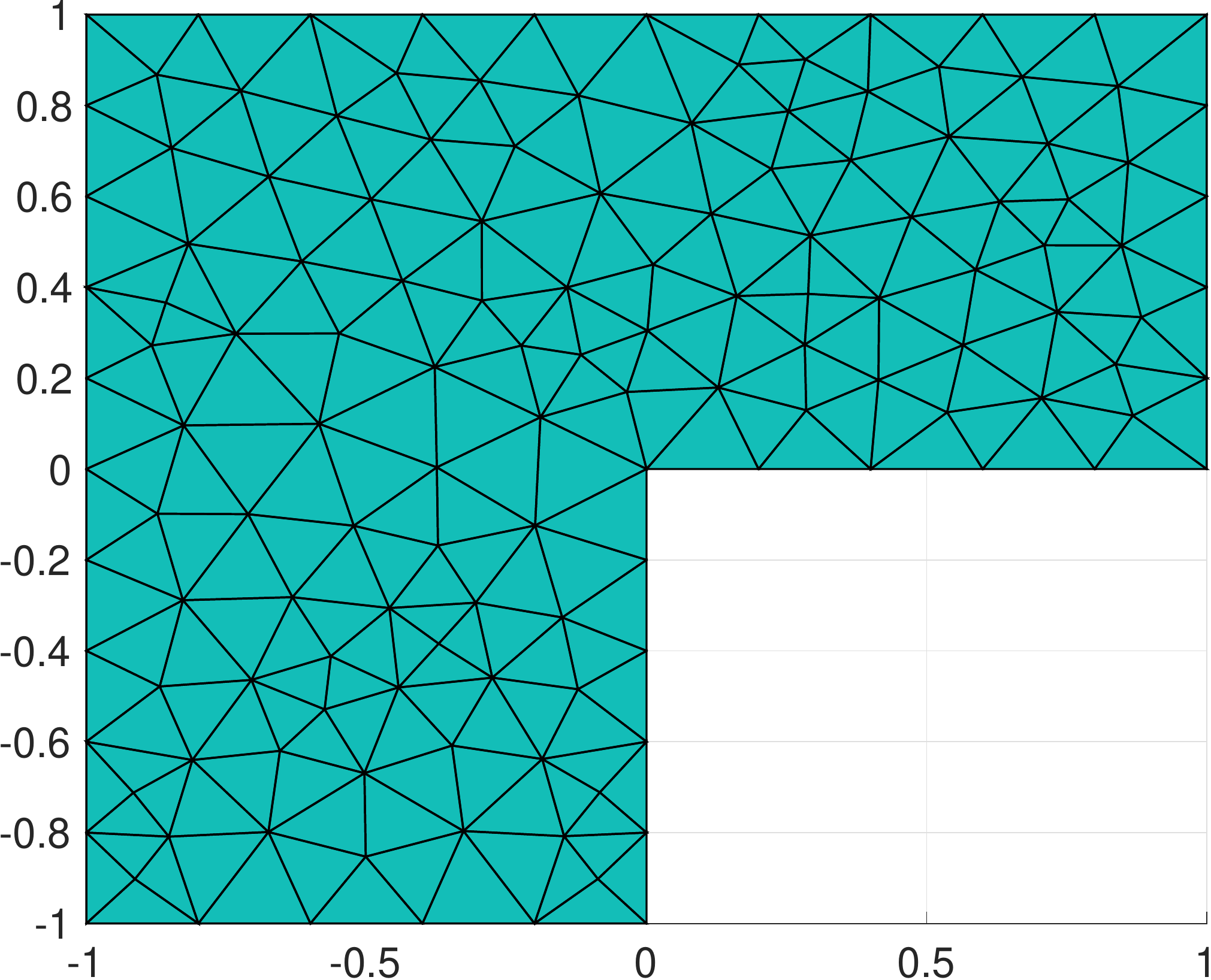}
\caption{Example~\ref{Sect:Num-3}: Illustration of intial mesh.}\label{fig:Lshape-mesh1}
\end{figure}
In this test, let the domain $\Omega = (-1,1)^2\backslash [0,1]\times[-1,0]$ and the exact solutions are chosen as follows:
\begin{eqnarray*}
\bu = \begin{pmatrix}
\sin(\pi x)\sin(\pi y),\\
\cos(\pi x)\cos(\pi y)
\end{pmatrix},\ p = r^{2/3}\sin(\frac{2\theta}{3}),
\end{eqnarray*}
where $r,\theta$ are in the polar coordinates. As we know, the velocity is smooth and the regularity of pressure is approximately $H^{1.67}$. We shall show the advantages of Algorithm~\ref{alg:WG} to Algorithm~\ref{alg:WG-2} by performing simulations with varying degrees in weak Galerkin finite elements. Let $\nu = 1$, we shall perform Algorithm~\ref{alg:WG} and Algorithm~\ref{alg:WG-2} to compare their corresponding numerical performance. The coarsest mesh is shown in Fig.~\ref{fig:Lshape-mesh1} and then the next level of mesh is derived by uniform refining the previous level of the mesh.

\begin{figure}[H]
\centering
\begin{tabular}{cc}
\includegraphics[width=.46\textwidth]{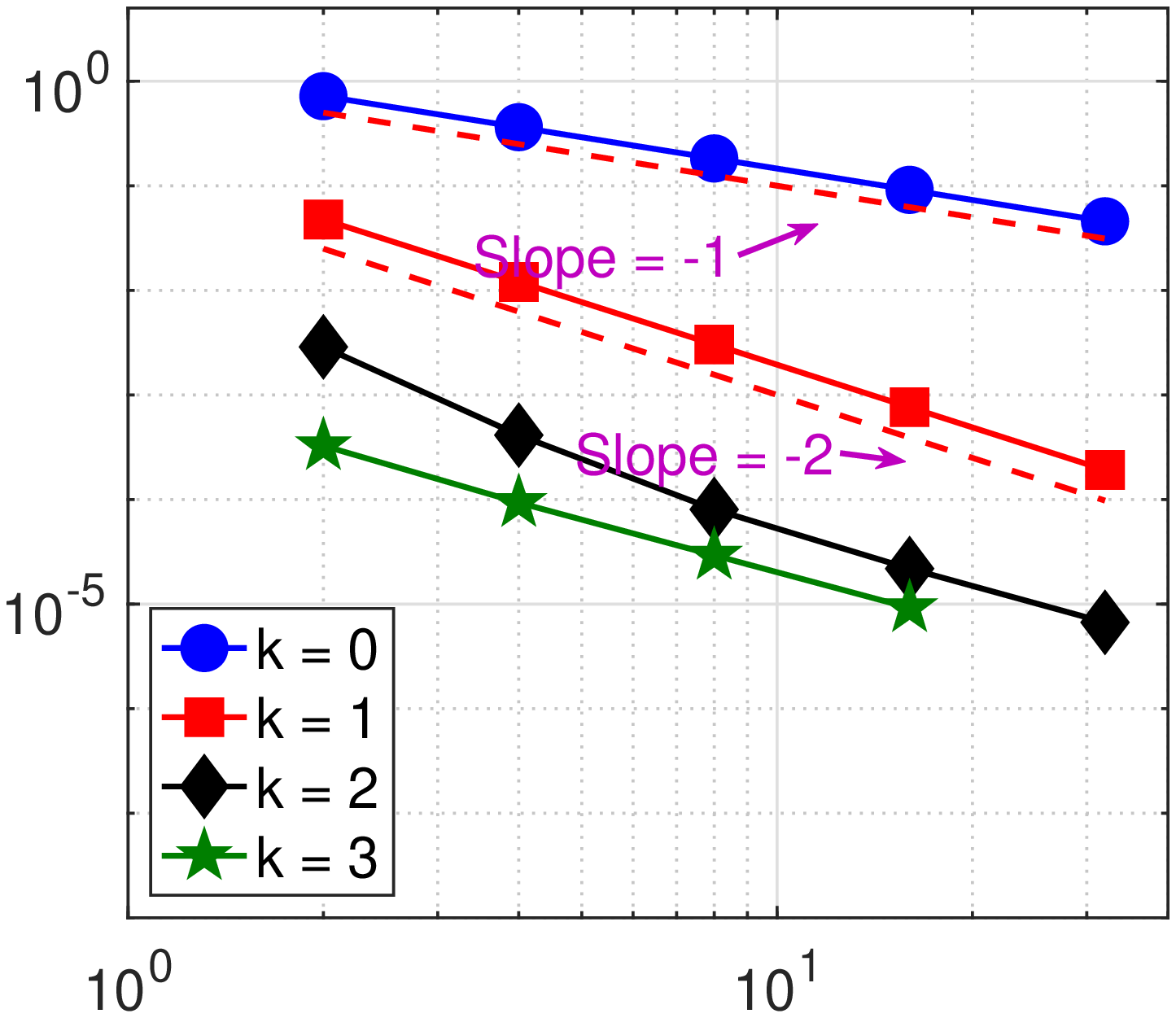}
&
\includegraphics[width=.46\textwidth]{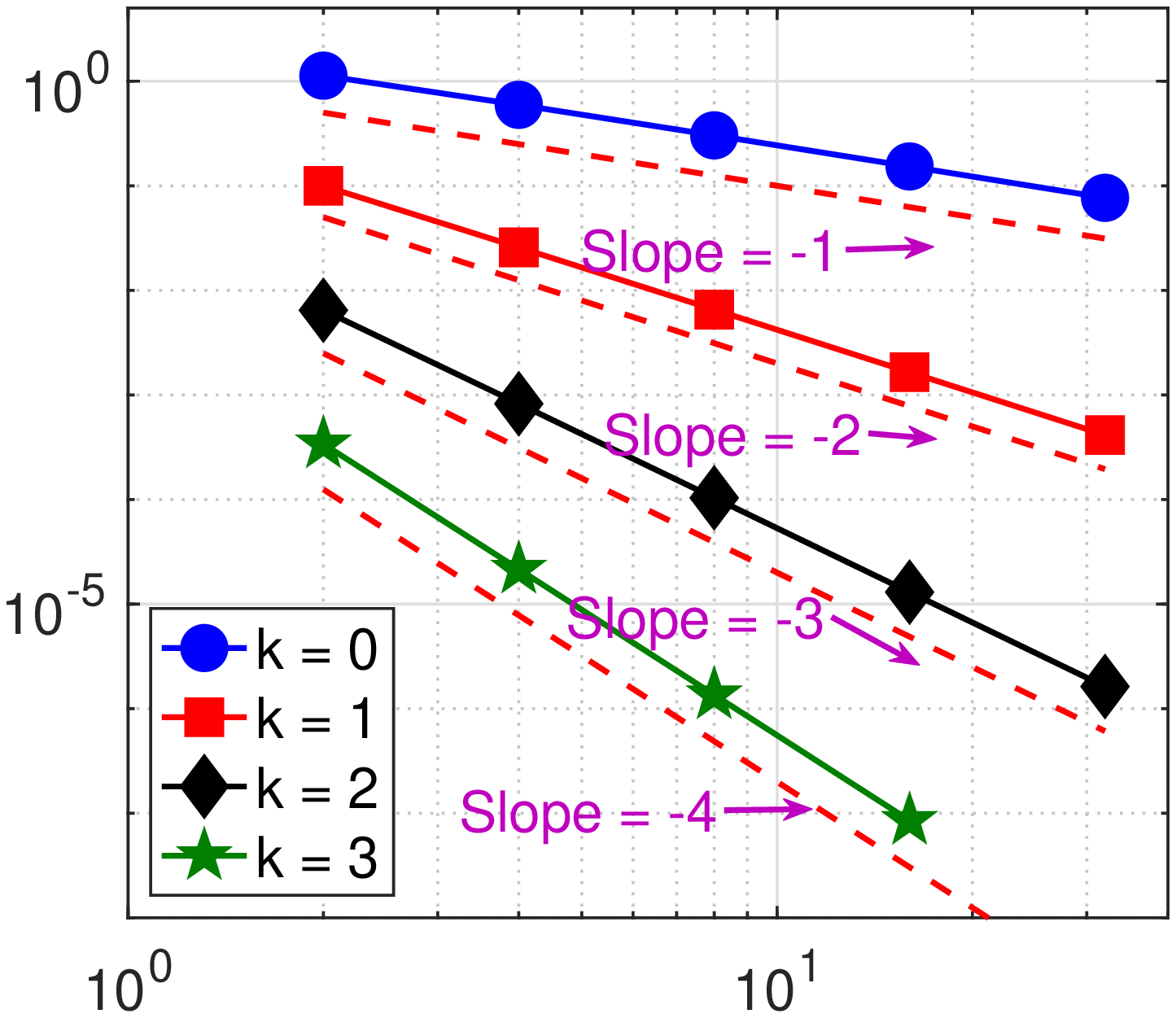}
\\
(a) & (b) 
\end{tabular}
\caption{Example~\ref{Sect:Num-3}: Plots of convergence results for $\3bar\be_h\3bar$: (a) WG Algorithm~\ref{alg:WG-2}; (b) WG Algorithm~\ref{alg:WG}.  }\label{Fig:Test-Lshape_uH1}
\end{figure}
\begin{figure}[H]
\centering
\begin{tabular}{cc}
\includegraphics[width=.46\textwidth]{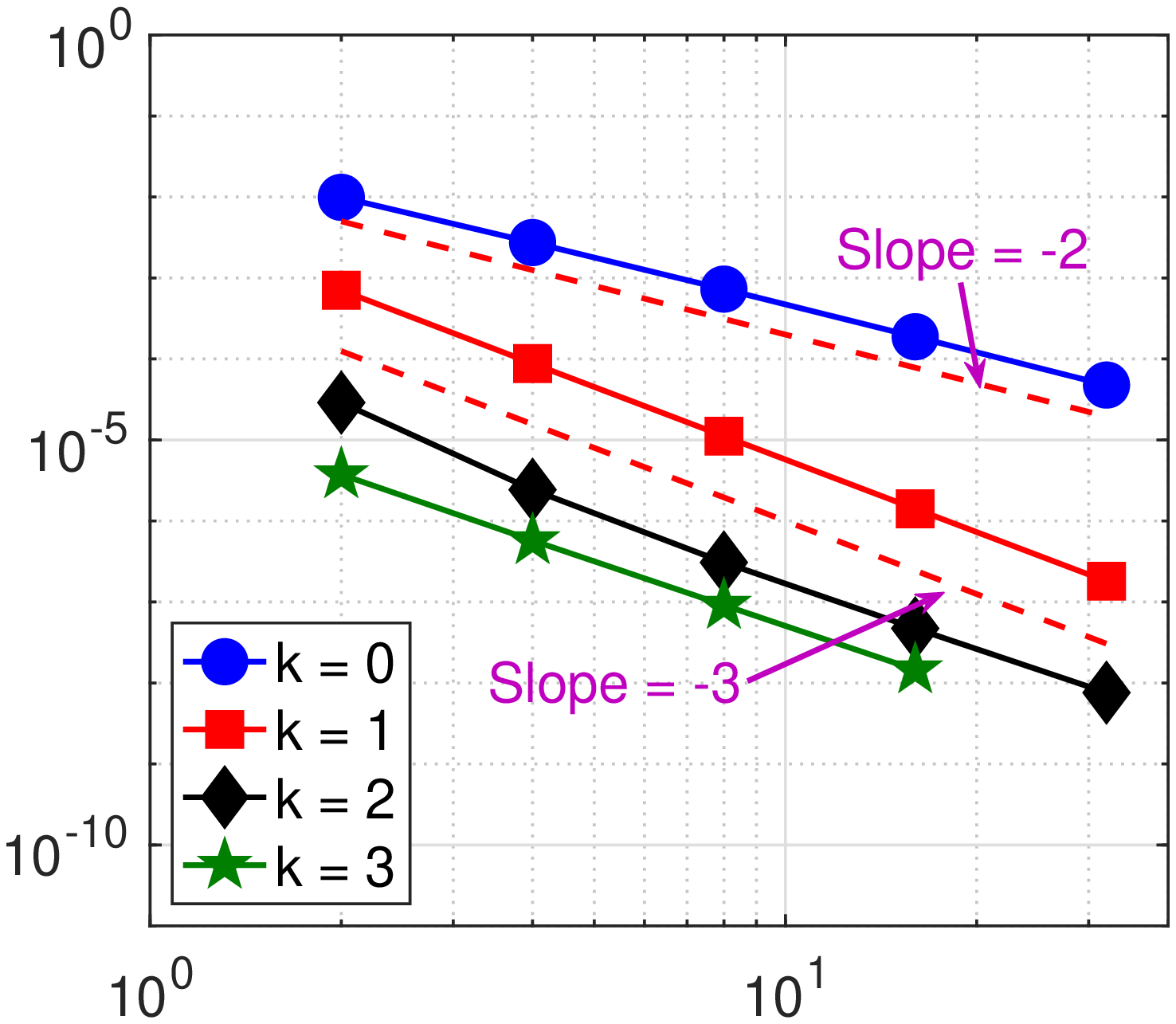}
&
\includegraphics[width=.46\textwidth]{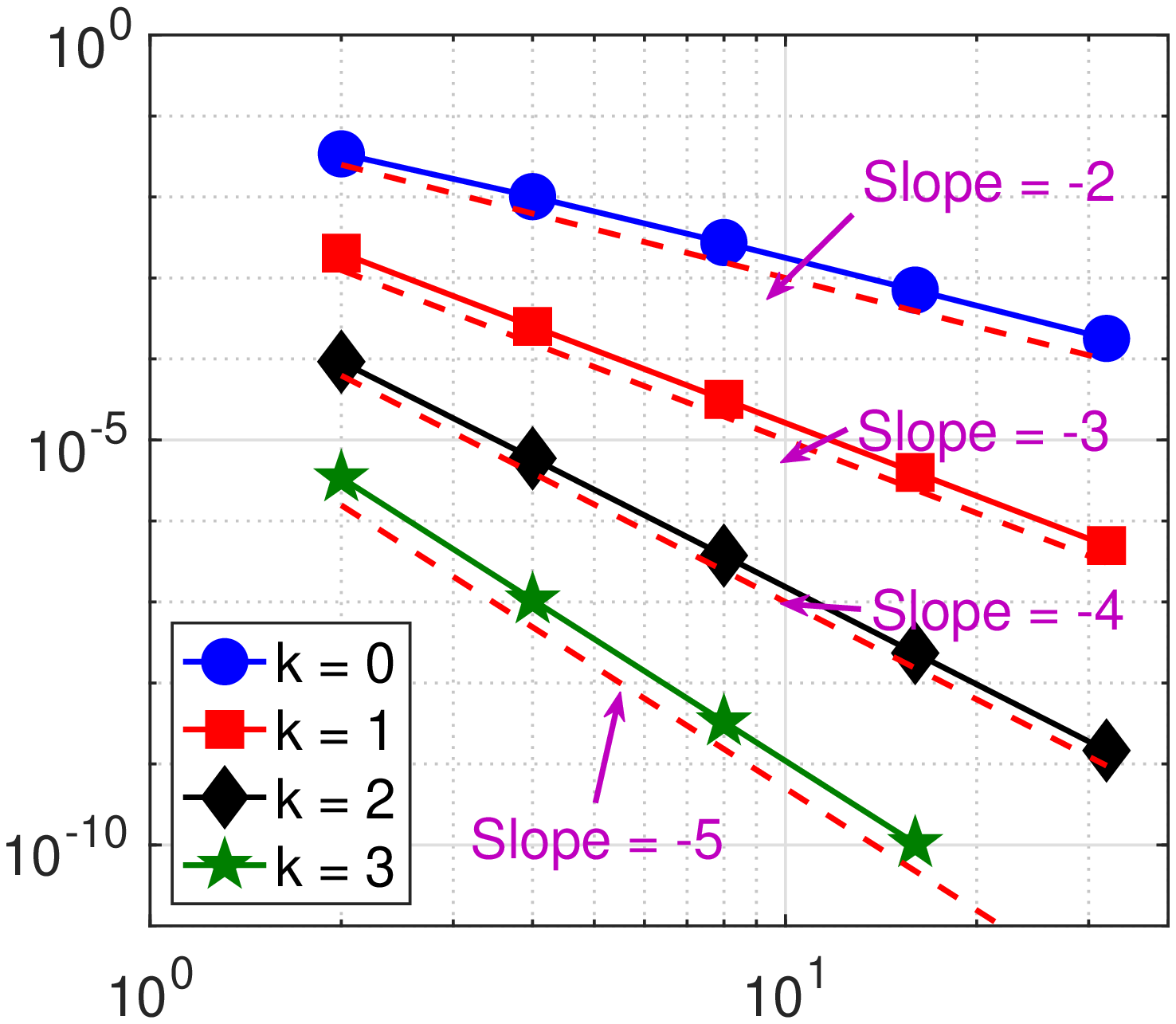}
\\
(a) & (b) 
\end{tabular}
\caption{Example~\ref{Sect:Num-3}: Plots of convergence results for $\|\be_h\|$: (a) WG Algorithm~\ref{alg:WG-2}; (b) WG Algorithm~\ref{alg:WG}.  }\label{Fig:Test-Lshape_uL2}
\end{figure}
\begin{figure}[H]
\centering
\begin{tabular}{cc}
\includegraphics[width=.46\textwidth]{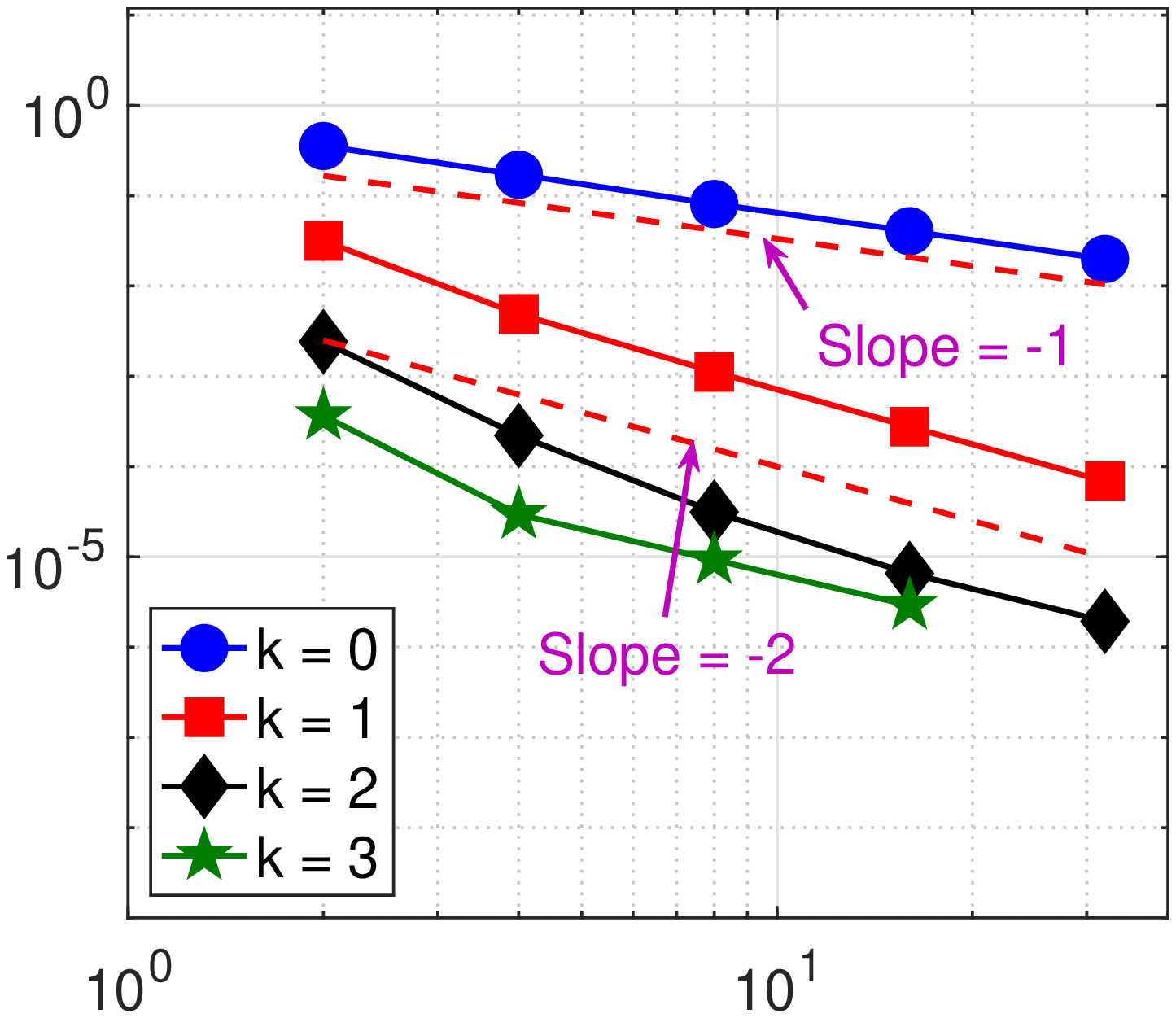}
&
\includegraphics[width=.46\textwidth]{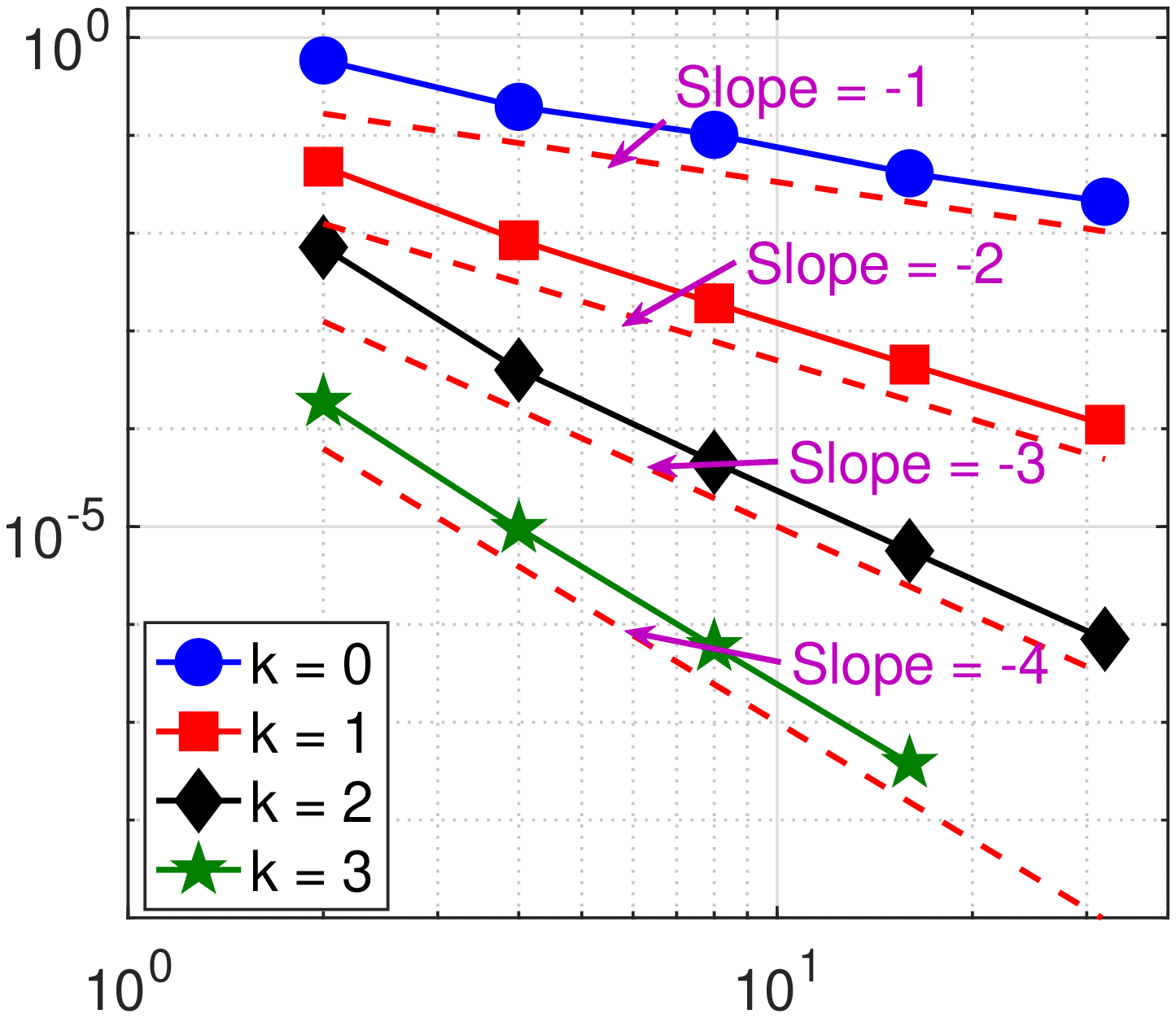}
\\
(a) & (b) 
\end{tabular}
\caption{Example~\ref{Sect:Num-3}: Plots of convergence results for $\|\epsilon_h\|$: (a) WG Algorithm~\ref{alg:WG-2}; (b) WG Algorithm~\ref{alg:WG}.  }\label{Fig:Test-Lshape_pL2}
\end{figure}

The error profiles and convergence results are plotted in Fig.~\ref{Fig:Test-Lshape_uH1}-Fig.~\ref{Fig:Test-Lshape_pL2}. As one can see from the left columns in these figures that the convergence rate of Algorithm~\ref{alg:WG-2} for $\3bar\be_h\3bar$ is limited by $\mathcal{O}(h^2)$, $\|\be_h\|$ is limited by $\mathcal{O}(h^3)$, and $\|\epsilon_h\|$ is limited by $\mathcal{O}(h^2)$ even for high polynomial degrees. This is because, the regularity of pressure will affect the simulation properties in Algorithm~\ref{alg:WG-2}. However, for all the simulations carried out by Algorithm~\ref{alg:WG}, we can achieve the optimal rate in convergence, which is $\mathcal{O}(h^{k+1})$, $\mathcal{O}(h^{k+2})$, and $\mathcal{O}(h^{k+1})$ for the errors measured in $\3bar\be_h\3bar$, $\|\be_h\|$, and $\|\epsilon_h\|$, respectively. This observation again validate the robustness of the proposed numerical scheme.

%\begin{figure}[H]
%\centering
%\begin{tabular}{ccc}
%\includegraphics[width=.32\textwidth]{./figure/Lshape_Type1_H1Err}
%&
%\includegraphics[width=.32\textwidth]{./figure/Lshape_Type1_uL2Err}
%&
%\includegraphics[width=.32\textwidth]{./figure/Lshape_Type1_pL2Err}
%\\
%(a) & (b) & (c)\\
%\includegraphics[width=.32\textwidth]{./figure/Lshape_Type2_H1Err}
%&
%\includegraphics[width=.32\textwidth]{./figure/Lshape_Type2_uL2Err}
%&
%\includegraphics[width=.32\textwidth]{./figure/Lshape_Type2_pL2Err}
%\\
%(d) & (e) &(f)
%\end{tabular}
%\caption{Example~\ref{Sect:Num-3}: Plots of WG approximation properties on mesh $h=1/40$ with: (a)$\bu_1$ by Algorithm 2; (a)$\bu_2$ by Algorithm 2; (c)$p$ by Algorithm 2; (d)$\bu_1$ by Algorithm 1; (e)$\bu_2$ by Algorithm 1; (c)$p$ by Algorithm 1;.  }%\label{Fig:Test-LidCavity}
%\end{figure}

%\subsection{Test 2: Vanishing external load}\label{Num-2}
%In this test, we consider the disk $\Omega_D$ where we compare the results obtained with discretization. The solutions are chosen in such a way that the pressures balance the nonlinear convective term yielding a vanishing external load $\bbf = 0$. 
%\subsubsection{Test A}
%We take $\nu = 1$ and the exact solution
%\begin{eqnarray*}
%\bu = \begin{pmatrix}
%-y\\x
%\end{pmatrix},\ p = -\frac{x^2+y^2}{2}+\frac{1}{4}.
%\end{eqnarray*}
%\subsubsection{Test B}
%We take $\nu = 1$ and the exact solution
%\begin{eqnarray*}
%\bu = 3\begin{pmatrix}
%x^2-y^2\\-2xy
%\end{pmatrix},\ p = 9\frac{(x^2+y^2)^2}{2}-\frac{3}{2}.
%\end{eqnarray*}

\subsection{Kovasznay Flow Benchmark}\label{Sect:Num-4}
Let $\Omega=(-0.5,1.5)\times(0,2)$ 
%with known analytical solution to assess the convergence of the numerical algorithm. This
and the exact solution is chosen as the Kovasznay flow solution \cite{Kovasznay1948}.
\begin{eqnarray*}
\bu_1(x,y) = 1-\exp(\lambda x)\cos(2\pi y),\\
\bu_2(x,y) = \frac{\lambda}{2\pi}\exp(\lambda x)\sin(2\pi y)\\
p(x,y) = -\frac{1}{2}\exp(2\lambda x)+\frac{\lambda}{2}(\exp(4\lambda)-1)
\end{eqnarray*}
with $\lambda = \text{Re}/2-\sqrt{\text{Re}^2/4+4\pi^2}$ and $\text{Re} = 1/(2\nu).$

\begin{figure}[H]
\centering
\begin{tabular}{cc}
\includegraphics[width=.45\textwidth]{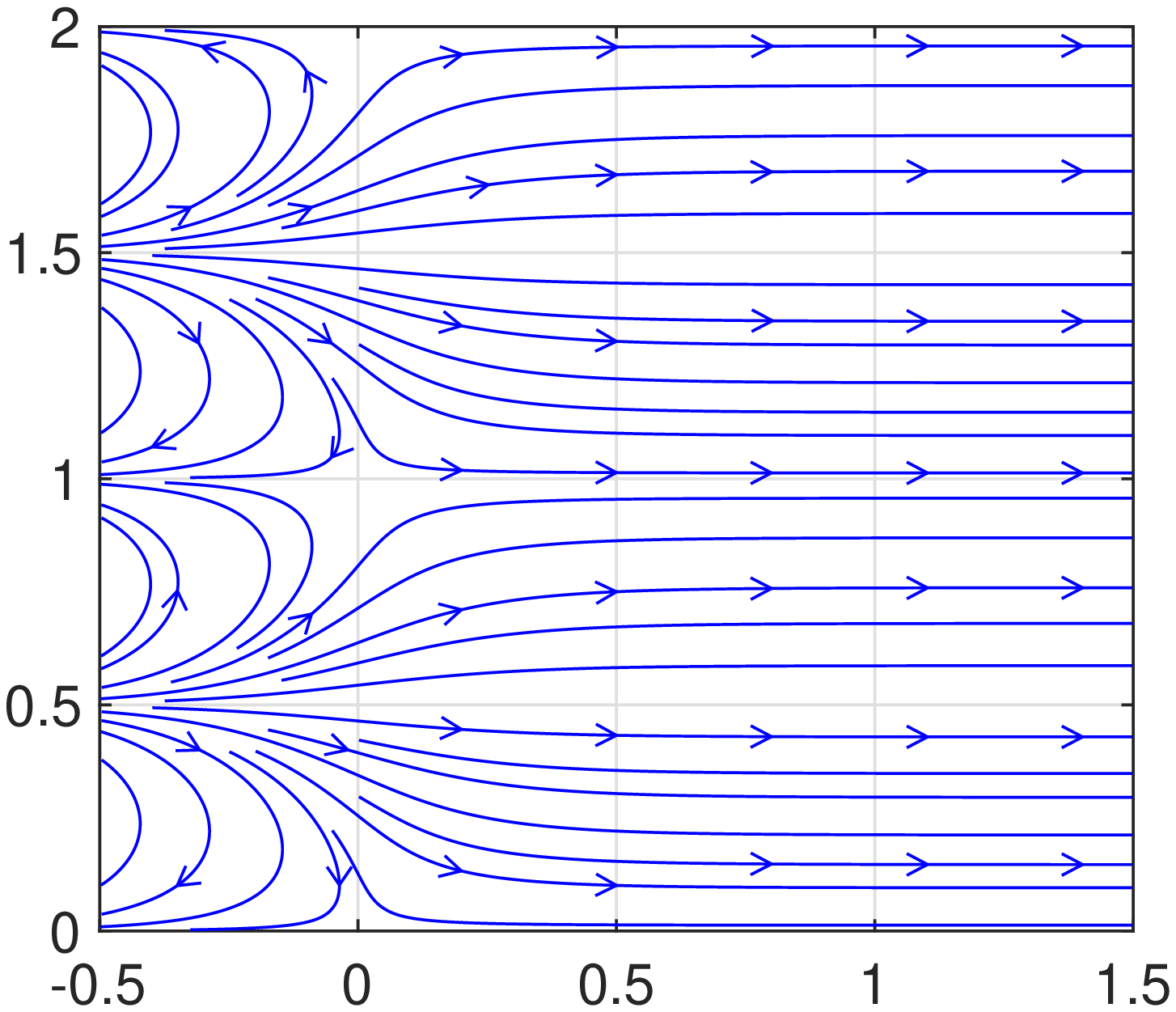}
&
\includegraphics[width=.45\textwidth]{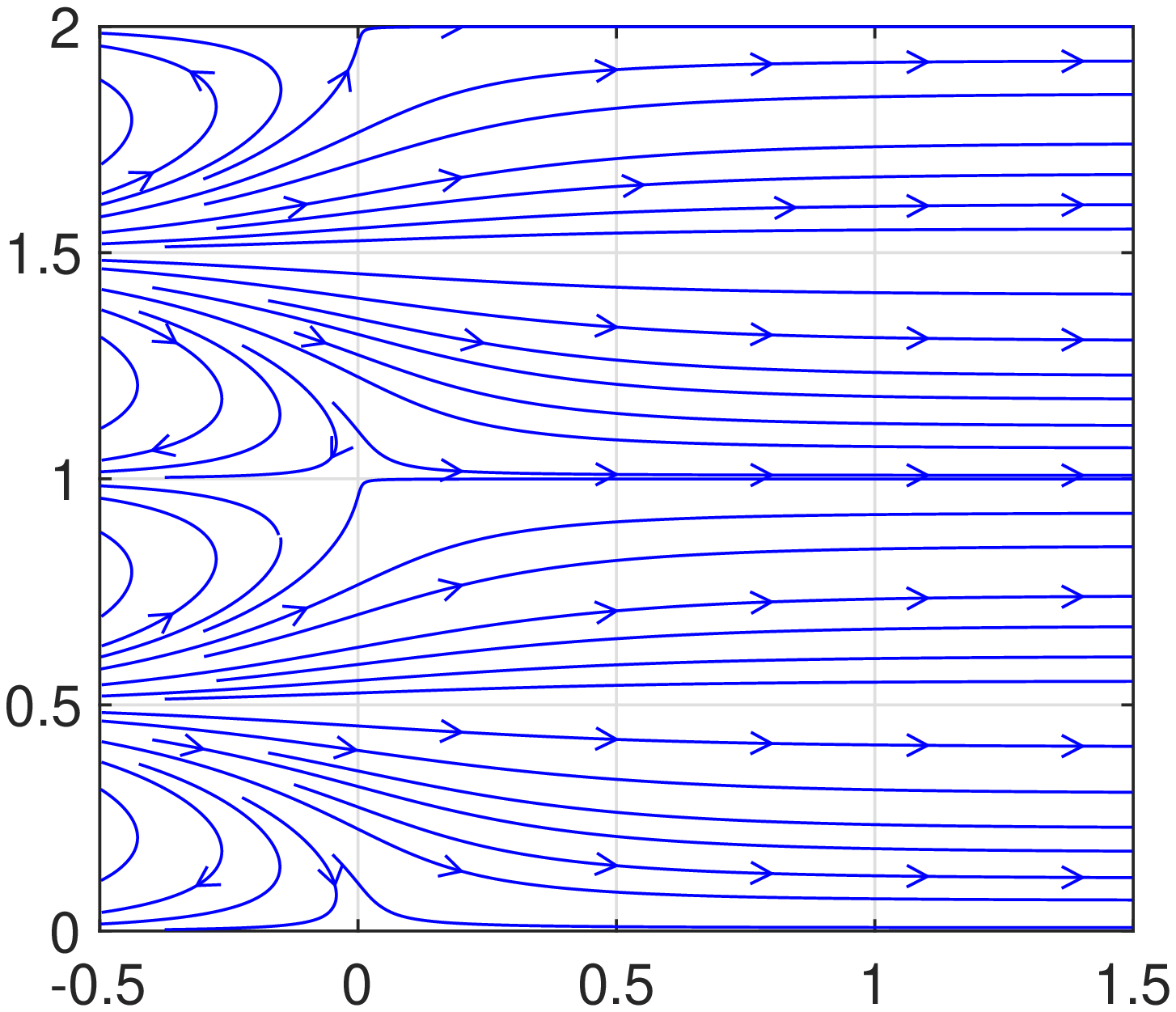}
\\
(a) & (b)
\end{tabular}
\caption{Example~\ref{Sect:Num-4}: Streamline plots of $\bu_h$ corresponding to: (a) $\nu = 1$; (b) $\nu = 1\text{E}$-1.  }\label{Fig:Test4-1}
\end{figure}

\begin{figure}[H]
\centering
\begin{tabular}{cc}
\includegraphics[width=.45\textwidth]{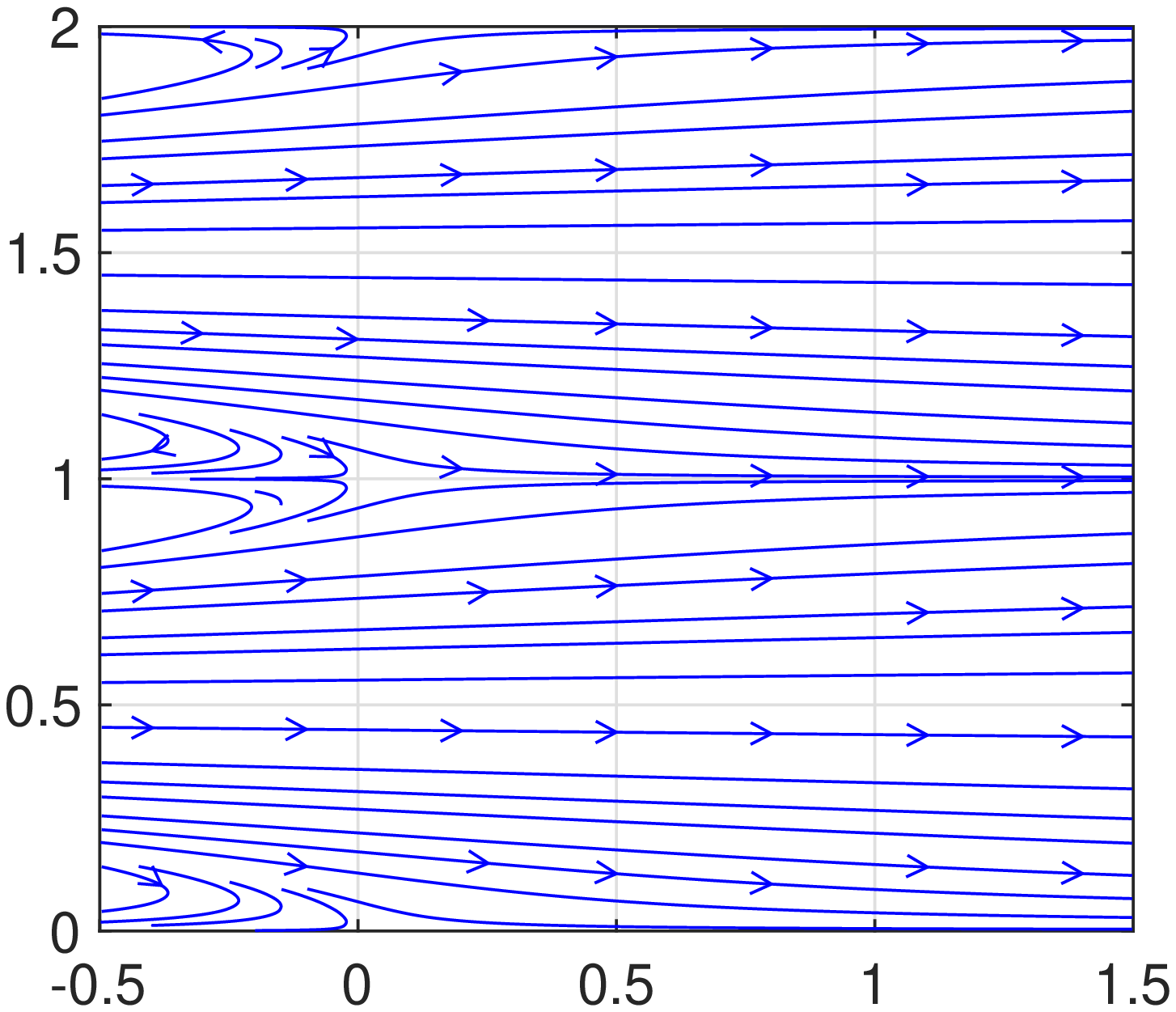}
&
\includegraphics[width=.45\textwidth]{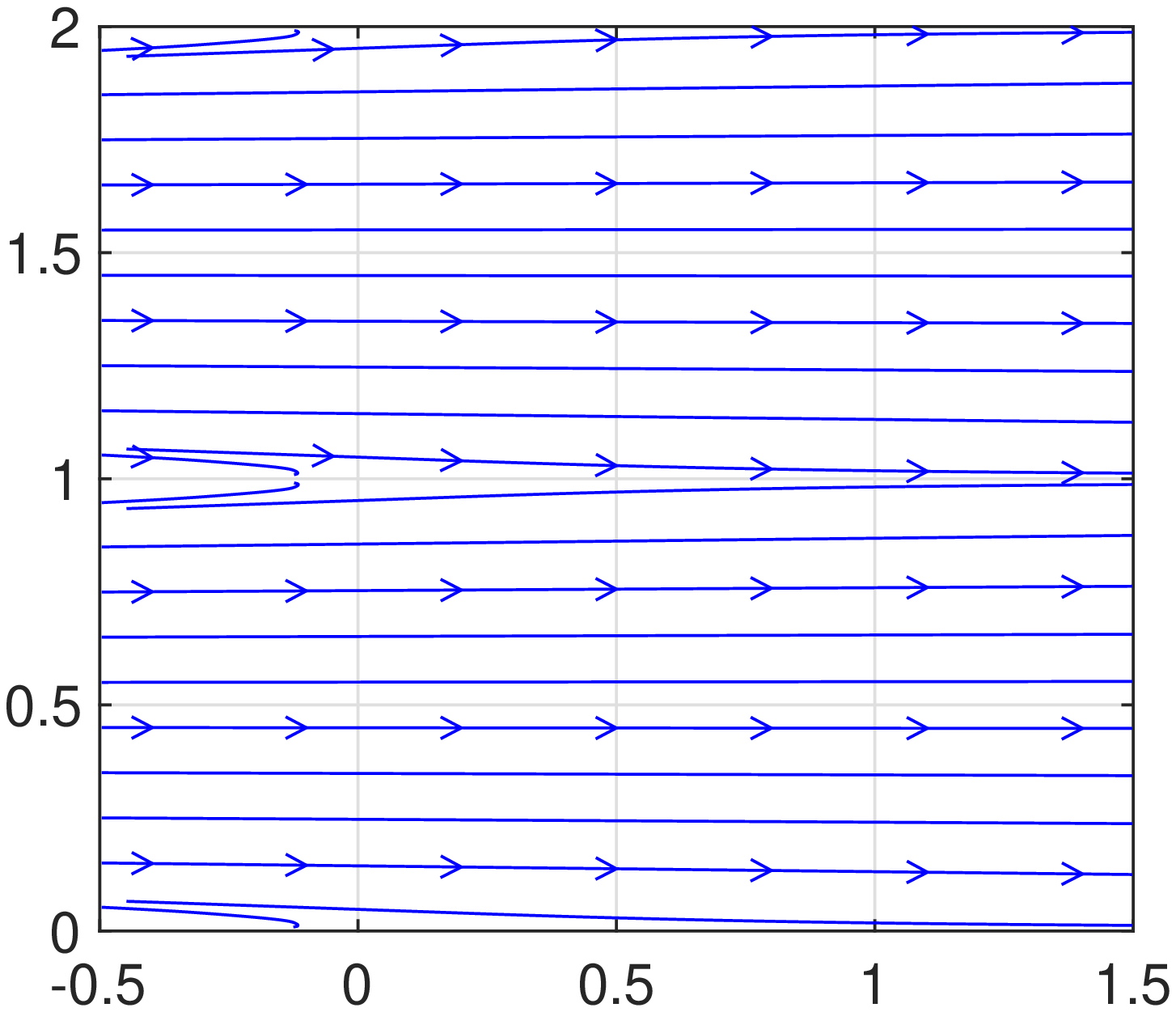}
\\
(a) & (b)
\end{tabular}
\caption{Example~\ref{Sect:Num-4}: Streamline plots of $\bu_h$ corresponding to: (a) $\nu = 1\text{E-}2$; (b) $\nu = 1\text{E}$-3.  }\label{Fig:Test4-2}
\end{figure}

We shall perform Algorithm~\ref{alg:WG} for various viscosity values, including $\nu = 1, $1E-1, 1E-2, and 1E-3. The corresponding numerical solutions in velocity are plotted in Fig.~\ref{Fig:Test4-1}-Fig.~\ref{Fig:Test4-2}. All the results match with the exact solutions very well.

%\begin{figure}[H]
%\centering
%\begin{tabular}{cc}
%\includegraphics[width=.45\textwidth]{./figure/ns8_nu1e-1}
%&
%\includegraphics[width=.45\textwidth]{./figure/ns8_nu25e-2}
%\\
%(a) & (b)
%\end{tabular}
%\caption{Example~\ref{Sect:Num-4}: Plots of numerical solution corresponding to: (a) $\nu = 1\text{E-}1$; (b) $\nu = 2.5\text{E}-2$.  }%\label{Fig:Test-LidCavity}
%\end{figure}

\subsection{Robustness of Irrotational body forces}\label{Sect:Num-5}
In this test, we shall demonstrate the robustness of the proposed method for large irrotational body forces. Let $\Omega = (0,1)^2$ and exact solutions are given by
\begin{eqnarray*}
\bu=\begin{pmatrix}
-y\\x
\end{pmatrix},\quad p=\lambda x^3+\frac{x^2+y^2}{2}-\frac{1}{4}.
\end{eqnarray*} 
It is easy to check that the force is purely irrotational, as computed below,
\begin{eqnarray*}
\bbf=\begin{pmatrix}
3\lambda x^2\\ 0
\end{pmatrix}.
\end{eqnarray*}
In the experiment, we take $\lambda = 10$ and $\lambda =1\text{E+}6$ and consider polynomial degrees $k = 0, 1,2$ to carry out Algorithm~\ref{alg:WG} and Algorithm~\ref{alg:WG-2}.

\begin{figure}[H]
\centering
\begin{tabular}{cc}
\includegraphics[width=.45\textwidth]{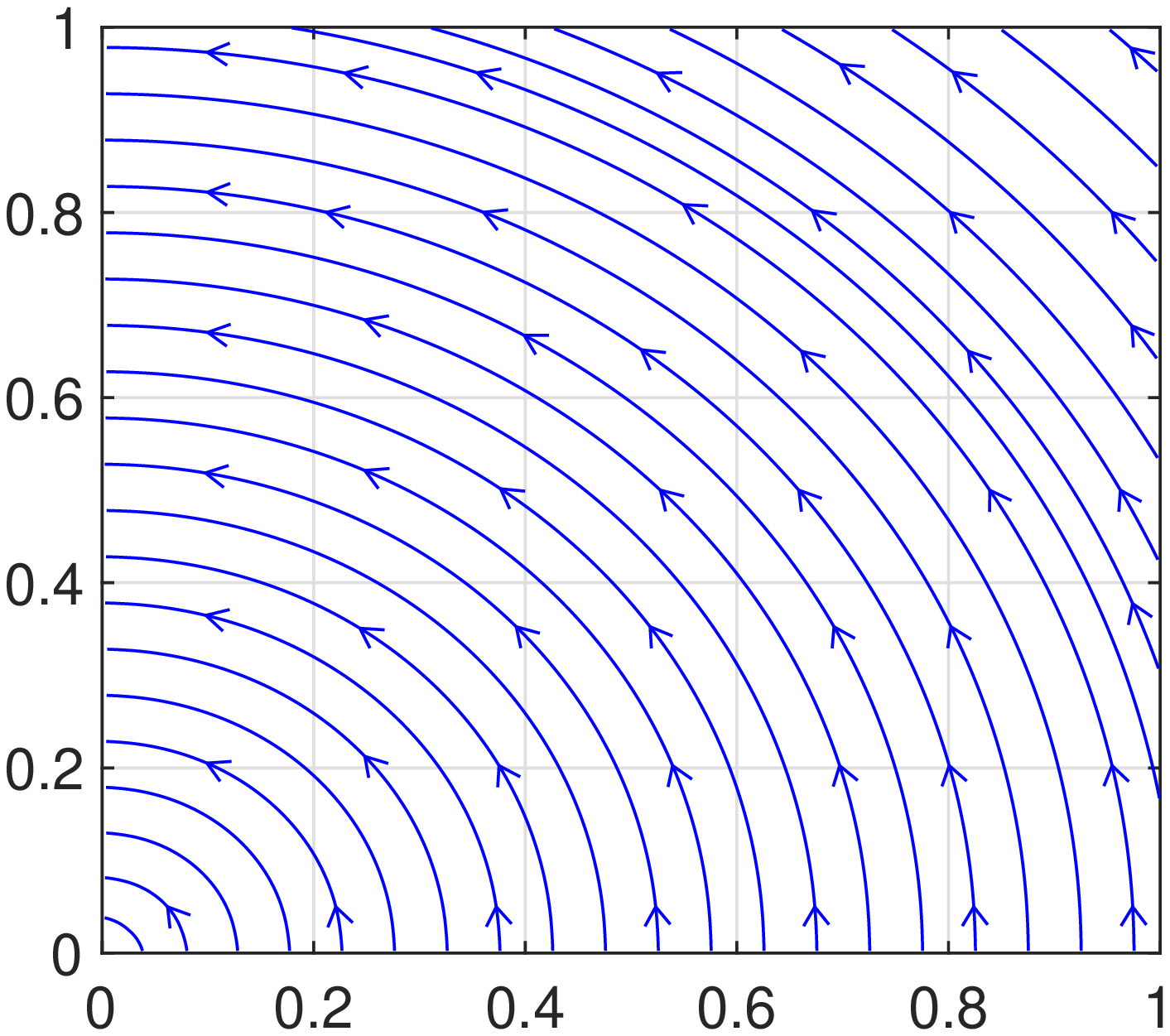}
&
\includegraphics[width=.45\textwidth]{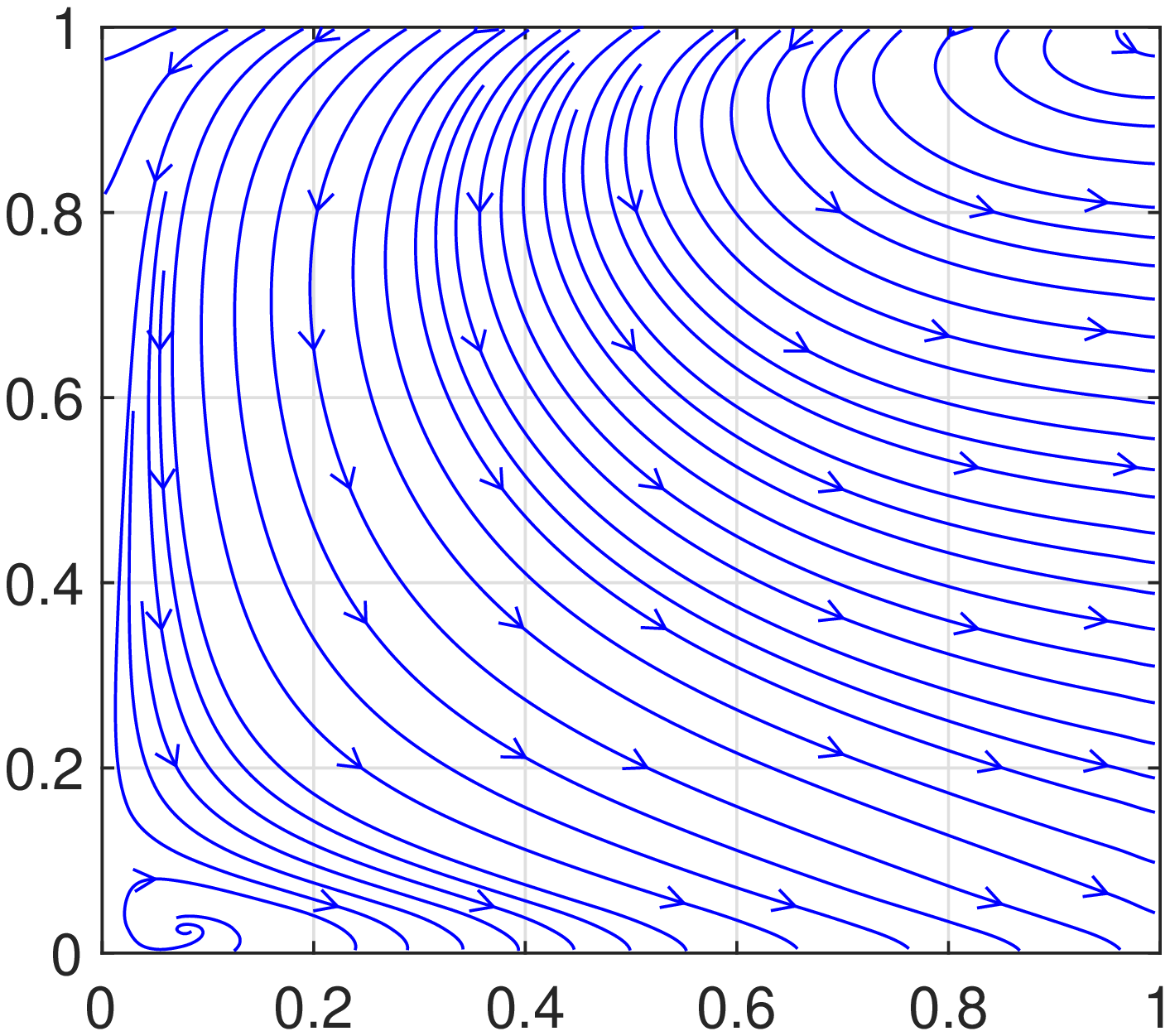}
\\
(a) & (b)
\end{tabular}
\caption{Example~\ref{Sect:Num-5}: Streamline Plots of WG Algorithm~\ref{alg:WG-2} on mesh $h=1/40$ and $k = 0$ with (a) $\lambda = 10$; (b) $\lambda = 1\text{E+}6$.  }\label{Fig:Test5-1}
\end{figure}

\begin{figure}[H]
\centering
\begin{tabular}{cc}
\includegraphics[width=.45\textwidth]{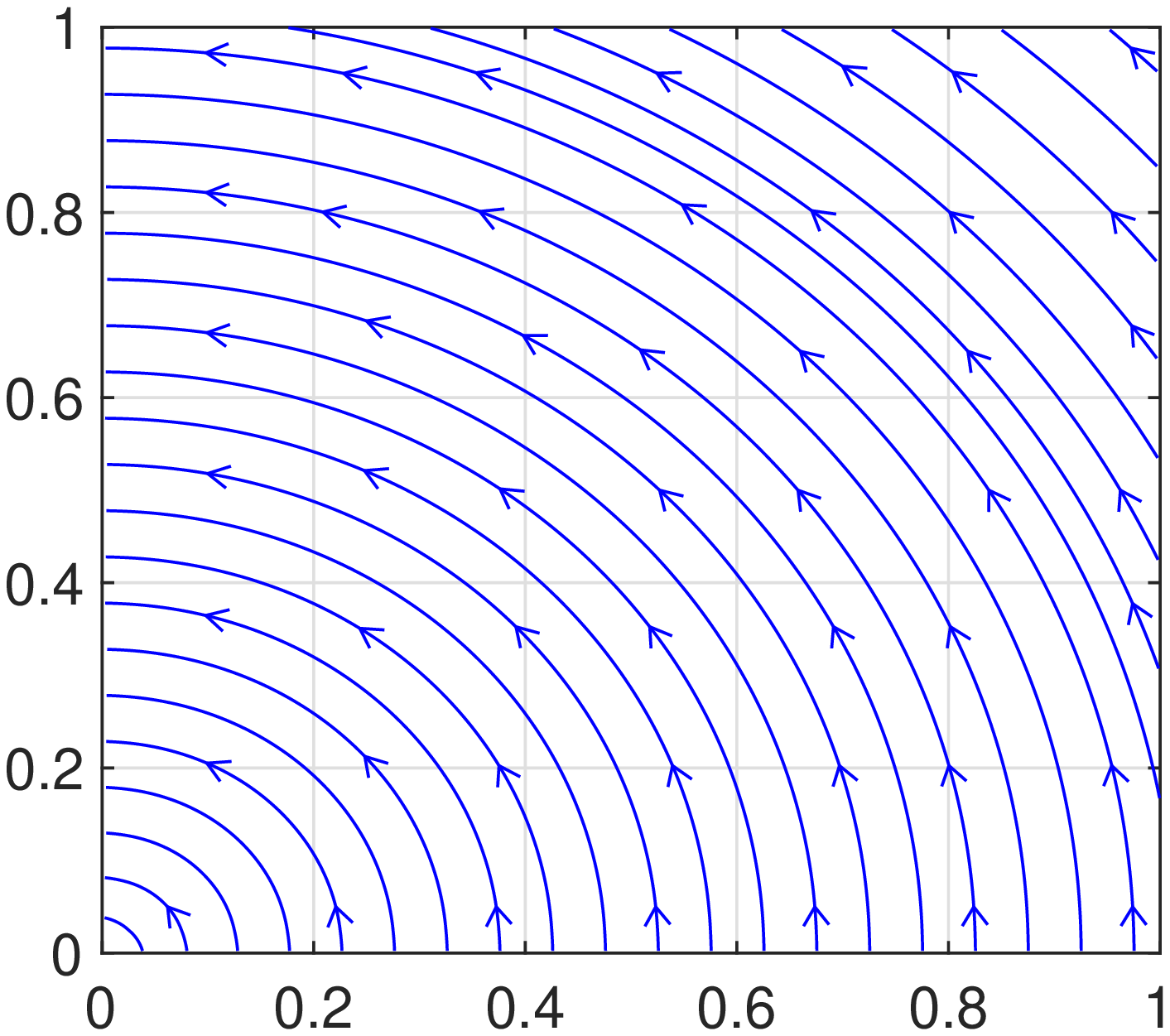}
&
\includegraphics[width=.45\textwidth]{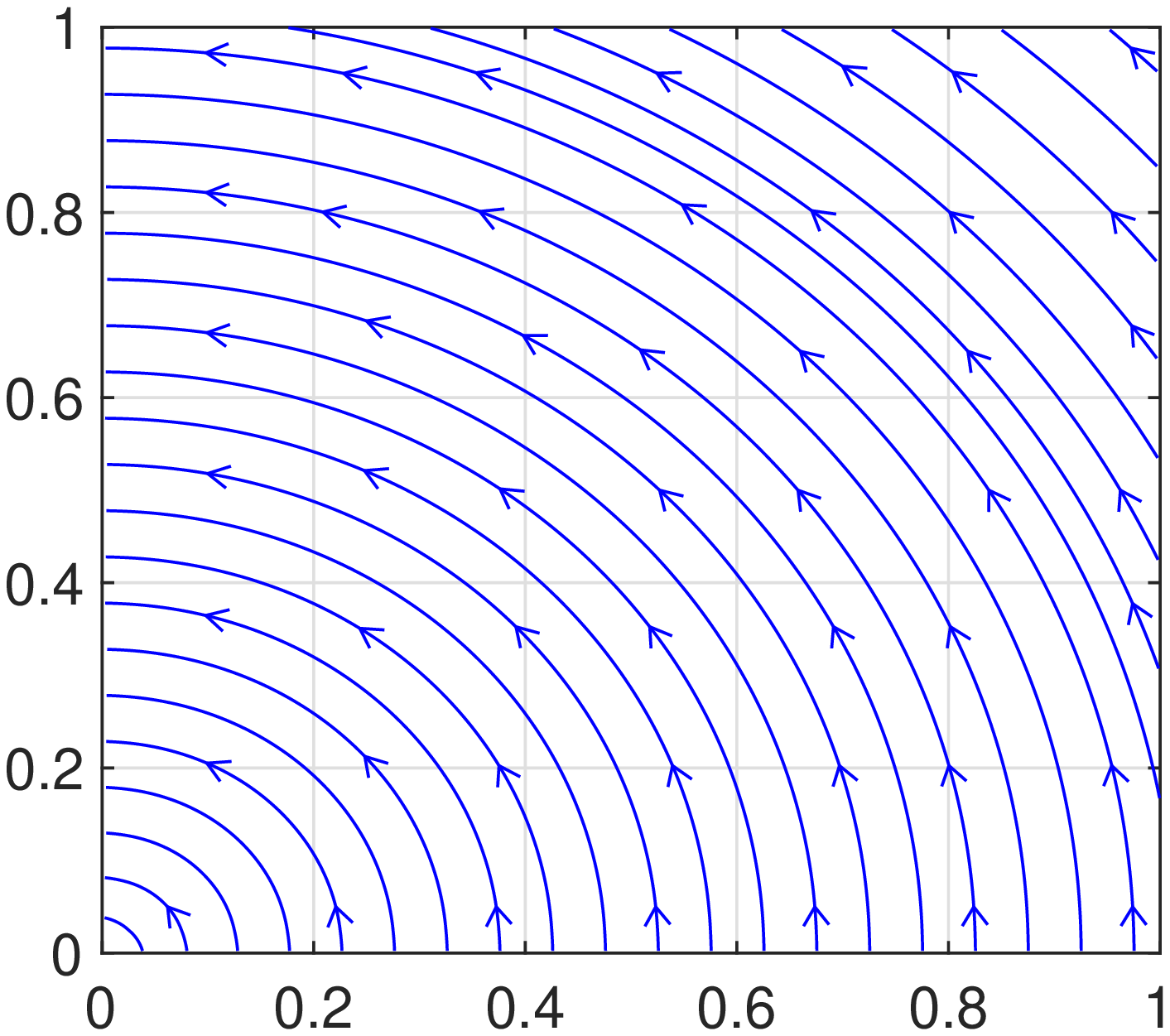}
\\
(a) & (b)
\end{tabular}
\caption{Example~\ref{Sect:Num-5}: Streamline Plots of WG Algorithm~\ref{alg:WG}  on mesh $h=1/40$ and $k = 0$  with (a) $\lambda = 10$; (b) $\lambda = 1\text{E+}6$.  }\label{Fig:Test5-2}
\end{figure}

The streamline plots are shown in Fig.~\ref{Fig:Test5-1}-\ref{Fig:Test5-2}. As one can observe that by increasing the value in $\lambda$ will change the numerical velocity noticeably for Algorithm~\ref{alg:WG-2}. Actually, as $\lambda = $1E6, the simulation by Algorithm~\ref{alg:WG-2} will produce a wrong solution on $h = 1/40$ and $k=0$. The numerical performance for Algorithm~\ref{alg:WG-2} is illustrated in Fig.~\ref{Fig:Test5-2}, from which we can notice the robustness with respect to the irrotational body forces. By comparing the two sub-figures in Fig.~\ref{Fig:Test5-2}, one can not find the difference in the streamline plot of velocity $\bu_h$. 

Next, the error profiles and convergence results are reported in Table~\ref{Tab:Num-5_WG2} and Table~\ref{Tab:Num-5_WG1}. One can notice the significant error increasing from $\lambda = 10$ to $\lambda = $1E6 in Table~\ref{Tab:Num-5_WG2}. However, we can notice that, by employing Algorithm~\ref{alg:WG}, the error for velocity is almost 0. It shows that even by the constant WG element, Algorithm~\ref{alg:WG} can produce nearly exact velocity and pressure simulation.

\begin{table}
\caption{Example~\ref{Sect:Num-4}: Error Profiles and Convergence Results for WG Algorithm~\ref{alg:WG-2}.}\label{Tab:Num-5_WG2}
\centering
\tabcolsep=3pt
\begin{tabular}{c|cc|cc|cc||cc|cc|cc}
\hline\hline
 &\multicolumn{6}{c||}{$\lambda = 10$} & \multicolumn{6}{c}{$\lambda = 1\text{E+}6$} \\
$1/h$	&$\3bar\be_h\3bar$	&Rate	&$\|\be_0\|$ &Rate		&$\|\epsilon_h\|$ &Rate	&$\3bar\be_h\3bar$	&Rate	&$\|\be_0\|$ &Rate		&$\|\epsilon_h\|$ &Rate	\\ \hline\hline
\multicolumn{13}{c}{$k = 0$}\\ \hline
16	&2.82E-1&	 	&6.89E-3&	 	&2.09E-1& &&&&&&\\						
32	&1.43E-1	&0.98	&1.78E-3	&1.95	&1.41E-1	&0.57		&1.24E+4	& 	&1.32E+2&	 	&1.77E+5&\\	
64	&7.17E-2	&0.99	&4.50E-4	&1.98	&8.40E-2	&0.74		&6.94E+3	&0.84	&4.11E+1	&1.69	&2.17E+4	&3.03\\
128	&3.59E-2	&1.00	&1.13E-4	&1.99	&4.66E-2	&0.85		&3.49E+3	&0.99	&1.07E+1	&1.94	&5.30E+3	&2.03\\ \hline
\multicolumn{13}{c}{$k = 1$}\\ \hline												
16	&5.37E-3&	 	&4.78E-5&	 	&9.14E-3&			&5.21E+2&	 	&4.64&	 	&8.38E+2&\\	
32	&1.36E-3	&1.99	&6.06E-6	&2.98	&2.30E-3	&1.99		&1.31E+2	&1.99	&5.85E-1	&2.99	&2.21E+2	&1.92\\
64	&3.40E-4	&1.99	&7.64E-7	&2.99	&5.76E-4	&2.00		&3.30E+1	&1.99	&7.37E-2	&2.99	&5.57E+1	&1.99\\
128	&8.52E-5	&2.00	&9.58E-8	&2.99	&1.44E-4	&2.00		&8.27	&2.00	&9.25E-3	&2.99	&1.39E+1	&2.00\\ \hline
\multicolumn{13}{c}{$k = 2$}\\ \hline													
16	&4.49E-5&	 	&2.40E-7&	 	&4.25E-5&			&4.49&	 	&2.40E-2&	 	&4.25&\\	
32	&5.65E-6	&2.99	&1.51E-8	&3.99	&5.32E-6	&3.00		&5.65E-1	&2.99	&1.51E-3	&3.99	&5.32E-1	&3.00\\
64	&7.08E-7	&3.00	&9.48E-10&3.99	&6.65E-7	&3.00		&7.08E-2	&3.00	&9.48E-5	&3.99	&6.65E-2	&3.00\\
128	&8.86E-8	&3.00	&5.94E-11	 &4.00	&8.39E-8	&2.99		&8.85E-3	&3.00	&5.93E-6	&4.00	&8.31E-3	&3.00
\\ \hline\hline
\end{tabular}
\end{table}

\begin{table}
\caption{Example~\ref{Sect:Num-4}: Error Profiles and Convergence Results for WG Algorithm~\ref{alg:WG}.}\label{Tab:Num-5_WG1}
\centering
\tabcolsep=3pt
\begin{tabular}{c|cc|cc|cc||cc|cc|cc}
\hline\hline
 &\multicolumn{6}{c||}{$\lambda = 10$} & \multicolumn{6}{c}{$\lambda = 1\text{E+}6$} \\
$1/h$	&$\3bar\be_h\3bar$	&Rate	&$\|\be_0\|$ &Rate		&$\|\epsilon_h\|$ &Rate	&$\3bar\be_h\3bar$	&Rate	&$\|\be_0\|$ &Rate		&$\|\epsilon_h\|$ &Rate	\\ \hline\hline
\multicolumn{13}{c}{$k = 0$}\\ \hline
16	&1.91e-13&- 	 	&4.37E-15&-	 	&1.72E-13&-		&1.91e-11&- 	 	&7.71E-13&-	 	&1.43E-09&-\\
32	&3.81e-13&- 		&1.35E-14&-	 	&5.37E-13&-		&1.61e-11&- 	 	&4.73E-13&-	 	&8.75E-10&-\\
64	&8.75e-13&- 		&5.55E-14&-	 	&2.18E-12&-		&1.71e-11&- 	 	&4.89E-13&-	 	&6.11E-10&-\\
128	&2.31e-12& -		&2.27E-13&-	 	&8.73E-12&-		&2.01e-11&- 	 	&6.49E-13&-	 	&1.63E-09&-\\ \hline
	\multicolumn{13}{c}{$k = 0$}\\ \hline									
16	&3.04E-13&-		&4.51E-15&-	 	&5.50E-13&-		&2.62E-11&-	 	&1.14E-12&-	 	&1.51E-09&-\\
32	&6.11E-13&-		&2.10E-14&-	 	&2.93E-12&-		&2.96E-11	&- 		&1.95E-12&-	 	&2.20E-09&-\\
64	&1.36E-12&-		&8.43E-14&-	 	&5.47E-12&-		&3.06E-11	&- 		&2.12E-12&-	 	&2.90E-09&-\\
128	&3.67E-12&-		&3.68E-13&-	 	&4.52E-11	&-		&5.54E-11&-	 	&6.02E-12&-	 	&6.02E-09&-\\ \hline
	\multicolumn{13}{c}{$k = 0$}\\ \hline								
16	&6.24E-13&-		&1.49E-14&-	 	&6.92E-12&-		&9.66E-11&-	 	&2.06E-12&-	 	&4.41E-09&-\\
32	&1.32E-12&-		&7.01E-14&-	 	&8.37E-12&-		&9.59E-11&-	 	&1.39E-12&-	 	&1.30E-08&-\\
64	&3.36E-12&-		&2.82E-13&-	 	&1.43E-11	&-		&9.68E-11	&- 		&1.04E-12&-	 	&1.80E-09&-\\
128	&1.01E-11	&-	        &1.15E-12	&- 	        &3.72E-10	&-		&9.77E-11	&- 		&1.45E-12&-	 	&2.30E-09&-
\\ \hline\hline
\end{tabular}
\end{table}

%\subsection{Test}
%Let $\omega=(0,1)^2$ and the load term $\bbf$ be chosen such that the analytical solution is
%\begin{eqnarray*}
%\bu(x,y) = \frac{1}{2}\begin{pmatrix}
%\sin^2(2\pi x)\sin(2\pi y)\cos(2\pi y)\\
%-\sin^2(2\pi y)\sin(2\pi x)\cos(2\pi x)
%\end{pmatrix}, \text{ and }p(x,y) = \pi^2\sin(2\pi x)\cos(2\pi y).
%\end{eqnarray*}
%
%%\subsection{Test}
%%Let $\omega=(0,1)^2$ and the load term $\bbf$ be chosen such that the analytical solution is
%%\begin{eqnarray*}
%%\bu(x,y) = \begin{pmatrix}
%%0.1 x^2(1-x)^2(4y^3-6y^2+2y)\\
%%0.1 y^2(1-y)^2(4x^3-6x^2+2x)
%%\end{pmatrix}, \text{ and }p(x,y) = x^3y^3-1/16.
%%\end{eqnarray*}
%%Note that in this test the velocity $\bu$ is not divergence free, i.e., $\nabla\cdot\bu = g\neq 0$. 
%
%
%
%\subsection{Test}
%\begin{eqnarray*}
%\bu = \begin{pmatrix}
%\cos(x)/\sqrt{x^2+y^2}\\
%\sin(x)/\sqrt{x^2+y^2}
%\end{pmatrix}
%\end{eqnarray*}

\subsection{Two-dimensional Lid-driven Cavity Flow}\label{Sect:Num-6}
\begin{figure}[H]
\centering
\begin{tabular}{ccc}
\includegraphics[width=.32\textwidth]{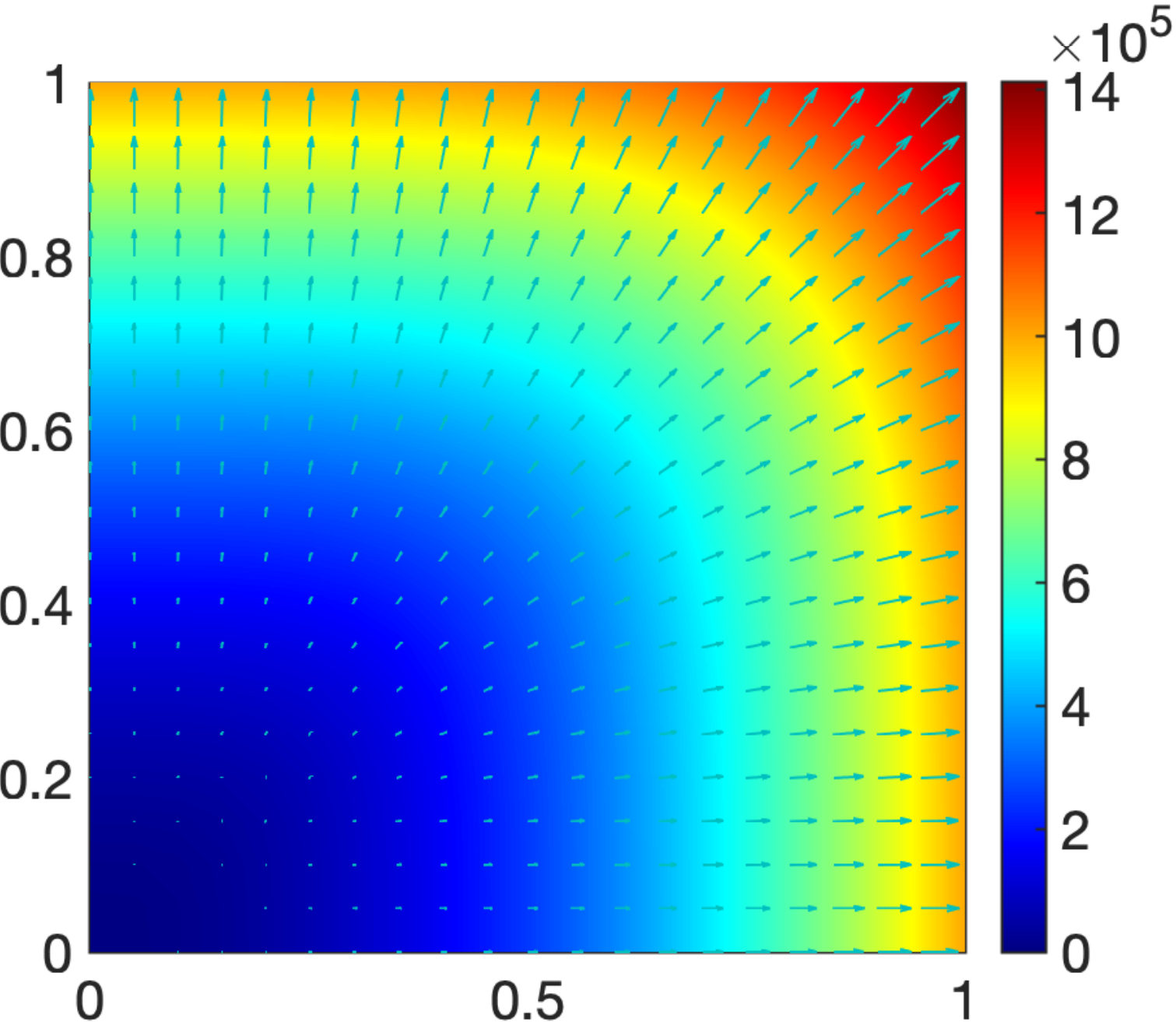}
&
\includegraphics[width=.32\textwidth]{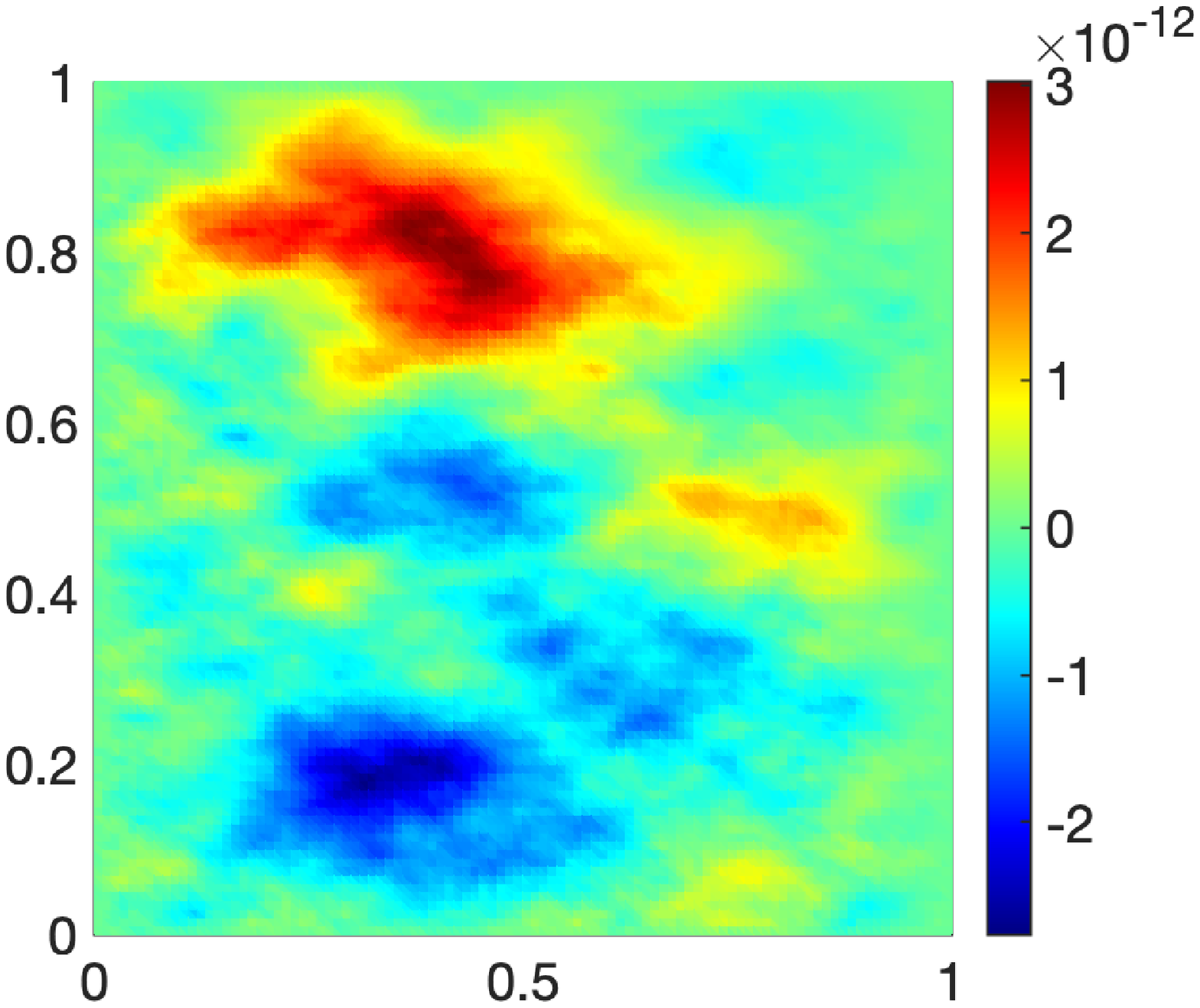}%Lambda1e6}
&
\includegraphics[width=.32\textwidth]{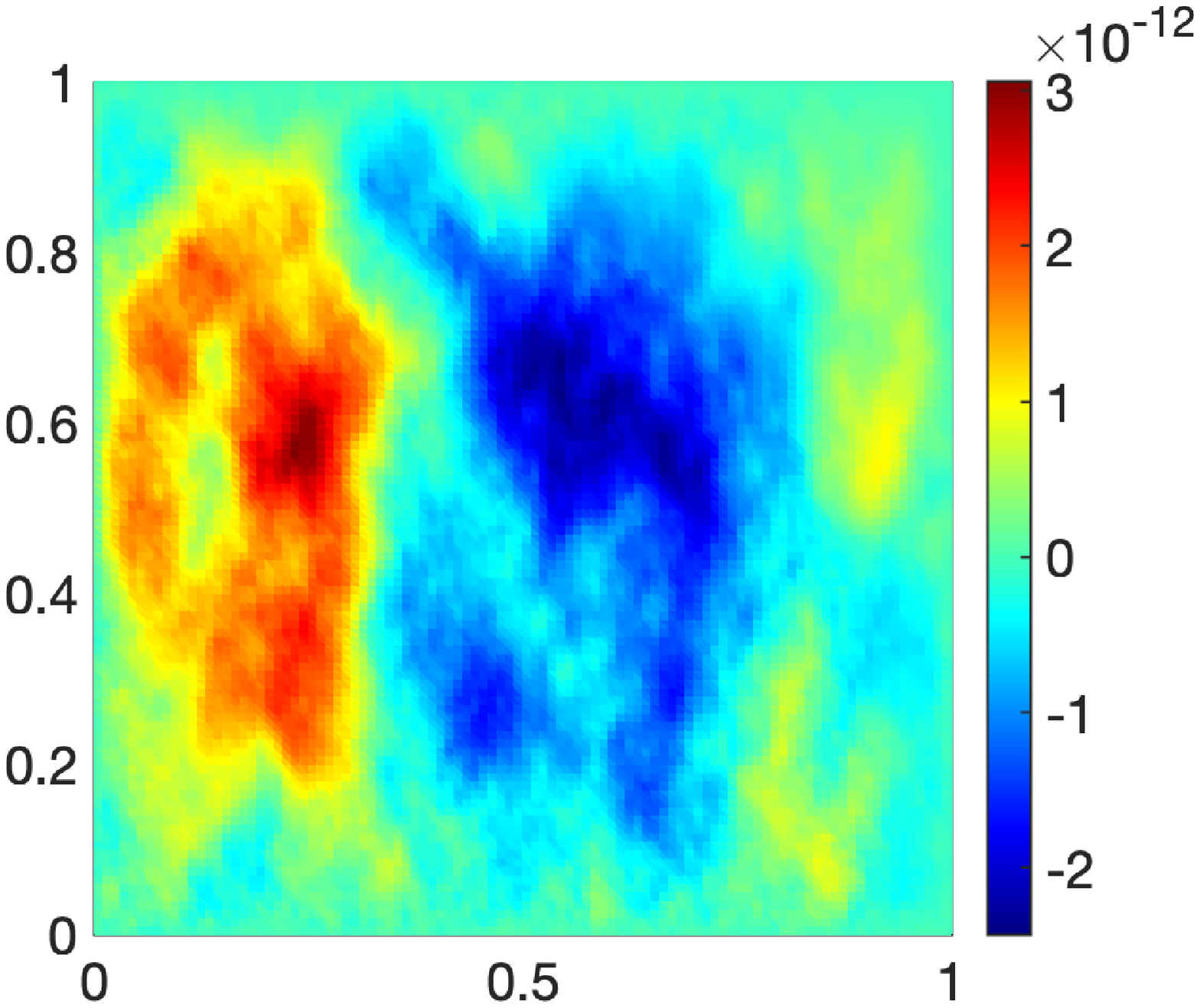}%Lid_Cavity}
\\
(a) & (b) &(c)
\end{tabular}
\caption{Example~\ref{Sect:Num-6}: (a). Plot of body force $\bbf_2$; (b) Difference of the first component of velocity for $\bbf_1$ and $\bbf_2$; (c) Difference of the second component of velocity for $\bbf_1$ and $\bbf_2$.  }\label{Fig:Test-LidCavity}
\end{figure}
In this example, we shall again validate the independence of the irrotational body force. 
Let $\Omega=(0,1)^2$ and $\nu = 1$. A lid-driven cavity flow is considered in this test. The Dirichlet boundary condition is given as
\begin{eqnarray*}
\bu|_{\partial\Omega} =\begin{cases}
(1,0)^\top,\text{ if }y=1,\\
(0,0)^\top,\text{ else}.
\end{cases}
\end{eqnarray*}
In this test, we let $\bbf_1 = 0$ and $\bbf_2 = \lambda\nabla\dfrac{1}{3}(x^3+y^3)$ (as shown in Fig.~\ref{Fig:Test-LidCavity}a) to perform the WG Algorithm~\ref{alg:WG}. It is easy to check that $\nabla\times\bbf_2 = 0$.

\begin{figure}[H]
\centering
\begin{tabular}{cc}
\includegraphics[width=.45\textwidth]{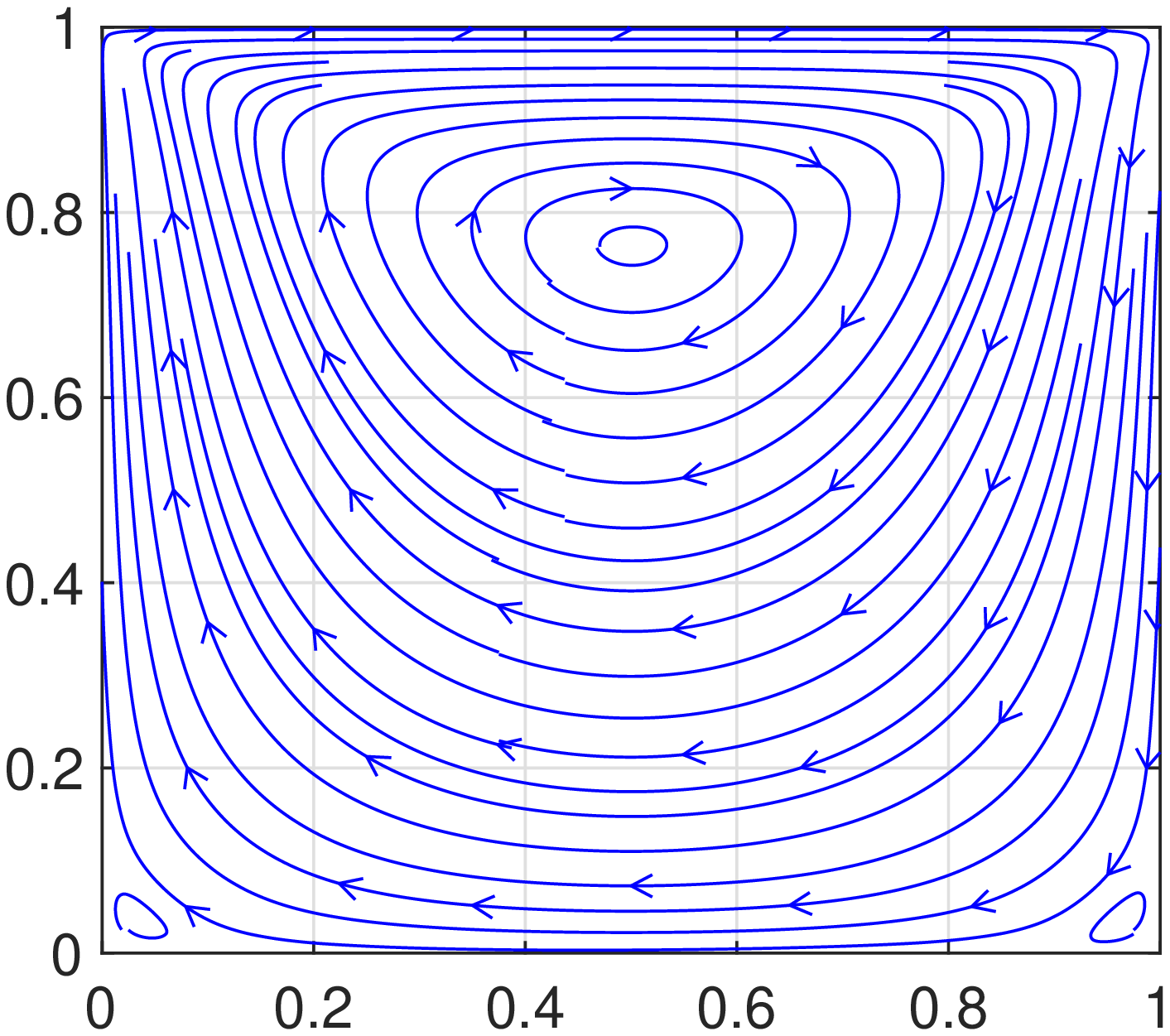}
&
\includegraphics[width=.45\textwidth]{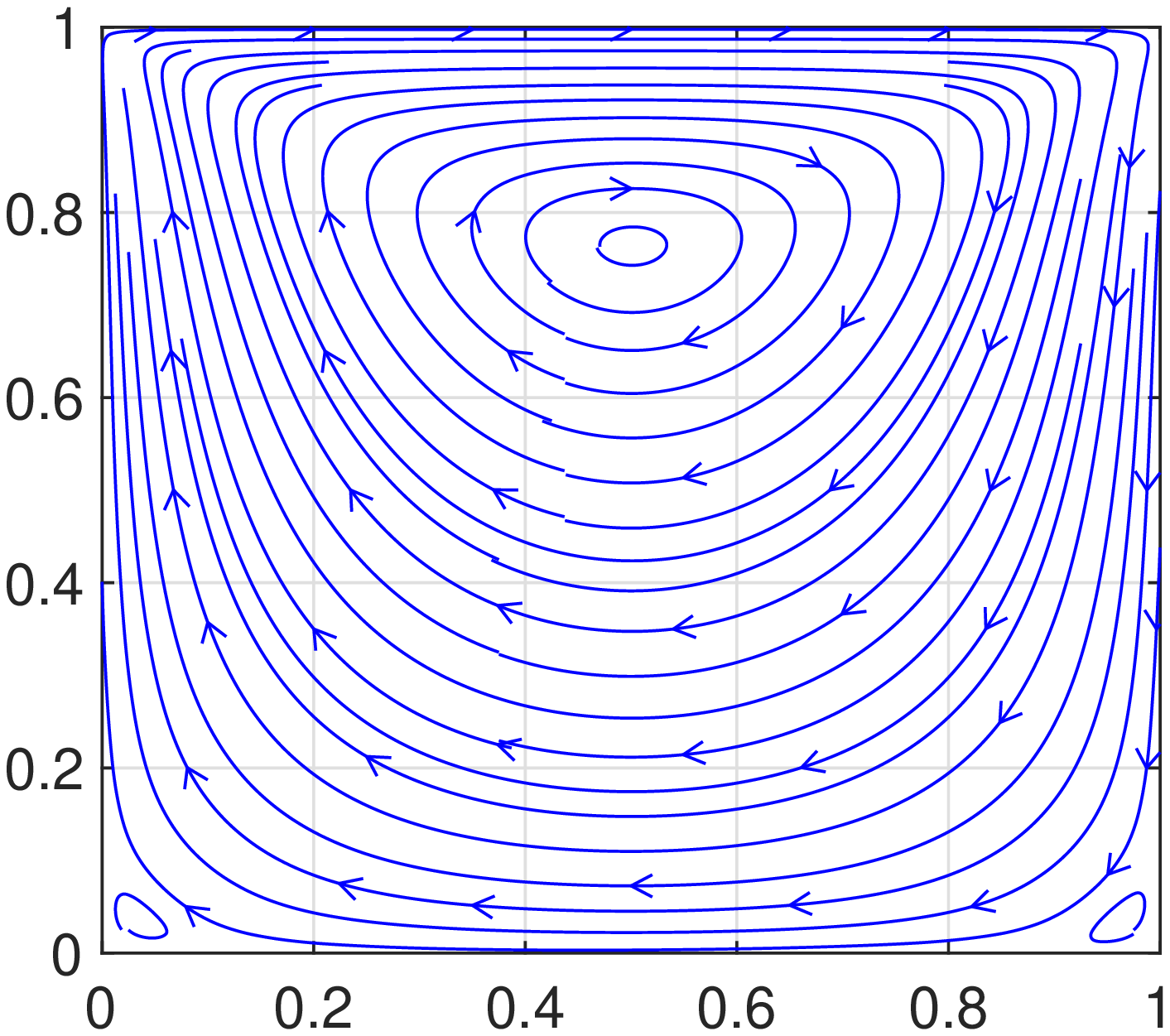}%Lambda1e6}
\\
(a) & (b) 
\end{tabular}
\caption{Example~\ref{Sect:Num-6}: Plots of streamline on mesh $h=0.01$ with  WG element $k=0$ for (a) $\bbf_1$; (b) $\bbf_2$.  }\label{Fig:Test-LidCavity2}
\end{figure}

By employing Algorithm~\ref{alg:WG} on the mesh with $h=1/100$ and $k=0,$
the difference between numerical solutions in velocity are plotted in Fig.~\ref{Fig:Test-LidCavity}b and Fig.~\ref{Fig:Test-LidCavity}c for $\bbf = \bbf_1$ and $\bbf = \bbf_2$. As one can observe from the plot, the difference is nearly zero, and thus validate our theoretical conclusions regarding the robustness with respect to irrotational body force. Then the streamlines corresponding to $\bbf_1$ and $\bbf_2$ are plotted in Fig.~\ref{Fig:Test-LidCavity2}, which again validate the invariance of irrotational body force.

\subsection{Incompressible Flow with Obstacle}\label{Sect:Num-7}
\begin{figure}[H]
\centering
\begin{tabular}{cc}
\includegraphics[width=.45\textwidth]{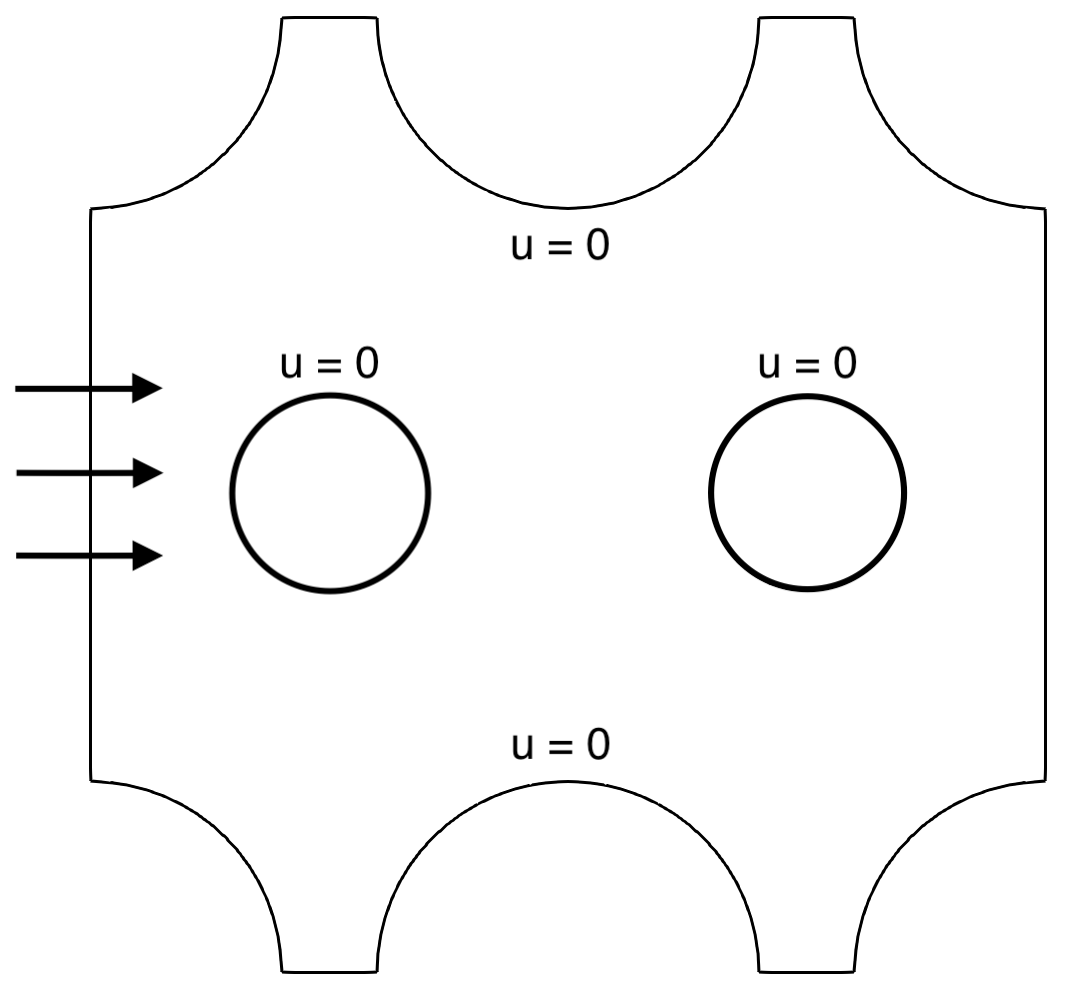}
&
\includegraphics[width=.45\textwidth]{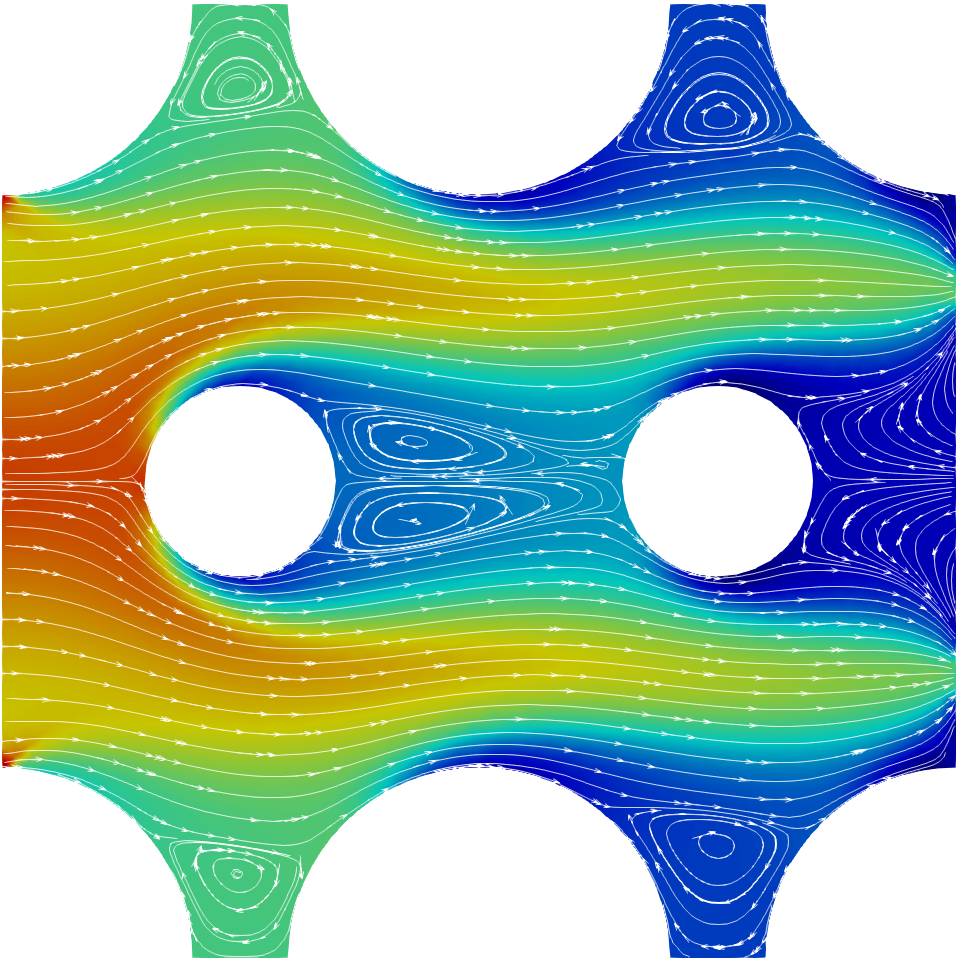}
\\
(a) & (b)
\end{tabular}
\caption{Example~\ref{Sect:Num-7}: (a). computational domain; (b). plot of streamline and pressure.  }\label{Fig:Num-7}
\end{figure}

In this test, we shall illustrate the numerical performance of Algorithm~\ref{alg:WG} for the incompressible flow with obstacle. The computational domain is plotted in Fig.~\ref{Fig:Num-7} (a). We assume the inflow boundary condition $\bu=(1,0)^\top$ on the left edge and outlet boundary condition on the left. The other boundary is assume to be wall boundary condition. 

Let $\nu = 1$, and we perform Algorithm~\ref{alg:WG} with $k = 0$. The numerical solution is plotted in Fig.~\ref{Fig:Num-7} (b). As the streamline plot for velocity, one can clear detect the vortexes in the simulation.

%\subsection{Test}
%\begin{figure}[H]
%\centering
%\begin{tabular}{cc}
%\includegraphics[width=.45\textwidth,height=.3\textwidth]{./figure/Domain2_mesh4_k0_nu1e-1v1e-1}
%&
%\includegraphics[width=.45\textwidth,height=.3\textwidth]{./figure/Domain2_mesh4_k1_nu1e-2v1e-1}
%\\
%(a) & (b)
%\end{tabular}
%\caption{Example~\ref{Num-Test5}: Plots of WG approximation properties (a) $\nu = 1\text{E-}1$ and $u_1 = 1\text{E}-1$; (b) $\nu = 1\text{E-}2$ and $u_1 = 1\text{E}-1$.  }%\label{Fig:Test-LidCavity}
%\end{figure}

%===============
% Conclusion
%===============
\section{Conclusion}\label{Sect:Con}
In this paper, we developed a pressure-robust weak Galerkin finite element methods for Navier-Stokes equations. By employing the velocity reconstruction operator in the body force assembling and the convective term, our algorithm can achieve the independence of pressure in the error estimate. Numerical tests in two dimensions have been carried out to validate the robustness of pressure and the irrotational body force, and thus confirm the advantages of the proposed approach. Three-dimensional numerical investigation will be carried out in the future. Besides,
the extension to the non-stationary Navier-Stokes equations and numerical scheme with upwind stabilization will be also investigated in the future. 

%%%%%%
% Ref
%%%%%%

\end{document}